\documentclass[final]{siamltex}
\usepackage[latin1]{inputenc}
\usepackage{amsmath}
\usepackage{amssymb}
\usepackage{mathrsfs}
\usepackage{epsfig,psfrag}
\usepackage{mathabx}
\usepackage{array}
\usepackage{showkeys}

\def\R{\mathbb{R}}
\def\T{\mathbb{T}}

\def\b0{\boldsymbol{0}}
\def\bA{\boldsymbol{A}}
\def\bB{\boldsymbol{B}}

\def\bbf{\boldsymbol{f}}
\def\bG{\boldsymbol{G}}
\def\bS{\boldsymbol{S}}
\def\tm{\widetilde{m}}
\def\tM{\widetilde{M}}
\def\bu{\boldsymbol{u}}
\def\bv{\boldsymbol{v}}
\def\bx{\boldsymbol{x}}
\def\bdelta{\boldsymbol{\delta}}
\def\bsigma{\boldsymbol{\sigma}}

\def\txi{\widetilde{\xi}}
\def\Re{\mathfrak{Re}}
\def\We{\mathfrak{We}}
\def\div{\mbox{\,{\textrm{div}}}}
\def\transp{\,^T\!}
\newcommand{\cqfd}{{\nobreak\hfil\penalty50\hskip2em\hbox{}\nobreak\hfil $\square$\quad\parfillskip=0pt\finalhyphendemerits=0\par\medskip}}

\newtheorem{remark}{Remark}[section]

\title{Mathematical existence results \\for the Doi-Edwards polymer model}

\author{Laurent Chupin\thanks{Universit\'e Blaise Pascal, Clermont-Ferrand II, Laboratoire de Math\'ematiques CNRS-UMR 6620, Campus des C\'ezeaux, F-63177 Aubi\`ere cedex, France ({\tt laurent.chupin@math.univ-bpclermont.fr}).}}

\begin{document}

\maketitle

\begin{abstract}
In this paper, we present some mathematical results on the Doi-Edwards model describing the dynamics of flexible polymers in melts and concentrated solutions.
This model, developed in the late 1970s, has been used and tested extensively in modeling and simulation of polymer flows.
From a mathematical point of view, the Doi-Edwards model consists in a strong coupling between the Navier-Stokes equations and a highly nonlinear constitutive law.
\par\noindent
The aim of this article is to provide a rigorous proof of the well-posedness of the Doi-Edwards model, namely it has a unique regular solution. We also prove, which is generally much more difficult for flows of viscoelastic type, that the solution is global in time in the two dimensional case, without any restriction on the smallness of the data.
\end{abstract}

\begin{keywords} 
Polymer, Viscoelastic flow, Doi-Edwards model, Navier-Stokes equations, global existence result.
\end{keywords}

\begin{AMS}
35A01, 35B45,  35Q35, 76A05, 76A10, 76D05
\end{AMS}

\pagestyle{myheadings}
\thispagestyle{plain}
\markboth{LAURENT CHUPIN}{Existence for the Doi-Edwards polymer model}

\section{Introduction}\label{part:introduction}

Numerous models exist for describing fluids with complex rheological properties. They generally are of great scientific interest and have a rich phenomenology. Their mathematical description remains challenging.  We are interested in this article to a model - the Doi-Edwards model - which was one of the foundation of the most recent physical theories but for which the mathematical theory remains very poor.\\[0.2cm]
M. Doi and S.F. Edwards wrote a series of papers~\cite{DoiEdwards78-1, DoiEdwards78-2, DoiEdwards78-3, DoiEdwards78-4} expanding the concept of reptation introduced by P.G. de Gennes in 1971. This approach was then taken up in a famous book in polymer physics in 1988, see~\cite{Doi}.
Since this model was derived, numerous studies have been carried out either from a physical point of view, either from a numerical point of view.
Moreover, several other models were born: simplified models using for instance the Independent Alignment approximation or the Currie approximation~\cite{Currie}, more complex models like the pom-pom model~\cite{Leish98}. Finally, much progress has been made on the modeling of both linear and branched polymers.
However, from a mathematical point of view, it seems that no justification was given even for the pioneering model.\\[0.2cm]
Nonetheless we can cite some recent theoretical papers on this subject, see~\cite{Ciuperca12, Heibig14}, in which the authors are only interested in specific cases: one dimensional shear flows under the independent alignment assumption in~\cite{Heibig14}, flows for which the coupling between the velocity and the stress is not taking into account, see~\cite{Ciuperca12}.
More generally, there seems to be a real challenge to obtain global existence in time for models of polymers. The most caricatural example is the Oldroyd model for which the question of global existence in dimension~$2$ remains an open question, see~\cite{Lions-Masmoudi-viscoelastique} for a partial answer.
However, there exists polymer models for which such results are proved. Thus, for the FENE type models, N. Masmoudi~\cite{Masmoudi13} proved a global existence result in dimension~$2$. Similarly, for integral fluid of type K-BKZ such results hold too (see~\cite{Chupin4}).\\[0.2cm]
The aim of this paper is to prove relevant mathematical results on this relevant physical problem. The first one is the following (a more specific version of this result is given by Theorems~\ref{th:local} and~\ref{th:unic}, page~\pageref{th:local}):\\[-0.2cm]
\begin{theorem}
There exists a time~$t^\star>0$ such that the Doi-Edwards model admits a unique strong solution on the interval time $[0,t^\star]$.
\end{theorem}
\par\noindent\\[-0.2cm]
By the expression "strong solution" we mean a sufficiently smooth solution so that each term of the system is well defined, as well as the initial conditions (corresponding to the time $t=0$).
The lifetime of the solution, {\it i.e.} the value of time~$t^\star$, is not easily quantifiable. In practice it is well known that for the Navier-Stokes equations - modeling newtonian behavior, the question of long time existence is still an open one. For the Newtonian fluids, the only existence results, for long time and for any data, correspond to the two dimensional case.
However, we have above pointed out the difficulty to get this kind of result even in~$2D$ for some viscoelastic fluids.
The major point of this paper is the proof indicating that the model of Doi-Edwards admits a strong solution for long time in~$2D$ (a more specific version of this result is given by Theorem~\ref{th:global}, page~\pageref{th:global}):\\[-0.2cm]
\begin{theorem}
For all time~$t^\star>0$, the Doi-Edwards model admits a unique strong solution on the interval time $[0,t^\star]$.
\end{theorem}
\par\noindent\\[-0.2cm]
The paper is organized as follows.
First - in Section~\ref{part:model}, we introduce the Doi-Edwards model specifying the physical meaning of each contribution. This second section ends by a dimensionless procedure that allows us to write the model with only three parameters (the Reynolds number, the Weissenberg number and the ratio between solvent viscosity and elastic viscosity).
In Section~\ref{part:framework}, we present the mathematical framework as well as the assumptions which are physically discussed. The main results are given at the end of this section.
Section~\ref{part:preliminaries} is devoted to fundamental preliminaries which correspond to some key points of the next proofs. The first two provide {\it a priori} bounds which will imply that that the stress defined in the Doi-Edwards model is automatically bounded. The third preliminary give a Gronwall lemma with two time variables.
The fourth preliminary result is about the maximum principle which can be applied many times to estimate the memory of the fluid.
The last preliminary result is about a Cauchy problem arising in the global existence proof.
The three last sections (\ref{part:proof1}, \ref{part:proof2} and~\ref{part:proof3}) are devoted to the proof of the three mains results , namely: the local existence result in Section~\ref{part:proof1}, the uniqueness result in Section~\ref{part:proof2} and the global existence result in Section~\ref{part:proof3}.
Some open questions are presented by way of conclusion.

\section{Governing equations}\label{part:model}

\subsection{Conservation laws}

In this paper we are interested in the flow of isothermal and incompressible fluids.
The incompressibility implies that the mass conservation is equivalent to the free-divergence of the velocity field.
The isothermal assumption implies that only one other conservation law suffices to describe the flow: it corresponds to the law of conservation of the momentum (Newton's second law of motion). This equation is written as a balance between the material derivative of the velocity and the divergence of the Cauchy stress tensor.
For a polymeric liquid, the equations of conservation can hence written as a system coupling the velocity field~$\bv$, the pressure~$p$ and the extra-stress tensor~$\bsigma$:
\begin{equation*}\label{1}
\left\{
\begin{aligned}
& \rho \, \mathrm d_t \bv + \nabla p - \eta_{\mathrm s} \Delta \bv = \div \, \bsigma, \\
& \div \, \bv = 0,
\end{aligned}
\right.
\end{equation*}
where~$\rho$ is the fluid density and~$\eta_{\mathrm s}$ the solvent viscosity. The notation~$d_t$ corresponds to the material derivative $d_t = \partial_t + \bv\cdot \nabla$.

\subsection{Constitutive equation}

A fundamental result of the Doi-Edwards theory is the expression for the stress tensor~$\bsigma$ which applies when the chains which compose the fluid are relaxed within their tubes.
More precisely, the stress can be deduce from a tensor~$\bS$ denoted the orientation order parameter of the chains.
Although the chain tension is permanently at the equilibrium value, the orientations become anisotropically distributed as a consequence of the flow, and a stress develops accordingly.
The stress is then modelized by (see~\cite[page 2056]{Doi80}):
\begin{equation*}\label{2}
\bsigma(t,\bx) = \frac{G_{\mathrm e}}{\ell} \int_{-\frac{\ell}{2}}^{\frac{\ell}{2}} \bS(t,\bx,s)\, \mathrm ds,
\end{equation*}
where $G_{\mathrm e}$ is a characteristic modulus and~$\ell$ is the equilibrium value of the contour length of the chains.
The quantity~$s$ is a arc-length coordinate along the primitive chain. 
\par
\begin{figure}[htbp]
\begin{center}
\includegraphics[height=4cm]{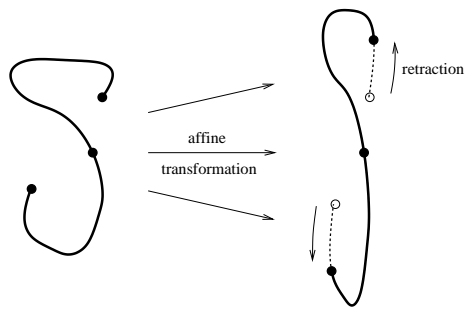}
\caption{(inspired from the book~\cite[page 2058]{Doi80}) --
When a macroscopic deformation is applied, a polymer chain is transformed into a new conformation. The new chain is on the curve which is the affine transformation of the initial curve. The new position of each segment is obtained by retraction, preserving the initial lengths.}\label{fig1}
\end{center}
\end{figure}

\noindent
To evaluate this tensor~$\bS$, M. Doi and S.F. Edwards write $\bS = \langle \bu \otimes \bu - \frac{1}{d}\bdelta \rangle$. Here, $\bu$ is a unit vector along the tangent to the primitive chain which depends on time~$t$, spatial position~$\bx$ and length~$s$, and~$\bdelta$ denotes the identity tensor. The entire~$d$ corresponds to the dimension of the spatial coordinates (in practice $d=2$ or $d=3$). The average is over the distribution of these vectors in the ensemble of chains, {\it i.e.}, more explicitly
\begin{equation*}\label{3}
\bS(t,\bx,s) = \int_{\mathbb S^{d-1}} f(\bu;t,\bx,s) \big( \bu \otimes \bu - \frac{1}{d}\bdelta \big) \, \mathrm d\bu,
\end{equation*}
with $f(\bu;t,\bx,s)$ given the orientation distribution function.
To obtain an expression for this distribution, the history of motion must be found. To this purpose, let us first recall the relevant aspects of the Doi-Edwards model. They are
\begin{itemize}
\item[$\checkmark$] The polymer moves randomly inside the tube executing one-dimensional Brownian motion. Moreover tube segments are randomly oriented when they are created and deform affinely thereafter;
\item[$\checkmark$] If the system is macroscopically deformed, the polymer conformation is also changed as presented on Figure~\ref{fig1};
\item[$\checkmark$] The macroscopic motion and Brownian motion coexist, independently from one another.
\end{itemize}
All these considerations being taken into account, it is possible to write (see~\cite[page 277]{Doi})
\begin{equation}\label{4}
\bS(t,\bx,s) = -\int_0^{+\infty} \partial_{T} K(t,T,\bx,s)\, \mathscr S(\bG(t,T,\bx))\, \mathrm dT.
\end{equation}
It makes appear the deformation gradient tensor~$\bG$ which depends not only on the current time~$t$ and spatial variable~$\bx$ but also on an other time~$T$. The time~$T$ allows to take into account all the history of the motion.
The deformation gradient~$\bG$ fulfills the differential equation (see~\cite{Chupin3}):
\begin{equation*}\label{4.5}
\mathrm d_t \bG + \partial_T \bG = \bG\cdot \nabla \bv.
\end{equation*}
Finally, the integral kernel~$K$ satisfies (the coefficient~$D_{\mathrm e}$ is a curvilinear diffusion coefficient):
\begin{equation}\label{5}
\mathrm d_t K + \partial_T K + \Big( \nabla \bv : \int_0^s \bS\, \mathrm ds \Big) \partial_s K - D_{\mathrm e} \, \partial^2_s K = 0.
\end{equation}
\begin{remark}
\begin{enumerate}
\item Note that the equation~\eqref{4} does not explicitly defined the orientation tensor since~$\bS$ again appears in the equation~\eqref{5}. It is one of the difficulties to obtain existence results.
\item The usual formulation use the time~$t$ and a other time~$t'$ in the past. Morally, the deformation gradient tensor~$\bG$ measures the deformation between these two times. In the present paper we select as independent variable the age $T=t-t'$, which is measured relative to the current time~$t$. This viewpoint is relatively classical in the numerical framework for integral models, see for instance~\cite{Hulsen01,Keunings,Tome08}.
\end{enumerate}
\end{remark}
\vspace{0.2cm}
\noindent
The model is closed with the expression of the function~$\mathscr S$, see~\cite[eq.~(7.115)]{Doi}:
\begin{equation}\label{6}
\mathscr S(\bG) = \frac{1}{\langle |\bG \cdot \bu| \rangle_0} \, \Big\langle \frac{(\bG \cdot \bu)\otimes (\bG \cdot \bu)}{|\bG \cdot \bu|} \Big\rangle_0 - \frac{1}{d} \bdelta,
\end{equation}
where the brackets $\langle\cdot\rangle_0$ correspond to the average over the isotropic distribution of unit vectors~$\bu\in \mathbb S^{d-1}$.

\subsection{Boundary and initial conditions}

The previous equations are supplemented by boundary and initial conditions.
Throughout this article we restrict to the case where the macroscopic field is assumed to be periodic. Thus the only condition that we impose on the unknowns~$\bv$, $p$, $K$ and~$\bG$ with respect to the variable~$\bx$ is to be periodic.
Clearly this ``simplification'' is purely mathematical and it will be interesting to treat a more physical case imposing, for instance, the value of the velocity at the macroscopic boundary.
\par\noindent
By definition of the integral kernel~$K$, we impose the following conditions:
\begin{equation}\label{7}
K(t,0,\bx,s) = 1
\quad \text{and} \quad
K(t,T,\bx,-\frac{\ell}{2}) = K(t,T,\bx,\frac{\ell}{2}) = 0.
\end{equation}
In the same way, the quantity~$\bG(t,T,\bx)$ which corresponds to the deformation gradient  from a past times $t-T$ to the current time~$t$ must naturally satisfies
\begin{equation}\label{8}
\bG(t,0,\bx) = \bdelta.
\end{equation}
The initial conditions correspond to the value that we impose at time $t=0$. We assume that we know the velocity at this initial time and at any point~$\bx$ of the domain. In the same way we assume that we know all the past of the flow before the initial time: we then know the value of~$\bG$ and~$K$ at $t=0$, for any age~$T$ and at any point~$\bx$.
To summarize, we assume that there exists an initial velocity~$\bv_0$, an initial deformation~$\bG_0$ and an initial function~$K_0$ such that
\begin{equation}\label{9}
\begin{aligned}
& \bv(0,\bx) = \bv_0(\bx), \\
& \bG(0,T,\bx) = \bG_0(T,\bx), \\
& K(0,T,\bx,s) = K_0(T,\bx,s).
\end{aligned}
\end{equation}

\subsection{Remark: The I.A. approximation}\label{subsec:IA}

A common approximation for the Doi-Edwards model, called the independent alignment approximation (I.A. approximation, see~\cite[section 7.7.2]{Doi}), is to neglect the transport term $\big( \nabla \bv : \int_0^s \bS\, \mathrm ds \big) \partial_s K$ in the equation~\eqref{5}, and also to simplify the expression~\eqref{6} of the function~$\mathscr S$ using 
\begin{equation*}\label{10}
\mathscr S^{(IA)}(\bG) = \Big\langle \frac{(\bG \cdot \bu)\otimes (\bG \cdot \bu)}{|\bG \cdot \bu|^2} \Big\rangle_0 - \frac{1}{d} \bdelta.
\end{equation*}
The I.A. approximation is actually quite popular in the rheology literature (see e.g., \cite{Larson88,Lin04,Marrucci99}) since the corresponding configurational equation for~$K$ can be explicitly solved using the Fourier series. The stress tensor~$\bsigma$ is then more simply given by (see~\cite[Equation~7.195]{Doi}):
\begin{equation}\label{11}
\begin{aligned}
& \bsigma^{(IA)}(t,\bx) = G_{\mathrm e} \int_0^{+\infty} m(T)\, \mathscr S^{(IA)}(\bG(t,T,\bx))\, \mathrm dT,\\
& \text{where} \quad m(T) = \sum_{p\, \mathrm{odd}} \frac{8\, D_{\mathrm e}}{\pi^2 \, \ell^2}\mathrm{exp}\Big(\frac{-T\, D_{\mathrm e}}{\ell^2}p^2 \Big).
\end{aligned}
\end{equation}
For such a model the global existence result is a consequence of a general result on viscoelastic flows with memory, see~\cite{Chupin4}.
Nevertheless, it is also well known that this approximation causes serious error in certain situations, this is clearly specified in the seminal book~\cite{Doi}.
More precisely, it is proved that I.A. predicts a negative Weissenberg effect (see~\cite{Hassager85}) while the version without I.A. predicts a positive Weissenberg effect (see~\cite{Marrucci88}).
To paraphrase M. Doi~\cite[page~2064]{Doi80}: "Mathematically [...] there seems no {\it a priori} reason why the term $\big( \nabla \bv : \int_0^s \bS\, \mathrm ds \big) \partial_s$ can be neglected compared with the term $D_{\mathrm e}\partial_s^2$".

\subsection{Dimensionless procedure}

In order to recover characteristic properties of the system, we use a nondimensionalization procedure. We denote by $L$ a characteristic macroscopic length, by $V$ a characteristic velocity of the flow. It is then natural to define a dimensionless coordinates~$\bx^\star$, a dimensionless velocity~$\bv^\star$ and a dimensionless time~$t^\star$ by the following relations
\begin{equation*}\label{12}
\bx = L \bx^\star, \qquad \bv = V\bv^\star, \qquad t = \frac{L}{V}t^\star.
\end{equation*}
For polymer flow, there exists also two microscopic characteristic sizes which correspond to the length~$\ell$ and to the diffusion coefficient~$D_{\mathrm e}$. They allow to define a dimensionless microscopic length~$s^\star$ and another dimensionless time~$T^\star$:
\begin{equation*}\label{13}
s = \ell \, s^\star, \qquad
T = \frac{\ell^2}{D_{\mathrm e}}\, T^\star.
\end{equation*}
Finally, in a dilute polymer solution, two viscosities naturally appear: the solvent viscosity~$\eta_{\mathrm s}$ and the elastic one defined using the characteristic modulus:~$\eta_{\mathrm e} = \frac{LG_{\mathrm e}}{V}$. If we denote by $\eta = \eta_{\mathrm s} + \eta_{\mathrm e}$ the total viscosity, then we defined the dimensionless pressure and stress as follows
\begin{equation*}\label{14}
p = \frac{\eta V}{L} p^\star, \qquad \bsigma = \frac{\eta V}{L} \bsigma^\star.
\end{equation*}
Taking into account all these new unknowns and new variables, the complete system reads (without $\star$ in the notations):
\begin{subequations} \label{system}
 \begin{align}
& \Re \, \mathrm d_t \bv + \nabla p - (1-\omega) \Delta \bv = \div \, \bsigma,\label{eq:1} \\[0.25cm]
& \div \, \bv = 0,\label{eq:2} \\
& \bsigma(t,\bx) = \omega \int_{-\frac{1}{2}}^{\frac{1}{2}} \bS(t,\bx,s)\, \mathrm ds,\label{eq:3}\\
& \bS(t,\bx,s) = -\int_0^{+\infty} \partial_{T} K(t,T,\bx,s)\, \mathscr S(\bG(t,T,\bx))\, \mathrm dT,\label{eq:4}\\
& \mathrm d_t \bG + \frac{1}{\We} \partial_T \bG = \bG\cdot \nabla \bv,\label{eq:5}\\
& \mathrm d_t K + \frac{1}{\We} \partial_T K + \Big( \nabla \bv: \int_0^s \bS \Big) \partial_s K - \frac{1}{\We} \partial^2_s K = 0.\label{eq:6}
\end{align}
\end{subequations}
In this set of equations, $\Re$ is the usual Reynolds number, $\omega$ stands for the viscosities ratio and~$\We$ is the Weissenberg number defined by the ratio between the macroscopic time and the microscopic time. More precisely we have
\begin{equation*}\label{15}
\Re = \frac{\rho V L}{\eta},
\qquad
\omega = \frac{\eta_{\mathrm e}}{\eta},
\qquad
\We = \frac{\ell^2/D_{\mathrm e}}{L/V}.
\end{equation*}
The functions~$\mathscr S$ is always defined by the relation~\eqref{6}.\\

\noindent
The goal of the rest of the paper is to analyze, from a mathematical point of view, the existence of a solution to the system~\eqref{system}.
More exactly, by given initial data $(\bv_0, \bG_0, K_0$), is there a triplet of functions $(\bv, \bG, K)$ which coincides with the data at initial time and such that the previous system holds for any future time?

\section{Mathematical framework, assumptions and main results}\label{part:framework}

\subsection{Notations}

The integer~$d$ stands for the spatial dimension of the flow. It will be equal to~$2$ or~$3$ in the first parts and exclusively equal to~$2$ in the Section~\ref{part:proof3} where we prove a global existence result.\\[-0.2cm]
\par\noindent
{\bf Notations for functional analysis} --
\begin{itemize}
\item[$-$] The $d$ dimensional torus is denoted~$\T$.
\item[$-$] For all real $n\geq 0$ and all integer $q\geq 1$, the set $W^{n,q}_{\bx}$ corresponds to the usual Sobolev spaces with respect to the space variable~$\bx\in \T$.
We classically denote $L^q_{\bx}=W^{0,q}_{\bx}$ and $H^n_{\bx}=W^{n,2}_{\bx}$ and we do not take into account the dimension in the notations, for instance the space $(W^{1,q}_{\bx})^3$ will be denoted $W^{1,q}_{\bx}$.
\item[$-$] All the norms will be denoted by index, like $\|\bv\|_{W^{1,q}_{\bx}}$.
\item[$-$] Since we are interested in the incompressible flows, we introduce
\begin{equation*}\label{15.5}
H_q = \{\bv \in L^q_{\bx} \, ; \, \div\, \bv=0 \}.
\end{equation*}
\item[$-$] The Stokes operator $A_q$ is introduced, with domain $\mathscr D(A_q) = W^{2,q}_{\bx} \cap H_q$, whereas we denote (see~\cite[Section 2.3]{Danchin} for some details on this space)
$$\mathscr D^r_q = \{\bv \in H_q \, ; \, \|\bv\|_{L^q_{\bx}} + \big( \int_0^{+\infty} \|A_q \mathrm e^{-tA_q}\bv \|_{L^q_{\bx}}^{r} \, \mathrm dt \big)^{1/r} < +\infty \}.$$
\item[$-$] The notations of kind $L^r(0,t^\star;X)$ denote the space of $r$-integrable functions on $(0,t^\star)$ with values in the space~$X$.
For instance $\bG \in L^r(0,t^\star;L^\infty_T L^q_{\bx})$ means that
\begin{equation*}\label{15.6}
\|\bG\|_{L^r(0,t^\star;L^\infty_T L^q_{\bx})}^r := \int_0^{t^\star} \sup_{T\in \R^+} \left( \int_\T |\bG(t,\bx,T)|^q \, \mathrm d\bx \right)^{\frac{r}{q}} \mathrm dt < +\infty.
\end{equation*}
\end{itemize}
\vspace{0.2cm}
\noindent
{\bf Notations for tensorial analysis} --
In System~\eqref{system}, the first equation~\eqref{eq:1} is a vectorial equation (the velocity $\bv$ is a function with values in $\R^d$). The equations~\eqref{eq:3}, \eqref{eq:4} and \eqref{eq:5} are tensorial equations (the stress~$\bsigma$, the orientation~$\bS$ and the deformation tensor~$\bG$ are functions with values in the set of the $2$-tensors).
In the following proofs, we need to work with the gradient of such $2$-tensors, that is with $3$-tensors, and even with $4$-tensors.
We introduce here some notations for tensors.
\begin{itemize}
\item[$-$] The set of linear applications on the $d$-dimensional space is denoted~$\mathcal L(\R^d)$.
\item[$-$] The products $\bA \otimes \bB$, $\bA \cdot \bB$ and~$\bA:\bB$ between two tensors of order~$p$ and~$q$ are respectively defined component by component by
\begin{equation*}
\begin{aligned}
& \big( \bA \otimes \bB \big)_{i_1,...i_p,j_1,...,j_q} = a_{i_1,...,i_p} \, b_{j_1,...,j_q} ~,\\
& \big( \bA \cdot \bB \big)_{i_1,...i_{p-1},j_2,...,j_q} = a_{i_1,...,i_{p-1},k} \, b_{k,j_2,...,j_q} ~,\\
& \big( \bA : \bB \big)_{i_1,...i_{p-2},j_3,...,j_q} = a_{i_1,...,i_{p-2},k,\ell} \, b_{k,\ell,j_{s+1},...,j_q} ~,
\end{aligned}
\end{equation*}
where we use the Einstein convention for the summations with indexes~$k$ and~$\ell$.
\item[$-$] Note also that all these products are inner products on the set of the $p$-tensors. It allows us to define a generalized Froebenius norm:
\begin{equation*}
|\bA|^2 = \sum_{i_1,...,i_p} a_{i_1,...,i_p}^2.
\end{equation*}
\end{itemize}
\noindent
We conclude this section introducing the constant~$C$. 
This constant stands for any constant depending on the data of the problem: initial conditions, physical parameters...
In some cases, informations will be given on the dependence of this constant (see for example Section~\ref{part:proof3} where we explain that this constant may depend on time~$t^\star$ but must remain bounded when~$t^\star$ is bounded).

\subsection{Assumptions}\label{sec:assumption}

The results proved in this article requires some assumptions about the data.
In addition to the assumptions on the regularity of the initial conditions that will be specified in each theorem statement, we will need some "natural" assumptions.
\par\noindent
$\checkmark$
The first assumption relates to the initial deformation~$\bG_0$:
\begin{equation}\label{30.1}
\exists \gamma>0 ~;~ \det \bG_0 \geq \gamma.
\end{equation}
We note that in many applications, the fluid is assumed to be initially quiescent. In that case, we have $\bG_0=\bdelta$ and $\det\bG_0=1$. Moreover, we will see (equation~\eqref{27} in the preliminary section~\ref{positive-deformation}) that the quantity $\det\bG$ is only convected by the flow. If the fluid is assumed to be at rest in the past (that is for $T$ large enough), then we always have  $\det\bG_0=1$. The assumption on the positiveness of~$\det\bG_0$ allows us consider, for instance, such cases.
\par\noindent
$\checkmark$
The second assumption relates to the initial memory $m_0=-\partial_T K_0$:
\begin{subequations} \label{assumptionK}
\begin{align}
& m_0 \geq 0, \label{30.2} \\
& \partial_T m_0 \leq 0. \label{30.3}
\end{align}
\end{subequations}
The assumption~\eqref{30.2} corresponds to the fact that the quantity $m_0$ describes the memory of the fluid, that is the weight that must have the quantity~$\mathscr S$ in the flow {\it via} the relation~\eqref{eq:4}. It is physically positive. The assumption~\eqref{30.3} indicates that the memory decreases with the age~$T$: It is linked to the principle of fading memory, see~\cite{Coleman}.
\par\noindent
In the integral models, that is to say when the memory is explicitly given in terms of age~$T$, it is a combination of exponentially decreasing functions (see for instance the Doi-Edwards model under the I.A. approximation, subsection~\ref{subsec:IA} and more precisely the expression~\eqref{11} of the memory).
Such decreasing behaviors will be prescribed in the functional spaces with exponential weight. For example, we will impose that there exists $\mu>0$ such that $m_0 \leq C\, \mathrm e^{-\We\, \mu T}$.

\subsection{Main results}

The first result concerns an existence result for strong solution. It is a local in time result:\\[-0.2cm]
\begin{theorem}[local existence]\label{th:local}
Let $r\in ]2,+\infty[$, $q\in]d,+\infty[$ and $\mu>0$.
\par\noindent
If the data $\bv_0$, $\bG_0$ and $K_0$ satisfy the assumptions~\eqref{30.1}, \eqref{30.2} and have the following regularity
\begin{equation*}\label{31}
\begin{aligned}
& \bv_0 \in \mathscr D^r_q, 
&& \bG_0\in L^\infty_T W^{1,q}_{\bx},
&& \partial_T \bG_0\in L^\infty_T L^q_{\bx}, \\
& \mathrm e^{\We\, \mu T} \partial_T K_0 \in L^\infty_T L^\infty_{\bx,s},
&& \mathrm e^{\We\, \mu T/2} \partial_T \nabla K_0 \in L^2_T L^q_{\bx,s},
&& \partial_s K_0 \in L^\infty_T L^2_{\bx,s},
\end{aligned}
\end{equation*}
then there exists $t^\star >0$ and a strong solution $(\bu,\bG,K)$ to System~\eqref{system} in $[0,t^\star]$, which satisfies the initial and boundary conditions~\eqref{7}, \eqref{8} and~\eqref{9}.
Moreover we have
\begin{equation*}\label{32}
\begin{array}{ll}
\bv \in L^r(0,t^\star;W^{2,q}_{\bx}), \phantom{totototototototototototototo}
& \partial_t \bv \in L^r(0,t^\star;L^q_{\bx}), \\
\bG \in L^\infty(0,t^\star,L^\infty_T W^{1,q}_{\bx}),
& \partial_s \bG, ~\partial_t  \bG \in L^r(0,t^\star;L^\infty_T L^q_{\bx}),\\
K, \partial_T K, \mathrm e^{\mu \, (\We \, T-t)} \partial_T K \in L^\infty(0,t^\star;L^\infty_T L^\infty_{\bx,s}),
& \partial_t K \in L^2(0,t^\star;L^\infty_T L^2_{\bx,s}),\\
\partial_T \nabla K \in L^\infty(0,t^\star;L^1_T L^q_{\bx,s}\cap L^2_T L^q_{\bx,s}),
& \partial_s K \in L^\infty(0,t^\star;L^\infty_T L^2_{\bx,s}),
\end{array}
\end{equation*}
and the memory~$m=-\partial_T K$ remains non negative.
\end{theorem}\\[-0.2cm]
\begin{remark}
In this article, we will not give any result on the pressure~$p$. In practice, the latter is regarded as a Lagrange multiplier associated to the divergence free constraint.
It can be solved using the Riesz transforms. More precisely, taking the divergence of the first equation~\eqref{eq:1} of System~\eqref{system} we use the periodic boundary conditions to have
\begin{equation}\label{pressure}
p=-(-\Delta)^{-1}\div\div\, (\bsigma-\bv\otimes\bv).
\end{equation}
From Theorem~\ref{th:local}, we can prove that the solutions of System~\eqref{system} discussed in this paper have $\bsigma-\bv \otimes \bv$ in $L^\infty(0,t^\star;L^2_{\bx})$. The pressure in the solution of~\eqref{system} is meant to be given by~\eqref{pressure}.
\end{remark}\\[0.2cm]
We will see during the proof of this theorem~\ref{th:local} that one of the key point is the behavior of the memory $m=-\partial_T K$ for large value of the age~$T$: if the memory is exponentially decreasing at $t=0$ (with respect to the age variable~$T$) then the solution will be exponentially decreasing for any time~$t>0$.
\par\noindent
In the same way, it is possible to prove that if the memory~$m$ is initially decreasing\footnote{Be careful not to confuse the terms "exponentially decreasing" and "decreasing". The first means that $m$ is bounded by a function of the form $\mathrm e^{-T}$, while the second means that $\partial_T m \leq 0$. Moreover the first is a global property, while the second is a local property.} (with respect to the age variable~$T$) then the solution will be decreasing for any time~$t>0$.
The proof - which is not given in the proof of Theorem~\ref{th:local} - consists in derivating twice the equation~\eqref{eq:6} with respect to~$T$, and next in applying the maximum principle (see the subsection~\ref{max-principle}, page~\pageref{max-principle}) to the function~$\partial_Tm$.
Therefore, the assumption~\eqref{30.3}, which is not necessary to obtain local existence, is also preserved in time.
\par\noindent
We will show that the solution obtained in Theorem~\ref{th:local} with this additional assumption~\eqref{30.3} is the only one in the class of regular solutions.
Precisely, the result reads as follows.\\[-0.2cm]
\begin{theorem}[uniqueness]\label{th:unic}
Let $t^\star>0$.
\par\noindent
Let $(\bu_1,\bG_1,K_1)$ and $(\bu_2,\bG_2,K_2)$ be two solutions to System~\eqref{system} satisfying the initial and boundary conditions~\eqref{7}, \eqref{8} and~\eqref{9}.
If we have, for $i\in \{1,2\}$,
\begin{equation}\label{33}
\begin{aligned}
& \nabla \bv_i \in L^2(0,t^\star;L^\infty_{\bx}), \\
& \bG_i \in L^2(0,t^\star;L^\infty_T (L^\infty_{\bx} \cap W^{1,d}_{\bx})), \\
& \partial_T K_i \in L^\infty(0,t^\star;L^\infty_T L^\infty_{\bx,s} \cap L^1_T L^\infty_{\bx,s}), \\
& \partial_T \nabla K_i \in L^\infty(0,t^\star;L^2_T L^d_{\bx,s}),
\end{aligned}
\end{equation}
and if each~$m_i=-\partial_T K_i$ is decreasing with respect to~$T$ then the two solutions coincide.
\end{theorem}\\[0.2cm]
Obviously, the solution obtained in Theorem~\ref{th:local} satisfies the regularity requested in~\eqref{33}. Combining Theorems~\ref{th:local} and~\ref{th:unic} we get a local and uniqueness result.
In the two-dimensional case, it is possible to show that the solution $(\bv,\bG,K)$ of problem~\eqref{system} exists for any time~$t^\star>0$. More precisely we have the following result:\\[-0.2cm]
\begin{theorem}[global existence in 2D]\label{th:global}
Let $r\in ]2,+\infty[$, $q\in]2,+\infty[$ and $\mu>0$.
\par\noindent
We assume that $\frac{1}{r} + \frac{1}{q} < \frac{1}{2}$ and that the data $\bv_0$, $\bG_0$ and $K_0$ satisfy the same assumptions that in Theorem~\ref{th:local}.
\par\noindent
Let $t^\star>0$ be arbitrary.
\par\noindent
There exists a constant $C$ depending only on the data with~$C$ bounded for bounded~$t^\star$, and a solution $(\bv,\bG,K)$ of~\eqref{system} satisfying the initial and boundary conditions~\eqref{7}, \eqref{8} and~\eqref{9} such that
\begin{equation}\label{33.1}
\begin{aligned}
& \|\nabla^2 \bv\|_{L^r(0,t^\star;L^q_{\bx})} \leq C, \qquad && \|\nabla \bv\|_{L^\infty(0,t^\star;L^\infty_{\bx})} \leq C, \\
& \|\nabla \bS\|_{L^r(0,t^\star;L^q_{\bx,s})} \leq C, \qquad && \|\bS\|_{L^\infty(0,t^\star;L^\infty_{\bx})} \leq C.
\end{aligned}
\end{equation}
\end{theorem}
The estimates~\eqref{33.1} announced in the theorem~\ref{th:global} above are sufficient to prove a global in time existence of a solution~$(\bv,\bG,K)$. In fact if~\eqref{33.1} holds then it is possible to prove - principally using the lemmas introduced in the proof of the local existence result - that the solution $(\bv,\bG,K)$ of~\eqref{system} at time $t^\star$ have the same regularity that at time~$0$. Applying the local result (theorem~\ref{th:local}), we deduce that the solution can not blow up in finite time.

\section{Preliminaries}\label{part:preliminaries}

In this section we give some results which will be using during the different proofs of the previous theorems.

\subsection{Some bounds for the function~$\mathscr S$}

One of the key points of the proof of global existence lies in the fact that the stress~$\bsigma$, defined by~\eqref{eq:3}-\eqref{eq:4} is bounded. The first result in this direction is the following result concerning the function~$\mathscr S$:\\[-0.2cm]
\begin{proposition}\label{propS}
The function~$\mathscr S$ defined by the relation~\eqref{6} is of class~$\mathscr C^1$ on $\mathcal L(\R^d) \setminus \{\b0 \}$ and satisfies the following properties:
\begin{subequations}
  \begin{align}
& \exists \mathscr S_\infty\geq 0 ~;~\forall \bG\in \mathcal L(\R^d) \setminus \{\b0 \} ~, ~\displaystyle |\mathscr S(\bG)|\leq \mathscr S_\infty; \label{16} \\[0.2cm]
& \exists \mathscr S_\infty' \geq 0 ~;~ \forall \bG\in \mathcal L(\R^d) \setminus \{\b0 \} ~, ~ |\bG| |\mathscr S'(\bG)|\leq \mathscr S_\infty'.\label{17}
\end{align}
\end{subequations}
\end{proposition}
\proof
Recall the definition~\eqref{6} of the function~$\mathscr S$:
\begin{equation*}\label{18}
\mathscr S(\bG) = \frac{1}{\langle |\bG \cdot \bu| \rangle_0} \, \Big\langle \frac{(\bG \cdot \bu)\otimes (\bG \cdot \bu)}{|\bG \cdot \bu|} \Big\rangle_0 - \frac{1}{d} \bdelta.
\end{equation*}
$\checkmark$ We first notice that the following inequality is obvious
\begin{equation}\label{19}
|\langle f \rangle_0| = \Big| \int_{\mathbb S^{d-1}} f(\bu) \, \mathrm d\bu \Big| \leq \int_{\mathbb S^{d-1}} |f(\bu)| \, \mathrm d\bu = \langle |f| \rangle_0,
\end{equation}
so that the first point of the proposition~\ref{propS} is a direct consequence of the inequality $|\mathbf A \otimes \mathbf B| \leq |\mathbf A|\,|\mathbf B|$ for all tensors~$\mathbf A$ and~$\mathbf B$, and of the relation $|\bdelta| = \sqrt{d}$. More precisely we obtain, for all $\bG\in \mathcal L(\R^d) \setminus \{\b0 \}$:
\begin{equation*}\label{20}
|\mathscr S(\bG)| \leq 1+\frac{1}{\sqrt{d}},
\end{equation*}
which corresponds to~\eqref{16}.
\par\noindent
$\checkmark$ For the second point~\eqref{17}, we write the tensor~$\mathscr S(\bG)$ component by component: for any $(k,\ell)\in \{1,...,d\}^2$,
$$\mathscr S(\bG)_{k\ell} = \frac{1}{\langle |\bG \cdot \bu| \rangle_0} \, \Big\langle \frac{(\bG \cdot \bu)_k\, (\bG \cdot \bu)_\ell}{|\bG \cdot \bu|} \Big\rangle_0 - \frac{1}{d} \bdelta_{k\ell}.$$
It is then not difficult to evaluate each components of the derivative\footnote{The application~$\mathscr S$ being defined on an open set of~$\mathcal L(\R^d)$ with values in~$\mathcal L(\R^d)$, its differential is an application with values in $\mathcal L(\mathcal L(\R^d),\mathcal L(\R^d))$. Consequently, $\mathscr S'(\bG)$ can be identify to a $4$-order tensor whose components are $\mathscr S'(\bG)_{i,j,k,\ell} = \partial_{\bG_{ij}} \mathscr S(\bG)_{k\ell}$.} tensor $\mathscr S'(\bG)$.
For instance for $(i,j,k,\ell)\in \{1,...,d\}^4$ we have
\begin{equation*}\label{21}
\partial_{\bG_{ij}} \big( |\bG \cdot \bu | \big) = \frac{(\bG \cdot \bu)_i \, \bu_j}{|\bG \cdot \bu |}.
\end{equation*}
More generally we obtain
\begin{equation}\label{22}
\begin{aligned}
\partial_{\bG_{ij}} \mathscr S(\bG)_{k\ell} 
=
& \frac{-1}{ \langle |\bG \cdot \bu| \rangle_0^2} \Big\langle \frac{(\bG \cdot \bu)_i \, \bu_j}{|\bG \cdot \bu|} \Big\rangle_0 \Big\langle \frac{(\bG \cdot \bu)_k \, (\bG \cdot \bu)_{\ell}}{|\bG \cdot \bu|} \Big\rangle_0 \\
& + \frac{1}{\langle |\bG \cdot \bu| \rangle_0} \Big\langle \frac{\bdelta_{ik} \bu_j (\bG \cdot \bu)_\ell + \bdelta_{i\ell} \bu_j (\bG \cdot \bu)_k}{|\bG \cdot \bu|} \Big\rangle_0 \\
& + \frac{-1}{\langle |\bG \cdot \bu| \rangle_0} \Big\langle \frac{(\bG \cdot \bu)_k (\bG \cdot \bu)_{\ell} (\bG \cdot \bu)_i \bu_j}{|\bG \cdot \bu|^3} \Big\rangle_0.
\end{aligned}
\end{equation}
Taking the norm in equality~\eqref{22} and using~\eqref{19}, we deduce
\begin{equation}\label{23}
|\mathscr S'(\bG)| \leq 2(1+\sqrt{d}) \frac{\langle |\bu| \rangle_0}{\langle |\bG \cdot \bu| \rangle_0}.
\end{equation}
We note that $\langle |\bu| \rangle_0 = \langle 1 \rangle_0 = \frac{2\pi^{d/2}}{\Gamma(d/2)}$. We next use the fact that the application $\bG \mapsto \langle |\bG \cdot \bu| \rangle_0$ is a norm on the finite dimensional space~$\mathcal L(\R^d)$. This norm is then equivalent to  $\bG \mapsto |\bG|$: there exists a constant~$C_d$ such that for all $\bG\in \mathcal L(\R^d)$ we have
\begin{equation*}\label{24}
\langle |\bG \cdot \bu| \rangle_0 \geq C_d \, |\bG|.
\end{equation*}
The inequality~\eqref{23} becomes
\begin{equation*}\label{25}
|\bG| |\mathscr S'(\bG)| \leq 2(1+\sqrt{d}) \frac{2\pi^{d/2}}{C_d \, \Gamma(d/2)},
\end{equation*}
that concludes the proof of the Proposition~\ref{propS}.
\cqfd

\subsection{Positive norm for the deformation~$\bG$}\label{positive-deformation}

In practice, we will prove that for any solutions $(\bv,\bG,K)$ to the System~\eqref{system} the deformation gradient tensor~$\bG$ has a positive norm. More precisely we have\\[-0.2cm]
\begin{lemma}\label{lemma7}
Let~$\bv$ be a free divergence vector field on~$\T$ and~$\bG$ be a solution to Equation~\eqref{eq:5} with initial conditions satisfy~$\det \bG|_{t=0} \geq \gamma > 0$ and $ \bG|_{T=0}=\bdelta$.
\par\noindent
There exists a constant $\widetilde \gamma$ such that for all $(t,T,\bx)\in (0,t^\star) \times \R^+ \times \T$ we have
\begin{equation*}\label{26}
|\bG(t,T,\bx)| \geq \widetilde \gamma >0.
\end{equation*}
\end{lemma}
\proof
A simple calculation shows that the quantity~$\det(\bG)$ satisfies
\begin{equation}\label{27}
\mathcal D \det(\bG) =  \div\, \bv \, \det(\bG)= 0,
\end{equation}
where~$\mathcal D$ refers to the one order derivating operator~$\mathcal D=d_t + \frac{1}{\We}\partial_T$.
The value~$\det(\bG)$ is then constant along the characteristic lines.
Since all the characteristic lines start from the lines $\{t=0\}$ or $\{T=0\}$ we deduce from the assumptions that $\det(\bG)\geq \min(\gamma,1)$ on $(0,t^\star) \times \R^+ \times \T$.
\par
\begin{figure}[htbp]
\begin{center}
\includegraphics[height=3cm]{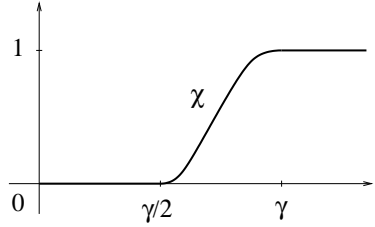}
\caption{An example of truncation function~$\chi$.}\label{fig2}
\end{center}
\end{figure}

\noindent
Due to the inequality of arithmetic and geometric means, we have 
\begin{equation*}\label{28}
|\bG|^2 = \mathrm{Tr}(\transp{\bG}\cdot\bG) \geq 2\sqrt{\det(\transp{\bG}\cdot\bG)} = 2|\det(\bG)| \geq 2\min(\gamma,1).
\end{equation*}
That concludes the proof of Lemma~\ref{lemma7} taking $\widetilde \gamma=\sqrt{2\min(\gamma,1)}$.
\cqfd
\noindent
For the local result we will use the fact that the derivative $\mathscr S'(\bG)$ is bounded. That is clearly not the case if we only use the Proposition~\ref{propS}. We then introduce the function
\begin{equation*}\label{29}
\widetilde{\mathscr S}(\bG) = \mathscr S(\bG) \, \chi(|\bG|),
\end{equation*}
where $\chi\in \mathscr C^\infty(\R^+,[0,1])$ satisfies $\chi|_{[0,\frac{\widetilde \gamma}{2}]}=0$, $\chi|_{[\widetilde \gamma,+\infty[}=1$ and $\chi'\leq \frac{3}{\widetilde \gamma}$ (see the Figure~\ref{fig2}).

\noindent
Using the Proposition~\ref{propS} it suffices to derivate~$\widetilde{\mathscr S}$:
\begin{equation*}\label{29.8}
\widetilde{\mathscr S}\, '(\bG) = \mathscr S'(\bG) \, \chi(|\bG|) + \mathscr S(\bG) \otimes \frac{\bG}{|\bG|}\, \chi'(|\bG|),
\end{equation*}
to deduce
\begin{equation*}\label{30}
\forall \bG\in \mathcal L(\R^d) \qquad
|\widetilde{\mathscr S}(\bG)|\leq \mathscr S_\infty
\quad \text{and} \quad
|\widetilde{\mathscr S}\; '(\bG)|\leq \widetilde{\mathscr S_\infty'},
\end{equation*}
where 
$\widetilde{\mathscr S_\infty'} = \frac{2}{\widetilde \gamma} \mathscr S_\infty' + \frac{3}{\widetilde \gamma} \mathscr S_\infty$. Using the Lemma~\ref{lemma7}, we have the following consequence:\\[-0.2cm]
\begin{lemma}\label{lemma8}
Under the assumption~\eqref{30.1}, the solution of the System~\eqref{system} is the same if we use $\widetilde{\mathscr S}$ instead of~$\mathscr S$.
In other words, we can assume that $\mathscr S'(\bG)$ is bounded on $(0,t^\star)\times \R^+ \times \T$.
\end{lemma}

\subsection{A Gronwall lemma with two times}\label{subsection-gronwall}

The choice that was made in this article is to use the current time~$t$ and another time~$T$ corresponding to the age of the flow. This choice simplifies the expression of the orientation tensor~$\bS$ since the t-dependence does not appear in the integral bounds (see~\eqref{eq:4}). The price to pay is that the derivatives in time in the evolution equations of~$K$ and~$\bG$ involve both~$t$ and~$T$. In this framework, the following lemma is the analog of the classical Gronwall lemma for functions depending only~$t$. Its proof is based on a change of variable that easily allows  to follow the characteristics. It is proved in~\cite{Chupin3}.
\begin{lemma}\label{lem:gronwall}
Let $f:\R^+ \mapsto \R^+$ a positive and locally integrable function.
If a function $y:\R^+ \! \times \R^+ \mapsto \R$ satisfies, for all $(t,T)\in (0,t^\star) \! \times \R^+$:
\begin{equation*}\label{eq:lem1}
\partial_t y(t,T) + \frac{1}{\We} \partial_T y(t,T) \leq f(t)\, y(t,T)
\end{equation*}
then we have, for all $(t,T)\in (0,t^\star) \! \times \R^+$:
\begin{equation*}\label{eq:lem2}
y(t,T) \leq  \zeta(t,T) \, \mathrm{exp}\bigg( \int_0^t f(t')\, \mathrm dt' \bigg),
\end{equation*}
where $\displaystyle \zeta(t,T) = \left\{
\begin{aligned}
& y\big( T-\frac{t}{\We},0\big) \quad \text{if $t\leq \We\, T$}, \\
& y(0,t-\We\, T) \quad \text{if $t > \We\, T$}.
\end{aligned}
\right.
$
\end{lemma}

\subsection{A maximum principle}\label{max-principle}

The memory of the fluid is described by the function~$K$. This function satisfies the equation~\eqref{eq:6}. The form of the equation~\eqref{eq:6} is particular in that it checks a "maximum principle". This will be used repeatedly in the following sections.
\begin{lemma}\label{lem:max}
If $\partial_s g \in L^1(0,t^\star; L^\infty_{\bx,s})$ and $\bv$ is free divergence on~$\T$ then the solution $f(t,T,\bx,s)$ to the following system
\begin{equation}\label{236}
\left\{
\begin{aligned}
& \mathrm d_t f + \frac{1}{\We} \partial_T f + g \, \partial_s f - \frac{1}{\We} \partial_s^2 f = 0, \\
& f\big|_{t=0} = f_0, \qquad
f\big|_{T=0} = f_1, \qquad
f\big|_{s=-\frac{1}{2}} = f\big|_{s=\frac{1}{2}} = 0,
\end{aligned}
\right.
\end{equation}
satisfies the following maximum principle on $(0,t^\star)\times \R^+ \times \T \times (-\frac{1}{2},\frac{1}{2})$:
\begin{equation*}\label{237}
\min\{ \inf_{T,\bx,s} f_0 , \inf_{t,\bx,s} f_1 \} \leq f \leq \max\{ \sup_{T,\bx,s} f_0 , \sup_{t,\bx,s} f_1 \}.
\end{equation*}
\end{lemma}
\proof
Considering $f-\min\{ \inf f_0 , \inf f_1 \}$ or $\max\{ \sup f_0 , \sup f_1 \}-f$ instead of~$f$, it suffices to show that the solution~$f$ is non negative if the data~$f_0$ and~$f_1$ are non negative.\\
We multiply the first equation of~\eqref{236} by the negative part~$f^-$ of~$f$, and we integrate with respect to~$\bx$ and~$s$.
It is important to notice that the function~$f$ does not necessary satisfy the boundary conditions $f\big|_{s=\frac{1}{2}} = f\big|_{s=-\frac{1}{2}} = 0$, but that its negative part~$f^-$ satisfies these conditions. Using integrations by parts, we deduce
\begin{equation*}\label{238}
\partial_t \|f^-\|_{L^2_{\bx,s}}^2 + \frac{1}{\We} \partial_T \|f^-\|_{L^2_{\bx,s}}^2 + \frac{2}{\We} \|\partial_s f^-\|_{L^2_{\bx,s}}^2 = \int_\T \int_{-\frac{1}{2}}^{\frac{1}{2}} \partial_s g \, |f^-|^2.
\end{equation*}
Since $\partial_s g \in L^1_t L^\infty_{\bx,s}$ we obtain
\begin{equation*}\label{239}
\partial_t \|f^-\|_{L^2_{\bx,s}}^2 + \frac{1}{\We} \partial_T \|f^-\|_{L^2_{\bx,s}}^2 \leq \|\partial_s g\|_{L^\infty_{\bx,s}} \|f^-\|_{L^2_{\bx,s}}^2.
\end{equation*}
The Gronwall lemma with two variables (see the lemma~\ref{lem:gronwall}) implies that
\begin{equation*}\label{240}
\|f^-\|_{L^2_{\bx,s}}^2(t,T) \leq \zeta(t,T) \, \mathrm{exp} \big(\|\partial_s g\|_{L^1(0,t; L^\infty_{\bx,s})} \big),
\end{equation*}
where the function~$\zeta$ only depends on the values of~$\|f^-\|_{L^2_{\bx,s}}^2$ on the boundaries~$\{t=0\}$ and $\{T=0\}$. In the case where the data~$f_0$ and $f_1$ are non negative we have $\zeta=0$ and we deduce that~$\|f^-\|_{L^2_{\bx,s}}^2=0$. We conclude that~$f$ is non negative too.
\cqfd

\subsection{A Cauchy problem involved in the proof of global existence}\label{subsection-ode}

In this section, we are interested in the following Cauchy problem
\begin{equation}\label{241}
\left\{
\begin{aligned}
& y"(x) = \xi_0 \, (y(x)+\xi_2)^k (y'(x))^2, \\
& y(0)=0, \quad y'(0)=\xi_1,
\end{aligned}
\right.
\end{equation}
where~$\xi_0$, $\xi_1$, $\xi_2$ and~$k$ are positive constants.\\
This problem will occur during the proof of global existence theorem, section~\ref{part:proof3}. More precisely, we will see that with a good choice of parameters~$\xi_0$, $\xi_1$, $\xi_2$ and~$k$, we can control $y(-\mathrm e^{\We \mu T} \partial_T K)$ for arbitrarily long time. This will allow to control the stress~$\bsigma$.\\
Although this equation~\eqref{241} is nonlinear and of order~$2$, we can explicitly give the solution.
We will see that its expression makes appear the function $F:\R^+ \longrightarrow \R^+$ defined by
\begin{equation}\label{241.5}
F(X) = \int_0^X  \mathrm e^{-\frac{\xi_0}{k+1} (x+\xi_2)^{k+1}}\, \mathrm dx.
\end{equation}
Clearly this function is one-to-one (increasing) from $\R^+$ to $[0,\ell[$, where the real $\ell$ denotes the limit:
\begin{equation*}\label{241.6}
\ell = \int_0^{+\infty}  \mathrm e^{-\frac{\xi_0}{k+1} (x+\xi_2)^{k+1}}\, \mathrm dx.
\end{equation*}

\begin{proposition}\label{prop:4}
The Cauchy problem~\eqref{241} admits a unique solution given by
\begin{equation}\label{242}
\begin{aligned}
& y(x)=F^{-1}(\xi_1 \mathrm e^{-\frac{\xi_0 \xi_2^{k+1}}{k+1}}x)
\qquad  \text{for all} \quad x\in \Big[0, \frac{\ell}{\xi_1} \mathrm e^{\frac{\xi_0 \xi_2^{k+1}}{k+1}}\Big[.
\end{aligned}
\end{equation}
\end{proposition}

\proof
The equation~\eqref{241} is a two order ordinary differential equation and it is possible to apply the Cauchy-Lipschitz theorem: there exists a unique local solution. Moreover, we clearly have $y">0$ and then $y'\geq \xi_1>0$. The  equation~\eqref{241} also write
\begin{equation*}\label{243}
\Big[ \ln(y') \Big]' = \Big[ \frac{\xi_0}{k+1} (y+\xi_2)^{k+1}  \Big]'.
\end{equation*}
Using boundary conditions given in~\eqref{241} we obtain the following first order ordinary differential equation:
\begin{equation*}\label{243.5}
y' = \xi_1 \mathrm e^{\frac{\xi_0 \xi_2^{k+1}}{k+1}} \mathrm e^{\frac{\xi_0}{k+1} (y+\xi_2)^{k+1}}.
\end{equation*}
Making appear the function $F$, we deduce that
\begin{equation*}\label{243.6}
\Big[F(y)\Big]' = \xi_1 \mathrm e^{\frac{\xi_0 \xi_2^{k+1}}{k+1}}.
\end{equation*}
Since $F(0)=0$ and $y(0)=0$ we integrate and deduce
\begin{equation*}\label{243.7}
F(y(x)) = \xi_1 \mathrm e^{\frac{\xi_0 \xi_2^{k+1}}{k+1}} x.
\end{equation*}
The expression~\eqref{242} given in Proposition~\ref{prop:4} follows since~$F$ is one-to-one.
\cqfd

\section{Proof of the local existence result: theorem~\ref{th:local}}\label{part:proof1}

\subsection{Strategy: a point fixed formulation}

In order to prove the local existence result, we rewrite the set of equations~\eqref{system} as a fixed point system.
More precisely, we consider the mapping (the spaces will be further given)
\begin{equation*}\label{35}
\Phi ~:~ (\overline{\bv} , \overline{\bG} , \overline{K}) \longmapsto ~(\bv,\bG,K),
\end{equation*}
defined as follows:
\par
$\checkmark$ {\bf Velocity problem} -- The velocity~$\bv(t,\bx)$ is the solution of the following Stokes problem
\begin{equation}\label{36}
\left\{
\begin{aligned}
& \Re \, \partial_t \bv + \nabla p - (1-\omega) \Delta \bv = \bbf, \\
& \div\, \bv = 0,\\
& \bv|_{t=0} = \bv_0,
\end{aligned}
\right.
\end{equation}
where the source term~$\bbf$ contains the nonlinear term of the Navier-Stokes equations and the term coupling velocity and stress, namely $\bbf = -\Re\, \overline{\bv} \cdot \nabla \overline{\bv} + \div\, \bsigma$. The stress~$\bsigma(t,\bx)$ is defined by
\begin{equation}\label{37}
\bsigma(t,\bx) = \omega \int_{-\frac{1}{2}}^{\frac{1}{2}} \bS(t,\bx,s)\, \mathrm ds,
\end{equation}
where the orientation tensor~$\bS(t,\bx,s)$ is given by
\begin{equation}\label{38}
\bS(t,\bx,s) = -\int_0^{+\infty} \partial_{T} \overline K(t,T,\bx,s)\, \mathscr S(\overline \bG(t,T,\bx))\, \mathrm dT.
\end{equation}

$\checkmark$ {\bf Deformation problem} -- The deformation gradient tensor~$\bG(t,T,\bx)$ is the solution of the linear equation
\begin{equation}\label{39}
\left\{
\begin{aligned}
& \partial_t \bG + \overline \bv \cdot \nabla \bG + \frac{1}{\We} \partial_T \bG = \bG\cdot \nabla \overline \bv,\\
& \bG|_{t=0} = \bG_0, \qquad \bG|_{T=0} = \bdelta.
\end{aligned}
\right.
\end{equation}

$\checkmark$ {\bf Memory problem} -- The scalar quantity~$K(t,T,\bx,s)$ is the solution of the following linear equation
\begin{equation}\label{40}
\left\{
\begin{aligned}
& \partial_t K + \overline \bv \cdot \nabla K + \frac{1}{\We} \partial_T K + g \, \partial_s K - \frac{1}{\We} \partial^2_s K = 0,\\
& K|_{s=-\frac{1}{2}} = K|_{s=\frac{1}{2}} = 0,\\
& K|_{t=0} = K_0, \qquad K|_{T=0} = 1,
\end{aligned}
\right.
\end{equation}
where~$g$, which only depends on~$t$, $\bx$ and~$s$, is given by
\begin{equation*}\label{41}
g(t,\bx,s) = \nabla \overline \bv(t,\bx) : \int_0^s \bS(t,\bx,s')\, \mathrm ds'.
\end{equation*}

The goal of the next subsections is to analyze these problems~\eqref{36}, \eqref{39} and~\eqref{40} independently. We will see at the end of this analyze that it is possible to apply the Schauder fixed point theorem for the function~$\Phi$ with adapted functional spaces in order to deduce the result of Theorem~\ref{th:local}.

\subsection{Estimates for the velocity~$\bv$ solution of a Stokes problem~\eqref{36}}

The results for the Stokes system~\eqref{36} are very numerous.
In this subsection we only recall, without proof (we can found a proof in~\cite{Giga}), a well known result for the time dependent Stokes problem.
In order to simplify expressions, we use the following norm on the velocity field:
\begin{equation*}\label{42}
\vvvert \bv \vvvert_1 := \|\bv\|_{L^r(0,t^\star;W^{2,q}_{\bx})} + \|\partial_t \bv\|_{L^r(0,t^\star;L^q_{\bx})}.
\end{equation*}
The result on the Stokes problem~\eqref{36} states as follows:\\[-0.2cm]

\begin{lemma}\label{lem:stokes}
Let $t^\star>0$, $r\in ]1,+\infty[$ and $q\in]1,+\infty[$.
\par\noindent
If $\bv_0 \in \mathscr D^r_q$ and $\bbf \in L^r(0,t^\star;L^q_{\bx})$ then there is a unique solution $\bv\in L^r(0,t^\star;\mathscr D(A_q))$ such that $\partial_t\bv \in L^r(0,t^\star;H_q)$ to System~\eqref{36}.
This solution satisfies
\begin{equation*}\label{43}
\vvvert \bv \vvvert_1 \leq F_1\big( \|\bbf\|_{L^r(0,t^\star;L^q_{\bx})} \big),
\end{equation*}
where the function~$F_1$ depends on~$r$, $q$, $\omega$, $\Re$ and the initial value~$\bv_0$. Moreover this function is continuous and nondecreasing on~$\R^+$.
\end{lemma}\\[0.2cm]
\noindent
In practice, the function~$F_1$ may be chosen as (see~\cite{Chupin3})
\begin{equation*}\label{44}
F_1(X) = \frac{C_1(r,q)}{1-\omega} \big( \Re \, \|\bv_0\|_{W^{2,q}_{\bx}} + X \big).
\end{equation*}

\subsection{Estimates for the deformation gradient~$\bG$, solution of~\eqref{39}}

The existence and regularity for the deformation gradient~$\bG$ is less classical.
As previously, we introduce an adapted norm, namely for the deformation gradient
\begin{equation}\label{45}
\vvvert \bG \vvvert_2 := \|\bG\|_{L^\infty(0,t^\star;L^\infty_T W^{1,q}_{\bx})} + \|\partial_t \bG\|_{L^r(0,t^\star;L^\infty_T L^q_{\bx})} + \|\partial_T \bG\|_{L^r(0,t^\star;L^\infty_T L^q_{\bx})}.
\end{equation}
\begin{lemma}\label{lem:1116}
Let $0<t^\star<1$, $r\in ]1,+\infty[$ and $q\in]d,+\infty[$.
\par\noindent
If $\bG_0\in L^\infty_T W^{1,q}_{\bx}$, $\partial_T\bG_0\in L^\infty_T L^q_{\bx}$ and $\overline \bv\in L^\infty(0,t^\star;W^{1,q}_{\bx}) \cap L^1(0,t^\star;W^{2,q}_{\bx})$ is free divergence then the problem~\eqref{39} admits a unique solution $\bG \in L^\infty(0,t^\star;L^\infty_T W^{1,q}_{\bx})$ such that~$\partial_t \bG$ and~$\partial_T \bG$ belongs to~$L^r(0,t^\star;L^\infty_T L^q_{\bx})$.
This solution satisfies
\begin{equation*}
\vvvert \bG \vvvert_2 \leq F_2\big( \|\overline \bv\|_{L^\infty(0,t^\star;W^{1,q}_{\bx})\cap L^1(0,t^\star;W^{2,q}_{\bx})} \big),
\end{equation*}
where the function~$F_2$ depends on~$r$, $q$, $\We$ and the initial value~$\bG_0$. Moreover this function is continuous and nondecreasing on~$\R^+$.
\end{lemma}\\[0.2cm]
\noindent
The proof of a very similar result is given in~\cite{Chupin3}. 
One of the differences is that we show here $\partial_t \bG \in L^r(0,t^\star;L^\infty_T L^q_{\bx})$ and not only $\partial_t \bG \in L^r(0,t^\star;L^q_{\bx})$ for any $T\in \R^+$. This difference takes its significance when we will give sense to initial conditions, see the remark~\ref{rem:0944}.\\
\proof
The existence of a unique solution to~\eqref{39} follows from the application of the method of characteristics (see~\cite[Appendix p.~26]{Fernandez-Guillen-Ortega}).
In practice, the following estimates will be made on regular solution~$\bG_n$ which approaches the solution~$\bG$ when a regular velocity field~$\bv_n$ approaches the velocity~$\bv$.
The regularity of these solutions~$\bG_n$ with respect to~$t$ and~$T$ comes from the Cauchy-Lipschitz theorem.
For sake of simplicity, we omit the indexes~''$n$''.
In the following proof, we refer to~\cite{Fernandez-Guillen-Ortega} for the passage to the limit $n\to +\infty$.
The rest of the proof of Lemma~\ref{lem:1116} is split into three parts: in the first one (see the subsection~\ref{part:estimates-1}) we obtain a first estimate concerning the regularity of~$\bG$, and in the subsection~\ref{part:estimates-3} we obtain the estimate for~$\partial_t \bG$.
This estimate requires an estimate for~$\partial_T \bG$, which is given in the subsection~\ref{part:estimates-2}.

\subsubsection{Estimate for the deformation gradient~$\bG$}\label{part:estimates-1}

Let $q>d$. We take the inner product of the equation~\eqref{39} by $q|\bG|^{q-2} \bG$, and next we integrate for $\bx\in \T$. Due to the incompressible condition $\div\, \overline \bv = 0$, we obtain
\begin{equation*}\label{46}
\begin{aligned}
\partial_t \|\bG\|_{L^q_{\bx}}^q + \frac{1}{\We} \partial_T \|\bG\|_{L^q_{\bx}}^q 
& = q \int_\T |\bG|^{q-2} (\bG\cdot \nabla \overline \bv) : \bG \\
& \leq q \|\nabla \overline \bv\|_{L^\infty_{\bx}} \|\bG\|_{L^q_{\bx}}^q.
\end{aligned}
\end{equation*}
Then we use the continuous injection $W^{1,q}_{\bx} \hookrightarrow L^\infty_{\bx}$, holds for $q>d$ and making appear a constant~$C_s$:
\begin{equation}\label{47}
\begin{aligned}
\partial_t \|\bG\|_{L^q_{\bx}}^q + \frac{1}{\We} \partial_T \|\bG\|_{L^q_{\bx}}^q
\leq
q\, C_s \|\nabla \overline \bv\|_{W^{1,q}_{\bx}} \|\bG\|_{L^q_{\bx}}^q.
\end{aligned}
\end{equation}
Now, we take the spatial gradient in~\eqref{39} and compute the inner product of both sides of the resulting equation with $q|\nabla \bG|^{q-2} \nabla \bG$ (we will note that this is a inner product on the $3$-tensor, defined by $\bA::\bB = a_{i,j,k}\,b_{i,j,k}$). After integrating for $\bx\in \T$ we obtain
\begin{equation*}\label{48}
\begin{aligned}
\partial_t \|\nabla \bG\|_{L^q_{\bx}}^q + \frac{1}{\We} \partial_T \|\nabla \bG\|_{L^q_{\bx}}^q
\leq & \; 2q \int_\T |\nabla \bG|^q |\nabla \overline \bv| + q \int_\T |\bG| |\nabla \bG|^{q-1} |\nabla^2\overline \bv|.
\end{aligned}
\end{equation*}
Using the Hölder inequality and the continuous injection $W^{1,q}_{\bx} \hookrightarrow L^\infty_{\bx}$ again, we deduce
\begin{equation}\label{49}
\partial_t  \|\nabla \bG\|_{L^q_{\bx}}^q + \frac{1}{\We} \partial_T \|\nabla \bG\|_{L^q_{\bx}}^q
\leq 3q\, C_s \|\nabla \overline \bv\|_{W^{1,q}_{\bx}} \|\bG\|_{W^{1,q}_{\bx}}^q.
\end{equation}
Adding this estimate~\eqref{49} with the estimate~\eqref{47}, we obtain
\begin{equation*}\label{50}
\partial_t \|\bG\|_{W^{1,q}_{\bx}} + \frac{1}{\We} \partial_T \|\bG\|_{W^{1,q}_{\bx}}
\leq 4C_s \|\nabla \overline \bv\|_{W^{1,q}_{\bx}} \|\bG\|_{W^{1,q}_{\bx}}.
\end{equation*}
Using the initial conditions we have
\begin{equation*}\label{51}
\|\bG\|_{W^{1,q}_{\bx}}\big|_{t=0} = \|\bG_0\|_{W^{1,q}_{\bx}}
\quad \text{and} \quad
\|\bG\|_{W^{1,q}_{\bx}}\big|_{T=0} = \sqrt{d},
\end{equation*}
so that the Gronwall type lemma~\ref{lem:gronwall} (see the preliminary subsection~\ref{subsection-gronwall}) implies that for $(t,T)\in (0,t^\star) \!\times \! \R^+$ we have
\begin{equation}\label{52}
\|\bG\|_{W^{1,q}_{\bx}}(t,T) \leq \zeta(t,T) \mathrm{exp} \Big( 4C_s \int_0^t \|\nabla \overline \bv\|_{W^{1,q}_{\bx}}\Big),
\end{equation}
where $\displaystyle \zeta(t,T) = \left\{
\begin{aligned}
\|\bG_0\|_{W^{1,q}_{\bx}}\big( T-\frac{t}{\We}\big) \quad & \text{if $t\leq \We\, T$}, \\
\sqrt{d} \hspace{1.7cm} & \text{if $t > \We\, T$}.
\end{aligned}
\right.
$
\par\noindent
The assumption $\bG_0\in L^\infty_T W^{1,q}_{\bx}$ implies $\zeta \in L^\infty(0,t^\star;L^\infty_T)$ with 
\begin{equation*}\label{53}
\|\zeta\|_{L^\infty(0,t^\star;L^\infty_T)}\leq \max\big\{\|\bG_0\|_{L^\infty_T W^{1,q}_{\bx}} , \sqrt{d} \big\}.
\end{equation*}
The relation~\eqref{52} now reads
\begin{equation}\label{54}
\begin{aligned}
& \|\bG\|_{L^\infty(0,t^\star;L^\infty_T W^{1,q}_{\bx})}
\leq
\|\zeta\|_{L^\infty(0,t^\star;L^\infty_T)} \, \mathrm{exp}\big( 4 C_s \|\nabla \overline \bv\|_{L^1(0,t^\star;W^{1,q}_{\bx})} \big).
\end{aligned}
\end{equation}

\subsubsection{Estimate for the age derivate~$\partial_T\bG$}\label{part:estimates-2}

We first remark that the derivative $\bG'=\partial_T \bG$ exactly satisfies the same PDE that $\bG$ (see the equation of~\eqref{39}; that is due to the fact that $\overline \bv$ does not depend on the variable~$T$).
We then deduce the same kind of estimate that~\eqref{47}:
\begin{equation*}\label{55}
\partial_t \|\bG'\|_{L^q_{\bx}} + \frac{1}{\We} \partial_T \|\bG'\|_{L^q_{\bx}}
\leq C_s \|\nabla \overline \bv\|_{W^{1,q}_{\bx}} \|\bG'\|_{L^q_{\bx}}.
\end{equation*}
But the initial conditions differ as follows: 
\begin{equation*}\label{56}
\bG'|_{t=0}=\partial_T\bG_0
\quad \text{and} \quad 
\bG'|_{T=0}=\We\, \nabla \overline \bv.
\end{equation*}
This last condition is obtained using $T=0$ in the equation~\eqref{39}.
Note that this result is valid because we are working on regular solutions $\bG_n$ (see the introduction of this proof) such that $\partial_t \bG_n$ is continuous at $T=0$.
From Lemma~\ref{lem:gronwall} given in the subsection~\ref{subsection-gronwall} we obtain for all $(t,T)\in (0,t^\star)\!\times \! \R^+$ the estimate
\begin{equation}\label{57}
\|\bG'\|_{L^q_{\bx}}(t,T)
\leq
\zeta'(t,T) \, \mathrm{exp} \Big( C_s \int_0^t \|\nabla \overline \bv\|_{W^{1,q}_{\bx}}\Big),
\end{equation}
where $ \displaystyle \zeta'(t,T)=\left\{
\begin{aligned}
\|\partial_T\bG_0\|_{L^q_{\bx}}\big( T-\frac{t}{\We} \big) \quad & \text{if $t\leq \We\, T$},\\
\We \, \|\nabla \overline \bv\|_{L^q_{\bx}}\big( t-\We\, T \big) \quad & \text{if $t > \We\, T$}.\\
\end{aligned}
\right.
$
\par\noindent
For each $t\in(0,t^\star)$ we estimate the $L^\infty_T$-norm of the function $T\mapsto \zeta'(t,T)$ as follows:
\begin{equation*}\label{58}
\| \zeta' \|_{L^\infty_T}(t) \leq \max( \We \, \| \nabla \overline \bv\|_{L^\infty(0,t;L^q_{\bx})} ,  \|\partial_T\bG_0\|_{L^\infty_T L^q_{\bx}} ).
\end{equation*}
Taking the norm in $L^r(0,t^\star)$ we deduce that $\zeta'\in L^r(0,t^\star; L^\infty_T)$ with
\begin{equation*}\label{59}
\|\zeta'\|_{L^r(0,t^\star; L^\infty_T)} \leq  \max( \We \, \| \nabla \overline \bv\|_{L^\infty(0,t^\star;L^q_{\bx})} ,  \|\partial_T\bG_0\|_{L^\infty_T L^q_{\bx}} ) \, {t^\star}^\frac{1}{r}.
\end{equation*}
By assumption we have $t^\star<1$ so that ${t^\star}^\frac{1}{r}<1$.
The relation~\eqref{57} now reads
\begin{equation}\label{60}
\|\bG'\|_{L^r(0,t^\star;L^\infty_T L^q_{\bx})}
\leq
\|\zeta'\|_{L^r(0,t^\star; L^\infty_T)} \, \mathrm{exp} \big( C_s \|\nabla \overline \bv\|_{L^1(0,t^\star;W^{1,p}_{\bx})} \big).
\end{equation}

\subsubsection{Estimate for the time derivate~$\partial_t\bG$}\label{part:estimates-3}

Isolating the term $\partial_t\bG$ in the equation~\eqref{39} we have
\begin{equation*}\label{61}
\begin{aligned}
\|\partial_t \bG\|_{L^q_{\bx}}
& \leq
\frac{1}{\We} \|\bG'\|_{L^q_{\bx}} + \|\overline \bv\|_{L^\infty_{\bx}} \|\nabla \bG\|_{L^q_{\bx}} + \|\bG\|_{L^\infty_{\bx}} \|\nabla \overline \bv\|_{L^q_{\bx}} \\
& \leq
\frac{1}{\We} \|\bG'\|_{L^q_{\bx}} + 2 C_s \|\overline \bv\|_{W^{1,q}_{\bx}} \|\bG\|_{W^{1,q}_{\bx}}.
\end{aligned}
\end{equation*}
Taking the $L^r(0,t^\star;L^\infty_T)$-norm for the variable~$(t,T)$, we obtain
\begin{equation*}\label{62}
\begin{aligned}
& \|\partial_t \bG\|_{L^r(0,t^\star;L^\infty_T L^q_{\bx})}
\leq \frac{1}{\We} \|\bG'\|_{L^r(0,t^\star;L^\infty_T L^q_{\bx})} \\
& \hspace{3.5cm} 
+ 2 C_s \|\overline \bv\|_{L^r(0,t^\star;W^{1,q}_{\bx})} \|\bG\|_{L^\infty(0,t^\star;L^\infty_T W^{1,q}_{\bx})}.
\end{aligned}
\end{equation*}
Using the previous estimates~\eqref{54} and~\eqref{60}, we deduce the result announced in Lemma~\ref{lem:1116}.
\cqfd

\subsection{Estimates for the memory function~$K$, solution of~\eqref{40}}

Unlike standard memory models\footnote{The term ``standard'' refers to models like K-BKZ in which the memory is a given function, usually on exponential type: $K(T)=\mathrm e^{-T}$.}, the estimate of the memory function in the Doi-Edwards model is one of the key point in the local existence proof.
In fact, the Doi-Edwards model is strongly based on the expression of the memory and on the equation that it satisfies.
We introduce the norm of the memory function that we will controlled in the next lemma:
\begin{equation}\label{63}
\begin{aligned}
\vvvert K \vvvert_3 := 
& \| K \|_{L^\infty(0,t^\star;L^\infty_T L^\infty_{\bx,s})} 
+ \| \partial_T K \|_{L^\infty(0,t^\star;L^\infty_T L^\infty_{\bx,s})} \\
& + \| \mathrm e^{\mu \, (\We \, T-t)} \partial_T K \|_{L^\infty(0,t^\star;L^\infty_T L^\infty_{\bx,s})}
+ \|\partial_T \nabla K\|_{L^\infty(0,t^\star;L^1_T L^q_{\bx,s}\cap L^2_T L^q_{\bx,s})} \\
& + \|\partial_t K \|_{L^2(0,t^\star;L^\infty_T L^2_{\bx,s})} + \|\partial_s K\|_{L^\infty(0,t^\star;L^\infty_T L^2_{\bx,s})}.
\end{aligned}
\end{equation}
\begin{lemma}\label{lem:1117}
Let $t^\star>0$ and $q\in]d,+\infty[$.
\par\noindent
If there exists $\mu>0$ such that $\mathrm e^{\We \mu T} \partial_T K_0 \in L^\infty_T L^\infty_{\bx,s}$ and $\mathrm e^{\We \mu T/2} \partial_T \nabla K_0 \in L^2_T L^q_{\bx,s}$,
\par\noindent
if $\partial_s K_0 \in L^\infty_T L^2_{\bx,s}$ satisfies $\partial_T K_0\leq 0$,
\par\noindent
if $\overline \bv \in L^2(0,t^\star;W^{1,\infty}_{\bx})$ is free divergence,
\par\noindent
if $\partial_s g \in L^2(0,t^\star;L^\infty_{\bx,s})$ and $\partial_s \nabla g \in L^2(0,t^\star;L^q_{\bx,s})$,
\par\noindent
then the problem~\eqref{40} admits a unique solution~$K$ satisfying $\partial_T K\leq 0$ and
\begin{equation*}\label{64}
\vvvert K \vvvert_3 \leq F_3\big( \|\partial_s g\|_{L^2(0,t^\star;L^\infty_{\bx,s})}, \|\partial_s \nabla g\|_{L^2(0,t^\star;L^q_{\bx,s})}, \| \overline \bv \|_{L^2(0,t^\star;W^{1,\infty}_{\bx})} \big),
\end{equation*}
where the function~$F_3$ depends on $q$, $\We$, $\mu$ and~$K_0$.
Moreover this function is continuous and nondecreasing in each of its variables.
\end{lemma}\\
\proof
The proof is composed of multiple steps since we need to control each term of the norm~$\vvvert K \vvvert_3$ introduced by~\eqref{63}.

\checkmark {\bf Step 1: control of $K$} -- Since $\partial_s g \in L^2(0,t^\star; L^\infty_{\bx,s})$ and we can directly apply the lemma~\ref{lem:max} given in the subsection~\ref{max-principle}. We deduce that for almost every $(t,T,\bx,s) \in (0,t^\star) \times \R^+ \times \T \times (-\frac{1}{2},\frac{1}{2})$ we have the estimate
\begin{equation}\label{65}
\min\{ 1, \inf K_0\} \leq K(t,T,\bx,s) \leq \max \{ 1,  \sup K_0\}.
\end{equation}
We remark that the assumption $\mathrm e^{\We \, \mu T} \partial_T K_0 \in L^\infty_T L^\infty_{\bx,s}$ implies a $L^\infty$-bound on~$K_0$:
\begin{equation}\label{66}
\begin{aligned}
\Big| \| K_0 \|_{L^\infty_{\bx,s}}(T) - 1 \Big| 
& \leq \Big| \int_0^T \| \partial_T K_0 \|_{L^\infty_{\bx,s}}(T') \, \mathrm dT' \Big| \\
& \leq \| \mathrm e^{\We \, \mu T} \partial_T K_0 \|_{L^\infty_T L^\infty_{\bx,s}} \int_0^T \mathrm e^{-\We \, \mu T'} \, \mathrm dT' \\
& \leq \frac{1}{\We \, \mu} \| \mathrm e^{\We \, \mu T} \partial_T K_0 \|_{L^\infty_T L^\infty_{\bx,s}}.
\end{aligned}
\end{equation}
To obtain the first inequality above we used the inequality $\partial_T \big( \| K_0 \|_{L^\infty_{\bx,s}} \big) \leq \| \partial_T K_0 \|_{L^\infty_{\bx,s}}$ which simply comes from to the inequality $\sup(b) - \sup(a) \leq \sup(b-a)$.
\par\noindent
As a consequence of~\eqref{66}, inequality~\eqref{65} effectively gives a $L^\infty$-bound on~$K$.\\

\checkmark {\bf Step 2: control of $\partial_TK$} -- Introducing $m=-\partial_T K$ and using the fact that~$g$ and~$\bv$ do not depend on the variable~$T$, we remark that Equation~\eqref{40} for~$K$ is translated into a similar equation for~$m$:
\begin{equation*}\label{67}
\mathrm d_t m + \frac{1}{\We} \partial_T m + g \, \partial_s m - \frac{1}{\We} \partial^2_s m = 0,
\end{equation*}
with the following boundary and initial conditions:
\begin{equation*}\label{68}
m\big|_{t=0} = -\partial_T K_0, \qquad
m\big|_{T=0} = 0, \qquad
m\big|_{s=\frac{1}{2}} = m\big|_{s=-\frac{1}{2}} = 0.
\end{equation*}
We can apply the lemma~\ref{lem:max} again to deduce that for almost every $(t,T,\bx,s) \in (0,t^\star) \times \R^+ \times \T \times (-\frac{1}{2},\frac{1}{2})$ we have
\begin{equation}\label{69}
0 \leq \inf (-\partial_T K_0) \leq m(t,T,\bx,s) \leq \sup (-\partial_T K_0).
\end{equation}
Obviously the $L^\infty$-bound on $\mathrm e^{\We \, \mu T} \partial_T K_0$ implies a $L^\infty$-bound on $\partial_T K_0$ so that the inequality~\eqref{69} corresponds to a $L^\infty$-bound on~$m=-\partial_T K$.
This estimate~\eqref{69} also proves that $m\geq 0$.\\

\checkmark {\bf Step 3: exponential control of $\partial_T K$} -- The exponential decreasing is obtained introducing $\widehat m = m \, \mathrm e^{\We \, \mu T} - C_0 \, \mathrm e^{\mu t}$ where $\displaystyle C_0 = \sup_{T,\bx,s}(-\partial_T K_0(T,\bx,s)\, \mathrm e^{\We \, \mu T})$. The quantity~$\widehat m$ satisfies
\begin{equation}\label{70}
\mathrm d_t \widehat m + \frac{1}{\We} \partial_T \widehat m - \mu\, \widehat m + g \partial_s \widehat m - \frac{1}{\We} \partial^2_s \widehat m = 0,
\end{equation}
with the following boundary and initial conditions (we note that we use the assumption $\partial_T K_0 \leq 0$):
\begin{equation*}\label{71}
\widehat m \big|_{t=0} \leq 0, \qquad
\widehat m \big|_{T=0} \leq 0, \qquad
\widehat m \big|_{s=-\frac{1}{2}} \leq 0, \qquad
\widehat m \big|_{s=\frac{1}{2}} \leq 0.
\end{equation*}
Testing the equation~\eqref{70} with the positive part of~$\widehat m$, we easily deduce that $\widehat m \leq 0$ (see the proof of Lemma~\ref{lem:max} for similar result).
This implies the desired bound for almost every $(t,T,\bx,s) \in (0,t^\star) \times \R^+ \times \T \times (-\frac{1}{2},\frac{1}{2})$:
\begin{equation}\label{72}
\mathrm e^{\mu \, (\We \, T-t)} m(t,T,\bx,s) \leq C_0.
\end{equation}

$\checkmark$ {\bf Step 4: control of $\partial_T \nabla K$} -- To have a bound on $\| \nabla m\|_{L^1_T L^q_{\bx,s}}$ we first note that if we introduce $\tm = \mathrm e^{\We \, \lambda T}m$, for any $\lambda>0$, then we have
\begin{equation}\label{73}
\begin{aligned}
\int_0^{+\infty} \| \nabla m\|_{L^q_{\bx,s}}(t,T) \, \mathrm dT
& =
\int_0^{+\infty} \mathrm e^{-\We \, \lambda T} \| \nabla \tm\|_{L^q_{\bx,s}}(t,T) \, \mathrm dT \\
& \leq
\Big( \int_0^{+\infty} \mathrm e^{-2\We \, \lambda T} \, \mathrm dT \Big)^{\frac{1}{2}} \Big( \int_0^{+\infty} \| \nabla \tm\|_{L^q_{\bx,s}}^2(t,T) \, \mathrm dT \Big)^{\frac{1}{2}} \\
& \leq \Big( \frac{1}{2 \, \We \, \lambda} \Big)^{\frac{1}{2}} \Big( \int_0^{+\infty} \| \nabla \tm\|_{L^q_{\bx,s}}^2(t,T) \, \mathrm dT \Big)^{\frac{1}{2}}.
\end{aligned}
\end{equation}
It is then sufficient to control~$\| \nabla \tm\|_{L^2_T L^q_{\bx,s}}$.
The quantity~$\tm$ satisfies
\begin{equation}\label{74}
\mathrm d_t \tm + \frac{1}{\We} \partial_T \tm  - \lambda \tm + g \partial_s \tm - \frac{1}{\We} \partial^2_s \tm = 0,
\end{equation}
with the following boundary and initial conditions:
\begin{equation*}\label{75}
\tm\big|_{t=0} = -\mathrm e^{\We \, \lambda T}\partial_T K_0, \qquad
\tm\big|_{T=0} = 0, \qquad
\tm\big|_{s=-\frac{1}{2}} = \tm\big|_{s=\frac{1}{2}} = 0.
\end{equation*}
Derivating the equation~\eqref{74} with respect to the spatial variable~$\bx$ we obtain\begin{equation}\label{76}
\mathrm d_t \nabla \tm + \nabla \overline \bv \cdot \nabla \tm + \frac{1}{\We} \partial_T \nabla \tm - \lambda \nabla \tm + \partial_s \tm \nabla g + g \, \partial_s \nabla \tm - \frac{1}{\We} \partial^2_s \nabla \tm = 0.
\end{equation}
We take the inner product of this equation~\eqref{76} by $q|\nabla \tm|^{q-2}\nabla \tm$, and we integrate with respect to $\bx\in \T$ and $s\in (-\frac{1}{2},\frac{1}{2})$. We deduce
\begin{equation*}\label{77}
\begin{aligned}
& \partial_t \|\nabla \tm\|_{L^q_{\bx,s}}^q + \frac{1}{\We} \partial_T \|\nabla \tm\|_{L^q_{\bx,s}}^q - q \, \lambda \|\nabla \tm\|_{L^q_{\bx,s}}^q \\
& + \frac{q}{\We} \int_{\T} \int_{-\frac{1}{2}}^{\frac{1}{2}} |\nabla \tm|^{q-2} |\partial_s \nabla \tm|^2 + \frac{q(q-2)}{\We} \int_{\T} \int_{-\frac{1}{2}}^{\frac{1}{2}} |\nabla \tm \cdot \partial_s \nabla \tm|^2 \, |\nabla \tm|^{q-4}\\
& \hspace{7cm} \leq A_0 + A_1 + A_2,
\end{aligned}
\end{equation*}
where
\begin{equation*}\label{78}
\begin{aligned}
& A_0 = - q \int_{\T} \int_{-\frac{1}{2}}^{\frac{1}{2}} (\nabla \overline \bv \cdot \nabla \tm) \cdot \nabla \tm \, |\nabla \tm|^{q-2}, \\
& A_1 = - q \int_{\T} \int_{-\frac{1}{2}}^{\frac{1}{2}} \Big( \nabla g \cdot \nabla \tm \Big) |\nabla \tm|^{q-2} \partial_s \tm, \\
& A_2 = - q \int_{\T} \int_{-\frac{1}{2}}^{\frac{1}{2}} \Big( \partial_s \nabla \tm \cdot \nabla \tm \Big) |\nabla \tm|^{q-2} g.
\end{aligned}
\end{equation*}
\checkmark
The first term $A_0$ is directly estimate as follows
\begin{equation}\label{79}
|A_0| \leq q \|\nabla \overline \bv\|_{L^\infty_{\bx}} \| \nabla \tm \|_{L^q_{\bx,s}}^q.
\end{equation}
\checkmark
The term $A_1$ is estimate as follows: we integrate by parts with respect to the variable~$s$ in order to write $A_1 = A_{11} + A_{12} + A_{13}$ with 
\begin{equation*}\label{80}
\begin{aligned}
& A_{11} = -q \int_{\T} \int_{-\frac{1}{2}}^{\frac{1}{2}} \Big( \partial_s \nabla g \cdot \nabla \tm \Big) |\nabla \tm|^{q-2} \tm, \\
& A_{12} = -q \int_{\T} \int_{-\frac{1}{2}}^{\frac{1}{2}} \Big( \nabla g \cdot \partial_s \nabla \tm \Big) |\nabla \tm|^{q-2} \tm, \\
& A_{13} = -q(q-2) \int_{\T} \int_{-\frac{1}{2}}^{\frac{1}{2}} \Big( \nabla g \cdot \nabla \tm \Big) \Big( \partial_s \nabla \tm \cdot \nabla \tm \Big) |\nabla \tm|^{q-4} \tm.
\end{aligned}
\end{equation*}
To estimate the term $A_{11}$ (which only depends on~$t$ and~$T$) we introduce the $L^\infty_{\bx,s}$ bound on~$\tm$:
\begin{equation*}\label{81}
|A_{11}| \leq q \| \tm \|_{L^\infty_{\bx,s}} \int_{\T} \int_{-\frac{1}{2}}^{\frac{1}{2}} |\partial_s \nabla g| \, |\nabla \tm|^{q-1}.
\end{equation*}
Due to the Hölder inequality we obtain
\begin{equation}\label{82}
|A_{11}| \leq q \| \tm \|_{L^\infty_{\bx,s}} \| \partial_s \nabla g \|_{L^q_{\bx,s}} \| \nabla \tm \|_{L^q_{\bx,s}}^{q-1}.
\end{equation}
In the same way, the term $A_{12}$ is treated as follows:
\begin{equation*}\label{83}
|A_{12}| \leq q \| \tm \|_{L^\infty_{\bx,s}} \int_{\T} \int_{-\frac{1}{2}}^{\frac{1}{2}}  |\partial_s \nabla \tm| \, |\nabla g| \, |\nabla \tm|^{q-2}.
\end{equation*}
Using the Young inequality and next the Hölder inequality, we write
\begin{equation}\label{84}
\begin{aligned}
|A_{12}| & \leq 
\frac{q}{2 \, \We} \int_{\T} \int_{-\frac{1}{2}}^{\frac{1}{2}} |\nabla \tm|^{q-2} |\partial_s \nabla \tm|^2
+ \frac{q\, \We}{2} \| \tm \|_{L^\infty_{\bx,s}}^2 \int_{\T} \int_{-\frac{1}{2}}^{\frac{1}{2}} | \nabla g |^2 \, | \nabla \tm |^{q-2} \\
& \leq \frac{q}{2 \, \We} \int_{\T} \int_{-\frac{1}{2}}^{\frac{1}{2}} |\nabla \tm|^{q-2} |\partial_s \nabla \tm|^2
+ \frac{q\, \We}{2} \| \tm \|_{L^\infty_{\bx,s}}^2 \| \nabla g \|_{L^q_{\bx,s}}^2 \| \nabla \tm \|_{L^q_{\bx,s}}^{q-2}.
\end{aligned}
\end{equation}
Similarly, using the Young inequality and the Hölder inequality, the term $A_{13}$ is controlled as follows
\begin{equation}\label{85}
\begin{aligned}
|A_{13}| & \leq \frac{q(q-2)}{2 \, \We} \int_{\T} \int_{-\frac{1}{2}}^{\frac{1}{2}} |\nabla \tm \cdot \partial_s \nabla \tm|^2 \, |\nabla \tm|^{q-4} \\
& \hspace{3cm} + \frac{q(q-2)\, \We}{2} \| \tm \|_{L^\infty_{\bx,s}}^2 \| \nabla g \|_{L^q_{\bx,s}}^2 \| \nabla \tm \|_{L^q_{\bx,s}}^{q-2}.
\end{aligned}
\end{equation}
\checkmark
Using an integration by parts with respect to the variable~$s$, the contribution~$A_2$ simply reads
\begin{equation*}\label{86}
A_2 = - \int_{\T} \int_{-\frac{1}{2}}^{\frac{1}{2}} |\nabla \tm|^q \, \partial_s g.
\end{equation*}
Hence we have the estimate
\begin{equation}\label{87}
|A_2| \leq \| \partial_s g\|_{L^\infty_{\bx,s}} \, \|\nabla \tm \|_{L^q_{\bx,s}}^q .
\end{equation}
\checkmark
All these estimates~\eqref{79}, \eqref{82}, \eqref{84}, \eqref{85} and~\eqref{87} imply
\begin{equation*}\label{88}
\begin{aligned}
\partial_t \|\nabla \tm\|_{L^q_{\bx,s}}^q + \frac{1}{\We}\partial_T \|\nabla \tm\|_{L^q_{\bx,s}}^q 
\leq \, & 
q \| \tm \|_{L^\infty_{\bx,s}} \| \partial_s \nabla g \|_{L^q_{\bx,s}} \|\nabla \tm\|_{L^q_{\bx,s}}^{q-1} \\
& + \frac{q(q-1)\, \We}{2} \| \tm \|_{L^\infty_{\bx,s}}^2 \| \nabla g \|_{L^q_{\bx,s}}^2 \|\nabla \tm\|_{L^q_{\bx,s}}^{q-2} \\
& + \big( \| \partial_s g\|_{L^\infty_{\bx,s}} + q \, \lambda + q \|\nabla \overline \bv\|_{L^\infty_{\bx}}\big) \|\nabla \tm\|_{L^q_{\bx,s}}^q.
\end{aligned}
\end{equation*}
Multiplying by $\frac{2}{q}\|\nabla \tm\|_{L^q_{\bx,s}}^{2-q}$ we also deduce
\begin{equation*}\label{89}
\begin{aligned}
\partial_t \|\nabla \tm\|_{L^q_{\bx,s}}^2 + \frac{1}{\We}\partial_T \|\nabla \tm\|_{L^q_{\bx,s}}^2
\! \leq \, &
2 \| \tm \|_{L^\infty_{\bx,s}} \| \partial_s \nabla g \|_{L^q_{\bx,s}} \|\nabla \tm\|_{L^q_{\bx,s}} \\
& + (q-1)\, \We \| \tm \|_{L^\infty_{\bx,s}}^2 \| \nabla g \|_{L^q_{\bx,s}}^2 \\
& + \frac{2}{q} \big( \| \partial_s g\|_{L^\infty_{\bx,s}} + q \, \lambda + q \|\nabla \overline \bv\|_{L^\infty_{\bx}}\big) \|\nabla \tm\|_{L^q_{\bx,s}}^2.
\end{aligned}
\end{equation*}
We integrate with respect to $T\in (0,+\infty)$.
Since $\tm|_{T=0}=0$ the contribution due to the term $\partial_T \|\nabla \tm\|_{L^q_{\bx,s}}^2$ is non negative. We obtain
\begin{equation*}\label{90}
\begin{aligned}
\partial_t \|\nabla \tm\|_{L^2_T L^q_{\bx,s}}^2 \leq \, &
2 \| \partial_s \nabla g \|_{L^q_{\bx,s}} \int_0^{+\infty} \| \tm \|_{L^\infty_{\bx,s}}  \|\nabla \tm\|_{L^q_{\bx,s}}  \, \mathrm dT \\
& + (q-1)\, \We \| \nabla g \|_{L^q_{\bx,s}}^2 \int_0^{+\infty} \| \tm \|_{L^\infty_{\bx,s}}^2 \, \mathrm dT \\
& + \frac{2}{q} \big( \| \partial_s g\|_{L^\infty_{\bx,s}} + q \, \lambda + q \|\nabla \overline \bv\|_{L^\infty_{\bx}} \big) \int_0^{+\infty}  \|\nabla \tm\|_{L^q_{\bx,s}}^2 \, \mathrm dT.
\end{aligned}
\end{equation*}
We now use the bound $\| \tm \|_{L^\infty_{\bx,s}}(t,T)\leq C_0 \, \mathrm e^{\mu t} \mathrm e^{-\We (\mu-\lambda)\, T}$ obtained in the third step of the present proof, to deduce
\begin{equation*}\label{91}
\begin{aligned}
\partial_t \|\nabla \tm\|_{L^2_T L^q_{\bx,s}}^2 \leq \, & 
2 C_0 \, \mathrm e^{\mu t} \| \partial_s \nabla g \|_{L^q_{\bx,s}} \int_0^{+\infty} \mathrm e^{- \We (\mu-\lambda) T} \|\nabla \tm\|_{L^q_{\bx,s}} \, \mathrm dT \\
& + C_0^2 \, (q-1)\, \We \, \mathrm e^{2 \, \mu t} \| \nabla g \|_{L^q_{\bx,s}}^2 \int_0^{+\infty} \mathrm e^{-2 \We (\mu-\lambda) T} \, \mathrm dT \\
& + \frac{2}{q} \big( \| \partial_s g\|_{L^\infty_{\bx,s}} + q \, \lambda + q \|\nabla \overline \bv\|_{L^\infty_{\bx}} \big) \int_0^{+\infty}  \|\nabla \tm\|_{L^q_{\bx,s}}^2 \, \mathrm dT.
\end{aligned}
\end{equation*}
Using the Cauchy-Schwarz inequality, we obtain for $0<\lambda<\mu$ the inequality
\begin{equation*}\label{92}
\begin{aligned}
\partial_t \|\nabla \tm\|_{L^2_T L^q_{\bx,s}}^2 \leq \, & 
2 C_0 \, \mathrm e^{\mu t} \| \partial_s \nabla g \|_{L^q_{\bx,s}} \Big( \frac{1}{2\, \We (\mu-\lambda)}\Big)^{\frac{1}{2}} \|\nabla \tm\|_{L^2_T L^q_{\bx,s}}\\
& + C_0^2 (q-1)\, \We \, \mathrm e^{2\mu t} \| \nabla g \|_{L^q_{\bx,s}}^2 \Big( \frac{1}{2\, \We (\mu-\lambda)}\Big) \\
& + \frac{2}{q} \big( \| \partial_s g\|_{L^\infty_{\bx,s}} + q \, \lambda + q \|\nabla \overline \bv\|_{L^\infty_{\bx}} \big) \|\nabla \tm\|_{L^2_T L^q_{\bx,s}}^2.
\end{aligned}
\end{equation*}
Using the Young inequality we deduce
\begin{equation}\label{93}
\begin{aligned}
\partial_t \|\nabla \tm\|_{L^2_T L^q_{\bx,s}}^2 \leq \, &
\frac{C_0^2}{\We (\mu-\lambda)} \, \mathrm e^{2 \mu t} \| \partial_s \nabla g \|_{L^q_{\bx,s}}^2 + \frac{1}{2} \|\nabla \tm\|_{L^2_T L^q_{\bx,s}}^2\\
& + \frac{C_0^2 \, (q-1)}{2\, (\mu-\lambda)} \mathrm e^{2\mu t} \| \nabla g \|_{L^q_{\bx,s}}^2 \\
& + \frac{2}{q} \big( \| \partial_s g\|_{L^\infty_{\bx,s}} + q \, \lambda + q \|\nabla \overline \bv\|_{L^\infty_{\bx}} \big) \|\nabla \tm\|_{L^2_T L^q_{\bx,s}}^2.
\end{aligned}
\end{equation}
We note that $\|\nabla g\|_{L^q_{\bx,s}} \leq \|\partial_s \nabla g\|_{L^q_{\bx,s}}$. Indeed, since $g|_{s=0}=0$ we have
\begin{equation*}\label{94}
\nabla g(t,\bx,s) = \int_0^s \partial_s \nabla g(t,\bx,s')\, \mathrm ds'.
\end{equation*}
By the triangular inequality and next by the Hölder inequality we conclude that
\begin{equation*}\label{95}
\|\nabla g\|_{L^q_{\bx,s}} \leq \int_{-\frac{1}{2}}^{\frac{1}{2}} \|\partial_s \nabla g\|_{L^q_{\bx}}(t,s')\, \mathrm ds' \leq \|\partial_s\nabla g\|_{L^q_{\bx,s}}.
\end{equation*}
By assumption we know that $\|\nabla \tm\|_{L^2_T L^q_{\bx,s}}(0) = \|\mathrm e^{\We \, \lambda \, T}\partial_T \nabla K_0\|_{L^2_T L^q_{\bx,s}}$ is bounded for $0<\lambda\leq \frac{\mu}{2}$. We then now choose $\lambda = \frac{\mu}{2}$ so that equation~\eqref{93} becomes
\begin{equation*}\label{97}
\partial_t \|\nabla \tm\|_{L^2_T L^q_{\bx,s}}^2 \leq a \, \|\nabla \tm\|_{L^2_T L^q_{\bx,s}}^2 + b,
\end{equation*}
with
\begin{equation*}\label{98}
\begin{aligned}
& a(t) = \frac{1}{2} + \frac{2}{q} \| \partial_s g\|_{L^\infty_{\bx,s}} + \mu + 2 \|\nabla \overline \bv\|_{L^\infty_{\bx}} \\
\text{and} \qquad
& b(t)=\frac{C_0^2}{\mu} \Big(\frac{2}{\We} + q - 1 \Big) \mathrm e^{2 \mu t} \| \partial_s \nabla g \|_{L^q_{\bx,s}}^2.
\end{aligned}
\end{equation*}
We conclude using the Gronwall lemma and the regularity which is assumed for~$g$ and~$\overline \bv$:
\begin{equation}\label{99}
\|\nabla \tm\|_{L^2_T L^q_{\bx,s}}^2 \leq C_0' \mathrm{exp}\Big( \int_0^t b \Big) + \int_0^t a(t')\, \mathrm{exp}\Big( \int_{t'}^t b \Big), 
\end{equation}
with
\begin{equation*}\label{100}
C_0' = \|\mathrm e^{\We \, \mu \, T/2}\partial_T \nabla K_0\|_{L^2_T L^q_{\bx,s}}.
\end{equation*}
The estimate~\eqref{99} shows that $\nabla \tm \in L^\infty(0,t^\star;L^2_T L^q_{\bx,s})$. From the inequality~\eqref{73}, it implies a bound on~$\nabla m$ in $L^\infty(0,t^\star;L^1_T L^q_{\bx,s})$.
We also remark that we could have directly obtained a bound on~$\nabla m$ in $L^\infty(0,t^\star;L^2_T L^q_{\bx,s})$ taking $\lambda=0$ starting form~\eqref{74}.

$\checkmark$ {\bf Step 5: control of $\partial_tK$ and $\partial_sK$} --  To get the bounds on~$\partial_tK$ and~$\partial_s K$ we use~$\partial_t K$ as test function in Equation~\eqref{40} satisfied by~$K$. We obtain
\begin{equation*}\label{101}
\begin{aligned}
\|\partial_t K\|_{L^2_{\bx,s}}^2 + \frac{1}{2\We} \partial_t \|\partial_s K\|_{L^2_{\bx,s}}^2 
& = - \int_\T\int_{-\frac{1}{2}}^{\frac{1}{2}} (\overline \bv \cdot \nabla K) \partial_t K \\
& \hspace{-2cm} - \int_\T\int_{-\frac{1}{2}}^{\frac{1}{2}} g \partial_s K \partial_t K
- \frac{1}{\We} \int_\T\int_{-\frac{1}{2}}^{\frac{1}{2}} \partial_T K \partial_t K \\
& \leq \| \overline \bv \|_{L^\frac{q-2}{2q}_{\bx}} \| \nabla K \|_{L^q_{\bx,s}} \| \partial_t K \|_{L^2_{\bx,s}} \\
& \hspace{-2cm} + \|g\|_{L^\infty_{\bx,s}} \|\partial_s K\|_{L^2_{\bx,s}} \|\partial_t K\|_{L^2_{\bx,s}}
+ \frac{1}{\We} \|\partial_T K\|_{L^2_{\bx,s}} \|\partial_t K\|_{L^2_{\bx,s}}.
\end{aligned}
\end{equation*}
Using the fact that $L^\infty_{\bx} \subset L^\frac{q-2}{2q}_{\bx}$ and using the Young inequality we deduce
\begin{equation*}\label{102}
\begin{aligned}
\|\partial_t K\|_{L^2_{\bx,s}}^2 + \frac{1}{\We} \partial_t \|\partial_s K\|_{L^2_{\bx,s}}^2 
& \leq 3 \|\overline \bv\|_{L^\infty_{\bx}}^2 \| \nabla K \|_{L^q_{\bx,s}}^2 \\
& \hspace{-1cm} + 3 \|g\|_{L^\infty_{\bx,s}}^2 \|\partial_s K\|_{L^2_{\bx,s}}^2
+ \frac{3}{\We^2} \|\partial_T K\|_{L^2_{\bx,s}}^2.
\end{aligned}
\end{equation*}
Taking the supremum with respect to the variable~$T\in (0,+\infty)$ we obtain\footnote{Note that we use the inequality $\partial_t \big( \sup_T F(t,T) \big) \leq \sup_T \big( \partial_t F(t,T) \big)$ which simply comes from to the inequality $\sup(b) - \sup(a) \leq \sup(b-a)$.}
\begin{equation}\label{103}
\begin{aligned}
\|\partial_t K\|_{L^\infty_T L^2_{\bx,s}}^2 + \frac{1}{\We} \partial_t \|\partial_s K\|_{L^\infty_T L^2_{\bx,s}}^2 
& \leq 3 \|\overline \bv\|_{L^\infty_{\bx}}^2 \|\nabla K\|_{L^\infty_T L^q_{\bx,s}}^2 \\
& \hspace{-1cm} + 3 \|g\|_{L^\infty_{\bx,s}}^2 \|\partial_s K\|_{L^\infty_T L^2_{\bx,s}}^2
+ \frac{3}{\We^2} \|\partial_T K\|_{L^\infty_T L^2_{\bx,s}}^2.
\end{aligned}
\end{equation}
By assumptions, we have $\|\overline \bv\|_{L^\infty_{\bx}}^2 \in L^1(0,t^\star)$ and $\|g\|_{L^\infty_{\bx,s}}^2 \in L^1(0,t^\star)$ (note that $g|_{s=0}=0$ so that $\|g\|_{L^\infty_{\bx,s}} \leq \|\partial_s g\|_{L^\infty_{\bx,s}}$).
Using the steps~2 and~4 respectively\footnote{Due to the condition $\nabla K\big|_{T=0} = 0$ the step~$4$ implies that $\|\nabla K\|_{L^q_{\bx,s}} \in L^\infty(0,t^\star;W^{1,1}_T)$.}, we know that $\|\partial_T K\|_{L^\infty_T L^2_{\bx,s}}^2$ and $\|\nabla K\|_{L^\infty_T L^q_{\bx,s}}^2$ belong in $L^\infty(0,t^\star)$. Consequently we can apply the Gronwall lemma to deduce from the inequality~\eqref{103} the following bound:
\begin{equation*}\label{104}
\|\partial_t K\|_{L^2(0,t^\star;L^\infty_T L^2_{\bx,s})}^2 + \|\partial_s K\|_{L^\infty(0,t^\star;L^\infty_T L^2_{\bx,s})}^2 \leq C,
\end{equation*}
where $C$ depends on the previous bounds.
\cqfd

\subsection{Proof of Theorem~\ref{th:local}}\label{part:estimates-4}

For any $t^\star>0$ we introduce the Banach space 
\begin{equation*}\label{105}
\begin{aligned}
\mathscr{B}(t^\star) = 
& L^r\big(0,t^\star;W^{1,q}(\T)\big) \\
& \times \mathcal C\big([0,t^\star] \! \times \! R^+ ;L^q(\T)\big) \\
& \quad  \times \mathcal C\big([0,t^\star]\! \times \! R^+ ;L^q(\T\! \times \! (-\tfrac{1}{2},\tfrac{1}{2}))\big)
\end{aligned}
\end{equation*}
and for any $R^\star>0$ the subset
\begin{equation*}\label{106}
\begin{aligned}
& \hspace{-1cm} \mathscr{H}(t^\star,R^\star)
=
\Big\{ (\overline{\bv} , \overline \bG, \overline{K}) \in \mathscr{B}(t^\star) ~ ; \\
& \overline \bv \in L^r(0,t^\star;\mathscr D(A_q)),
\quad \partial_t \overline \bv \in L^r(0,t^\star;H_q), \\
& \overline \bG \in L^\infty(0,t^\star,L^\infty_T W^{1,q}_{\bx}),
\quad  \partial_t \overline \bG, ~\partial_T \overline \bG \in L^r(0,t^\star;L^\infty_T L^q_{\bx}),\\
& \overline K, \partial_T \overline K \in L^\infty(0,t^\star;L^\infty_T L^\infty_{\bx,s}),
\quad  \partial_t \overline K \in L^2(0,t^\star;L^\infty_T L^q_{\bx,s}),\\
& \partial_T \nabla \overline K \in L^\infty(0,t^\star;L^1_T L^q_{\bx,s}\cap L^2_T L^q_{\bx,s}),
\quad  \partial_s \overline K \in L^\infty(0,t^\star;L^\infty_T L^q_{\bx,s}),\\
& \overline \bv|_{t=0} = \bv_0,
\quad \overline \bG|_{t=0} = \bG_0,
\quad \overline \bG|_{T=0} = \bdelta, \\
& \overline K|_{t=0} = K_0,
\quad \overline K|_{T=0} = 1,
\quad \overline K|_{s=-\frac{1}{2}} = \overline K|_{s=\frac{1}{2}} = 0, \\
& \vvvert \overline \bv \vvvert_1 \leq R^\star, \qquad
\vvvert \overline \bG \vvvert_2 \leq R^\star, \qquad
\vvvert \overline K \vvvert_3 \leq R^\star ~ \Big\},
\end{aligned}
\end{equation*}
where we recall that the norms~$\vvvert\cdot\vvvert_1$, $\vvvert\cdot\vvvert_2$ and~$\vvvert\cdot\vvvert_3$ are defined by~\eqref{42}, \eqref{45} and~\eqref{63} respectively.
\par\noindent
It is important to remark that such a set is non-empty, for instance if $R^\star$ is large enough. More precisely, if
\begin{equation}\label{107}
R^\star \geq \max ( F_1(0), \vvvert \bG_0 \vvvert_2,  \vvvert K_0 \vvvert_3 )
\end{equation}
then for any $t^\star>0$ we can build a velocity field $\bv^\star$ such that $(\bv^\star,\bG_0,K_0) \in \mathscr{H}(t^\star,R^\star)$, see an example of construction in~\cite[p.6]{Fernandez-Guillen-Ortega}.\\[-0.2cm]
\begin{remark}\label{rem:0944}
If $(\bv,\bG,K)\in \mathscr{H}(t^\star,R^\star)$ for some $t^\star$ and~$R^\star$ then the velocity field~$\bv$, the tensor~$\bG$ and the function~$K$ are continuous with respect to the time~$t$, the age~$T$ and the length~$s$.
In fact, these continuity properties follow from the Sobolev injections of kind $W^{1,\alpha}(0,t^\star;X) \subset \mathcal C([0,t^\star];X)$, hold for $\alpha>1$.
Moreover, they make sense of the initial conditions $\bv|_{t=0} = \bv_0$, $\bG|_{t=0} = \bG_0$ and~$\bG|_{T=0} = \bdelta$.
\par\noindent
More precisely, if $(\bv,\bG,K)\in \mathscr{H}(t^\star,R^\star)$ then we have
\begin{equation*}\label{107.5}
\begin{aligned}
& \bv|_{t=0} = \bv_0 && \text{in} \quad H_q,\\
& \bG|_{t=0} = \bG_0 && \text{in} \quad L^\infty_T L^q_{\bx},\\
& K|_{t=0} = K_0 && \text{in} \quad L^\infty_T L^q_{\bx,s}.
\end{aligned}
\end{equation*}
Noting\footnote{This result is  based on the following inequality $\sup_T \big( \int_0^{t^\star} f(t,T) \, \mathrm dt\big) \leq \int_0^{t^\star} \sup_T f(t,T) \, \mathrm dt$ for all positive function~$f\in L^1(0,t^\star;L^\infty_T)$.} that $L^r(0,t^\star;L^\infty_T L^q_{\bx}) \subset L^\infty(\R^+;L^r(0,t^\star;L^q_{\bx}))$ we also deduce that
\begin{equation*}\label{107.6}
\begin{aligned}
& \bG|_{T=0} = \bdelta && \text{in} \quad L^r(0,t^\star;L^q_{\bx}),\\
& K|_{T=0} = 1 && \text{in} \quad L^\infty(0,t^\star;L^\infty_{\bx,s}).
\end{aligned}
\end{equation*}
Similarly, we have for almost every $(t,T)\in (0,t^\star)\times\R^+$ the relation $K\in L^\infty_{\bx,s}$ and $\partial_s K \in L^q_{\bx,s}$. We deduce that
\begin{equation*}\label{107.7}
K|_{s=-\frac{1}{2}} = K|_{s=\frac{1}{2}} = 0 \quad \text{in} \quad L^\infty(0,t^\star;L^\infty_TL^q_{\bx}).
\end{equation*}
\end{remark}

\noindent
We consider the mapping
\begin{equation*}\label{108}
\begin{aligned}
\Phi ~:~ \mathscr{H}(t^\star,R^\star) & \longrightarrow ~ \mathscr{B}(t^\star) \\
(\overline{\bv} , \overline{\bG} , \overline{K}) ~~~~~ & \longmapsto ~(\bv,\bG,K),
\end{aligned}
\end{equation*}
where $\bv$ is the unique solution of the Stokes problem~\eqref{36} with
\begin{equation}\label{109}
\bbf = -\Re\, \overline{\bv} \cdot \nabla \overline{\bv} + \div\, \bsigma,
\end{equation}
where $\bG$ solves the problem~\eqref{39} depending on $\overline \bv$, and where~$K$ is given as the solution of~\eqref{40} with
\begin{equation}\label{110}
\partial_s g = \nabla \overline \bv : \bS.
\end{equation}
\begin{lemma}\label{lem:1648}
If $(\overline \bv,\overline \bG,\overline K)\in \mathscr{H}(t^\star,R^\star)$ and $r>2$ then we have
\begin{equation}\label{111}
\begin{aligned}
& \|\bbf\|_{L^r(0,t^\star;L^q_{\bx})} 
+ \|\overline \bv\|_{L^\infty(0,t^\star;W^{1,q}_{\bx})\cap L^1(0,t^\star;W^{2,q}_{\bx}) \cap L^2(0,t^\star;W^{1,\infty}_{\bx})} \\
& \hspace{2cm} \|\partial_s g\|_{L^2(0,t^\star;L^\infty_{\bx,s})} + \|\partial_s \nabla g\|_{L^2(0,t^\star;L^q_{\bx,s})} \leq G(t^\star, R^\star),
\end{aligned}
\end{equation}
where~$G$ is a continuous function vanishing for $t^\star = 0$.
\end{lemma}\\
\proof
We will estimate each term of the left hand side of~\eqref{111}.
\par
$\checkmark$ {\bf Step 1: Control of the source term~$\bbf$} --
The source term~$\bbf$ being given by the relation~\eqref{109}, we directly have
\begin{equation*}\label{112}
\begin{aligned}
& \|\bbf\|_{L^r(0,t^\star;L^q_{\bx})} \leq \Re \|\overline \bv \cdot \nabla \overline \bv\|_{L^r(0,t^\star;L^q_{\bx})} + \|\div \, \bsigma\|_{L^r(0,t^\star;L^q_{\bx})}.
\end{aligned}
\end{equation*}
We follow the ideas of~\cite{Fernandez-Guillen-Ortega}, generalized by~\cite{Chupin3} to the $d$-dimensional case, in order to treat the bilinear term $\overline \bv \cdot \nabla \overline \bv$. We have
\begin{equation}\label{113}
\|\overline \bv \cdot \nabla \overline \bv\|_{L^r(0,t^\star;L^q_{\bx})} \leq C\, {t^\star}^{\frac{q-d}{2rq}} {R^\star}^{\frac{q+d}{2q}} \|\bv_0\|_{L^p_{\bx}}^{\frac{3q-d}{2q}}
+ C\, {t^\star}^{\frac{3q-d}{2q}-\frac{1}{r}} {R^\star}^2.
\end{equation}
The last term $\div \, \bsigma$ is controlled using the orientation tensor~$\bS$, see~\eqref{37}:
\begin{equation*}\label{114}
\| \div \, \bsigma\|_{L^r(0,t^\star;L^q_{\bx})} \leq \omega \|\nabla \bS\|_{L^r(0,t^\star;L^q_{\bx,s})}.
\end{equation*}
By the definition of~$\bS$, see~\eqref{38}, the gradient $\nabla \bS$ reads (making appears $\overline m = -\partial_T \overline K$)
\begin{equation*}\label{115}
\begin{aligned}
\nabla \bS(t,\bx,s)
& = \int_0^{+\infty} \nabla \overline m(t,T,\bx,s) \otimes \mathscr S( \overline \bG(t,T,\bx))\, \mathrm dT
\\
& \qquad +
\int_0^{+\infty} \overline m(t,T,\bx,s) \, \mathscr S'( \overline \bG(t,T,\bx)) : \nabla \overline \bG(t,T,\bx) \, \mathrm dT.
\end{aligned}
\end{equation*}
Its $L^q_{\bx,s}$-norm can be estimate as follows
\begin{equation}\label{116}
\begin{aligned}
\|\nabla \bS\|_{L^q_{\bx,s}}(t) 
& \leq \mathscr S_\infty \int_0^{+\infty} \| \nabla \overline m\|_{L^q_{\bx,s}}(t,T) \, \mathrm dT \\
& \qquad + \widetilde{\mathscr S_\infty'} \int_0^{+\infty} \| \overline m \|_{L^\infty_{\bx,s}}(t,T) \| \nabla \overline \bG \|_{L^q_{\bx}}(t,T) \, \mathrm dT.
\end{aligned}
\end{equation}
We can read this inequality~\eqref{116} as
\begin{equation*}\label{117}
\|\nabla \bS\|_{L^q_{\bx,s}} \leq \mathscr S_\infty \| \nabla \overline m\|_{L^1_T L^q_{\bx,s}} + \mathscr S_\infty' \| \overline m \|_{L^1_T L^\infty_{\bx,s}} \| \nabla \overline \bG \|_{L^\infty_T L^q_{\bx}}.
\end{equation*}
Note that since $\| \mathrm e^{\mu(\We \, T-t)} \overline m \|_{L^\infty(0,t^\star;L^\infty_T L^\infty_{\bx,s})} \leq R^\star$ we have, for all $t\in (0,t^\star)$:
\begin{equation*}\label{118}
\begin{aligned}
\| \overline m \|_{L^1_T L^\infty_{\bx,s}}(t)
& = \int_0^\infty \| \overline m \|_{L^\infty_{\bx,s}}(t,T)\, \mathrm dT \\
& \leq R^\star \int_0^\infty \mathrm e^{-\mu(\We \, T-t)} \, \mathrm dT \\
& \leq \frac{R^\star \, \mathrm e^{\mu \, t}}{\mu \, \We}.
\end{aligned}
\end{equation*}
We deduce that for all $t\in (0,t^\star)$ we have
\begin{equation}\label{119}
\|\nabla \bS\|_{L^q_{\bx,s}}(t) 
\leq \mathscr S_\infty R^\star + \mathscr S_\infty' \frac{{R^\star}^2}{\We \, \mu} \mathrm e^{\mu t}.
\end{equation}
Taking the $L^r$-norm we finally obtain
\begin{equation}\label{120}
\| \div \, \bsigma\|_{L^r(0,t^\star;L^q_{\bx})} \leq C \big( {t^\star}^{\frac{1}{r}} R^\star + (\mathrm e^{r \mu t^\star} - 1)^{\frac{1}{r}} {R^\star}^2 \big) ,
\end{equation}
where the constant~$C$ does not depend on~$t^\star$ nor~$R^\star$.

$\checkmark$ {\bf Step 2: Control of the velocity~$\overline \bv$} --
We use the following classical result~\cite{Simon}: if $X\subset Y$ are two Banach spaces with compact injection then the following injection is continuous:
\begin{equation*}\label{121}
\big\{ \overline \bv \in L^2(0,t^\star;X) ~;~ \partial_t \overline \bv \in L^2(0,t^\star;Y) \big\} \subset \mathscr C(0,t^\star;[X,Y]_{\frac{1}{2}}).
\end{equation*}
Using $X=W^{2,q}_{\bx}$ and $Y=L^q_{\bx}$ we obtain
\begin{equation}\label{122}
\begin{aligned}
\| \overline \bv \|_{L^\infty(0,t^\star;W^{1,q}_{\bx})}
& \leq C \, \| \overline \bv \|_{L^2(0,t^\star;W^{2,q}_{\bx})}^{\frac{1}{2}} \| \partial_t \overline \bv \|_{L^2(0,t^\star;L^q_{\bx})}^{\frac{1}{2}} \\
& \leq C \, {t^\star}^{\frac{r-2}{2r}} \| \overline \bv \|_{L^r(0,t^\star;W^{2,q}_{\bx})}^{\frac{1}{2}} \| \partial_t \overline \bv \|_{L^r(0,t^\star;L^q_{\bx})}^{\frac{1}{2}} \\
& \leq C \, {t^\star}^{\frac{r-2}{2r}} R^\star.
\end{aligned}
\end{equation}
More simply (in fact using the Hölder inequality and the continuous Sobolev embedding $W^{2,q}_{\bx} \hookrightarrow W^{1,\infty}_{\bx}$) we have
\begin{equation}\label{123}
\begin{aligned}
& \| \overline \bv \|_{L^1(0,t^\star;W^{2,q}_{\bx})} \leq {t^\star}^{\frac{r-1}{r}} \| \overline \bv \|_{L^r(0,t^\star;W^{2,q}_{\bx})} \leq {t^\star}^{\frac{r-1}{r}} R^\star, \\
&  \| \overline \bv \|_{L^2(0,t^\star;W^{1,\infty}_{\bx})} \leq C \, {t^\star}^{\frac{r-2}{2r}} \| \overline \bv \|_{L^r(0,t^\star;W^{2,q}_{\bx})} \leq C \, {t^\star}^{\frac{r-2}{2r}} R^\star.
\end{aligned}
\end{equation}


$\checkmark$ {\bf Step 3: Control of the source term~$g$} --
Using the definition~\eqref{110} of~$\partial_s g$ we directly obtain
\begin{equation}\label{123.5}
\| \partial_s g\|_{L^2(0,t^\star;L^\infty_{\bx,s})} \leq \| \nabla \overline \bv\|_{L^2(0,t^\star;L^\infty_{\bx})} \| \bS\|_{L^\infty(0,t^\star;L^\infty_{\bx,s})}.
\end{equation}
The velocity contribution is controlled using the Sobolev embedding $W^{1,q}_{\bx} \hookrightarrow L^\infty_{\bx}$:
\begin{equation}\label{123.6}
\| \nabla \overline \bv\|_{L^2(0,t^\star;L^\infty_{\bx})}
\leq
C \| \overline \bv\|_{L^2(0,t^\star;W^{2,q}_{\bx})}
\leq
C {t^\star}^{\frac{r-2}{2r}} \| \overline \bv\|_{L^r(0,t^\star;W^{2,q}_{\bx})}
\leq
C {t^\star}^{\frac{r-2}{2r}} R^\star.
\end{equation}
The contribution of~$\bS$ in the inequality~\eqref{123.5} is estimate using its definition~\eqref{38} and using the inequality $|\partial_T K| \leq R^\star \mathrm e^{\mu t}\mathrm e^{-\mu \We T}$ which comes from the bound~$\vvvert \overline K \vvvert_3 \leq R^\star$:
\begin{equation}\label{123.7}
\| \bS\|_{L^\infty(0,t^\star;L^\infty_{\bx,s})} \leq C \, R^\star \, \mathrm e^{\mu t^\star}.
\end{equation}
The estimates~\eqref{123.6} and~\eqref{123.7} allows to write the inequality~\eqref{123.5} as follows
\begin{equation}\label{123.8}
\| \partial_s g\|_{L^2(0,t^\star;L^\infty_{\bx,s})} \leq C {t^\star}^{\frac{r-2}{2r}} {R^\star}^2 \, \mathrm e^{\mu t^\star}.
\end{equation}
We complete this step considering~$\partial_s \nabla g$. We have
\begin{equation*}\label{124}
\begin{aligned}
\| \partial_s \nabla g\|_{L^2(0,t^\star;L^q_{\bx,s})}
\leq
& \| \nabla \overline \bv\|_{L^2(0,t^\star;L^\infty_{\bx})}  \| \nabla \bS\|_{L^\infty(0,t^\star;L^q_{\bx,s})} \\
& +  \| \nabla^2 \overline \bv\|_{L^2(0,t^\star;L^q_{\bx})}  \| \bS\|_{L^\infty(0,t^\star;L^\infty_{\bx,s})}.
\end{aligned}
\end{equation*}
We use estimates~\eqref{123.6} and~\eqref{123.7}, supplemented by estimates (the first comes from to estimate~\eqref{119}, the second is a direct consequence of the bound~$\vvvert \overline \bv \vvvert_1 \leq R^\star$):
\begin{equation*}\label{125}
\begin{aligned}
& \|\nabla \bS\|_{L^\infty(0,t^\star;L^q_{\bx,s})} \leq C \, \big( R^\star + {R^\star}^2 \, \mathrm e^{\mu t^\star} \big), \\
& \| \nabla^2 \overline \bv\|_{L^2(0,t^\star;L^q_{\bx})} \leq C t^{\frac{r-2}{2r}} \| \overline \bv\|_{L^r(0,t^\star;W^{2,q}_{\bx})} \leq C t^{\frac{r-2}{2r}} R^\star,
\end{aligned}
\end{equation*}
to deduce
\begin{equation}\label{126}
\| \partial_s \nabla g\|_{L^2(0,t^\star;L^q_{\bx,s})} \leq C \, {t^\star}^{\frac{r-2}{2r}} {R^\star}^2 \big( 1 + \mathrm e^{\mu t^\star} (1+ R^\star) \big).
\end{equation}
Gathering the estimates \eqref{113}, \eqref{120}, \eqref{122}, \eqref{123}, \eqref{123.8} and~\eqref{126} we conclude the proof of lemma~\ref{lem:1648}. 
\cqfd
{\bf $\Phi$-Invariant subset} --
Consequently, using the lemmas~\ref{lem:stokes}, \ref{lem:1116}, \ref{lem:1117} and~\ref{lem:1648}, we deduce that if $(\overline \bv,\overline \bG,\overline K)\in \mathscr{H}(t^\star,R^\star)$ then its image $(\bv,\bG,K) = \Phi(\overline \bv,\overline \bG,\overline K)$ satisfies
\begin{equation*}\label{127}
\begin{aligned}
& \vvvert \bv \vvvert_1 \leq F_1\big( G(t^\star, R^\star) \big), \\
& \vvvert \bG \vvvert_2 \leq F_2\big( G(t^\star, R^\star) \big), \\
& \vvvert K \vvvert_3 \leq F_3\big( G(t^\star, R^\star) , G(t^\star, R^\star) \big).
\end{aligned}
\end{equation*}
To ensure that the ball~$\mathscr B(t^\star)$ is invariant under the action of~$\Phi$, we must find $t^\star$ and~$R^\star$ such that
\begin{equation}\label{128}
\begin{aligned}
& F_1\big( G(t^\star, R^\star) \big) \leq R^\star, \\
& F_2\big( G(t^\star, R^\star) \big) \leq R^\star, \\
& F_3\big( G(t^\star, R^\star) , G(t^\star, R^\star) \big) \leq R^\star.
\end{aligned}
\end{equation}
Notice that if we choose $t^\star=0$ then the previous inequalities~\eqref{128} become (the function~$G$ vanishes for $t^\star=0$)
\begin{equation*}\label{129}
\begin{aligned}
& F_1\big( 0 \big) \leq R^\star, \\
& F_2\big( 0 \big) \leq R^\star, \\
& F_3\big( 0 , 0 \big) \leq R^\star.
\end{aligned}
\end{equation*}
By continuity argument, taking
\begin{equation*}\label{130}
R^\star = \max \big\{ 2 F_1(0) \, , \, F_2(0) +  \vvvert \bG_0 \vvvert_2 \, , \, F_3(0,0) + \vvvert K_0 \vvvert_3 \big\},
\end{equation*}
we deduce that there exists $t^\star>0$ such that~\eqref{128} holds. For such a choice we have the inclusion $\Phi(\mathscr B(t^\star))\subset \mathscr B(t^\star)$.

\noindent
Note that at the same time , we chose~$R^\star$ so that the inequalities~\eqref{107} hold, this insures that $\mathscr{H}(t^\star,R^\star) \neq \emptyset$.
Moreover the function $\Phi$ is continuous and $\mathscr{H}(t^\star,R^\star)$ is a convex compact subset of $\mathscr{B}(t^\star)$, see~\cite{Guillope-Saut3} for similar properties.
We conclude the proof using the Schauder's theorem.
\cqfd

\section{Proof of the uniqueness result}\label{part:proof2}

This section is devoted to the proof of the theorem~\ref{th:unic}.
As usual to prove an uniqueness result we take the difference of the two solutions indexed by~$1$ and~$2$ such that
\begin{equation}\label{131}
\begin{aligned}
& \nabla \bv_i \in L^2(0,t^\star;L^\infty_{\bx}), \\
& \bG_i \in L^2(0,t^\star;L^\infty_T (L^\infty_{\bx} \cap W^{1,d}_{\bx})), \\
& \partial_T K_i \in L^\infty(0,t^\star;L^\infty_T L^\infty_{\bx,s} \cap L^1_T L^\infty_{\bx,s}), \\
& \partial_T \nabla K_i \in L^\infty(0,t^\star;L^2_T L^d_{\bx,s}).
\end{aligned}
\end{equation}
The vector $\bv=\bv_1-\bv_2$, the scalars $p=p_1-p_2$, $K=K_1-K_2$ and the tensor $\bG=\bG_1-\bG_2$ satisfy the following:
\begin{subequations} \label{systeme-unic}
  \begin{align}
& \Re ( \partial_t \bv + \bv_1\cdot \nabla \bv + \bv\cdot \nabla \bv_2 ) + \nabla p - (1-\omega) \Delta \bv = \div \, \bsigma,\label{eq-unic:1} \\
& \div \, \bv = 0,\label{eq-unic:2} \\
& \bsigma(t,\bx) = \omega \int_{-\frac{1}{2}}^{\frac{1}{2}} \bS(t,\bx,s)\, \mathrm ds,\label{eq-unic:3}\\
& \bS(t,\bx,s) = \int_0^{+\infty} \!\! \Big( -\partial_{T} K_1\, \big(\mathscr S(\bG_1) - \mathscr S(\bG_2) \big) - \partial_{T} K \, \mathscr S(\bG_2) \Big) \mathrm dT,\label{eq-unic:4}\\
& \partial_t \bG + \bv_1 \cdot \nabla \bG + \bv \cdot \nabla \bG_2 + \frac{1}{\We} \partial_T \bG = \bG_1 \cdot \nabla \bv + \bG \cdot \nabla \bv_2,\label{eq-unic:5}\\
& \partial_t K + \bv_1\cdot \nabla K + \bv\cdot \nabla K_2
+ \frac{1}{\We} \partial_T K - \frac{1}{\We} \partial^2_s K + \Big( \nabla \bv_1: \int_0^s \bS_1 \Big) \partial_s K \nonumber \\
& \hspace{3.5cm} + \Big( \nabla \bv_1: \int_0^s \bS \Big) \partial_s K_2 + \Big( \nabla \bv: \int_0^s \bS_2 \Big) \partial_s K_2 = 0.\label{eq-unic:6}
  \end{align}
\end{subequations}
together with zero initial and boundary conditions.
The uniqueness proof consists in demonstrate that $\bv = \b0$, $\bG=\b0$ and $K=0$.
We will initially provide estimates on these three quantities.
More exactly, we introduce the following quantities
\begin{equation}\label{131.5}
\begin{aligned}
& W(t)=\|\nabla \bv\|_{L^2_{\bx}}^2, \\
& X(t)=\Re\,\|\bv\|_{L^2_{\bx}}^2, \\
& Y(t) = \int_0^{+\infty} \| m_1\|_{L^\infty_{\bx,s}}^2 \| \bG \|_{L^2_{\bx,s}}^2 \, \mathrm dT, \\
& Z(t) = \int_0^{+\infty}  \| m \|_{L^2_{\bx,s}}^2 \, \mathrm dT,
\end{aligned}
\end{equation}
and we will obtained some relations between them.\\[0.3cm]
{\bf Velocity estimate --}
Taking the inner product of the equation~\eqref{eq-unic:1} by~$\bv$ in $L^2_{\bx}$, we obtain
\begin{equation*}\label{132}
\frac{\Re}{2} d_t \|\bv\|_{L^2_{\bx}}^2 + (1-\omega) \|\nabla \bv\|_{L^2_{\bx}}^2
=
- \int_\T \bsigma \cdot \nabla \bv - \Re \int_\T (\bv \cdot \nabla \bv_2) \cdot \bv.
\end{equation*}
From the Cauchy-Schwarz inequality and the Young inequality, we obtain
\begin{equation}\label{133}
\Re \, d_t \|\bv\|_{L^2_{\bx}}^2 + (1-\omega) \|\nabla \bv\|_{L^2_{\bx}}^2
\leq \, \frac{1}{1-\omega}\|\bsigma\|_{L^2_{\bx}}^2 + 2 \Re \|\nabla \bv_2\|_{L^\infty_{\bx}} \|\bv\|_{L^2_{\bx}}^2.
\end{equation}
Introducing $a_1(t) = 2\|\nabla \bv_2\|_{L^\infty_{\bx}}$ this estimate~\eqref{133} reads using the notations~$X$ and~$W$ introduced below (see~\eqref{131.5}):
\begin{equation}\label{134}
X' + (1-\omega) W
\leq
\frac{1}{1-\omega} \|\bsigma\|_{L^2_{\bx}}^2  + a_1 X.
\end{equation}
It is important to notice that, due to the assumptions~\eqref{131}, we have $a_1\in L^1(0,t^\star)$.\\

{\bf Stress tensor estimate --}
From the definition of the stress tensor~$\bsigma$ in the System~\eqref{systeme-unic} we have $\| \bsigma \|_{L^2_{\bx}} \leq \omega \| \bS \|_{L^2_{\bx,s}}$ and due to the definition of $\bS$ (see equation~\eqref{eq-unic:4}), we obtain
\begin{equation*}\label{135}
\begin{aligned}
\| \bsigma \|_{L^2_{\bx}}^2 \leq 
& \omega^2 \int_0^{+\infty} \| m_1\|_{L^\infty_{\bx,s}}^2 \| \mathscr S(\bG_1) - \mathscr S(\bG_2) \|_{L^2_{\bx,s}}^2 \\
& +  \omega^2 \int_0^{+\infty}  \| m \|_{L^2_{\bx,s}}^2 \| \mathscr S(\bG_2) \|_{L^\infty_{\bx,s}}^2 \mathrm dT,
\end{aligned}
\end{equation*}
where $m_i=-\partial_T K_i$, $i\in \{1,2\}$ and $m=m_1-m_2$.
By assumption, the function~$\mathscr S$ and~$\mathscr S'$ are bounded by~$\mathscr S_\infty$ and~$\widetilde{\mathscr S_\infty'}$ respectively (in practice, recall that we work with the function~$\widetilde{\mathscr S}$, see the preliminary section~\ref{positive-deformation}). In particular we have $|\mathscr S(\bG_1)-\mathscr S(\bG_2)|\leq \widetilde{\mathscr S_\infty}' |\bG_1-\bG_2|$.
In term of $Y$ and~$Z$ (see their definition given by~\eqref{131.5}) we deduce
\begin{equation}\label{136}
\| \bsigma \|_{L^2_{\bx}}^2 \leq \omega^2 \widetilde{\mathscr S_\infty'}^2 Y + \omega^2 \mathscr S_\infty^2 Z.
\end{equation}

{\bf Deformation gradient estimate --}
Taking the inner product of the equation~\eqref{eq-unic:5} by $\bG$ in $L^2_{\bx}$, we obtain
\begin{equation*}\label{138}
\begin{aligned}
\frac{1}{2} \partial_t \|\bG\|_{L^2_{\bx}}^2 + \frac{1}{2\We} \partial_T \|\bG\|_{L^2_{\bx}}^2 & = \int_\T (\bG_1 \cdot \nabla \bv)\cdot \bG \\
& + \int_\T (\bG \cdot \nabla \bv_2)\cdot \bG - \int_\T (\bv \cdot \nabla \bG_2)\cdot \bG.
\end{aligned}
\end{equation*}
Using the Cauchy-Schwarz and the Hölder inequalities, we have the estimate
\begin{equation*}\label{139}
\begin{aligned}
\frac{1}{2} \partial_t \|\bG\|_{L^2_{\bx}}^2 + \frac{1}{2\We} \partial_T \|\bG\|_{L^2_{\bx}}^2
\leq
& \|\bG_1\|_{L^\infty_{\bx}} \|\nabla \bv\|_{L^2_{\bx}} \|\bG\|_{L^2_{\bx}} \\
& + \|\nabla \bv_2\|_{L^\infty_{\bx}} \|\bG\|_{L^2_{\bx}}^2 \\
& + \|\bv\|_{L^{\frac{2d}{d-2}}_{\bx}} \|\nabla \bG_2\|_{L^d_{\bx}} \|\bG\|_{L^2_{\bx}}.
\end{aligned}
\end{equation*}
Due to the Sobolev continuous injection $H^1_{\bx} \hookrightarrow L^{\frac{2d}{d-2}}_{\bx}$ and the Young inequality, we obtain:
\begin{equation}\label{140}
\partial_t \|\bG\|_{L^2_{\bx}}^2 + \frac{1}{\We} \partial_T \|\bG\|_{L^2_{\bx}}^2
\leq
\|\bv\|_{H^1_{\bx}}^2 + a_2(t) \|\bG\|_{L^2_{\bx}}^2,
\end{equation}
where the function 
\begin{equation*}\label{141}
a_2(t) = 2 \|\bG_1\|_{L^\infty_T L^\infty_{\bx}}^2 + 2\|\nabla \bv_2\|_{L^\infty_{\bx}} + C^2 \, \|\nabla \bG_2\|_{L^\infty_T L^d_{\bx}}^2,
\end{equation*}
and where the constant~$C$ depends on~$q$ and~$d$.
Multiplying this estimate~\eqref{140} by $\| m_1\|_{L^\infty_{\bx,s}}^2$ where $m_1 = -\partial_{T} K_1$, and integrating for $T\in (0,+\infty)$ we obtain
\begin{equation}\label{142}
Y'(t) + \frac{1}{\We} \underbrace{\int_0^{+\infty}  \| m_1\|_{L^\infty_{\bx,s}}^2 \partial_T \|\bG\|_{L^2_{\bx}}^2 \, \mathrm dT}_{\mathscr I}
\leq
C \big( X(t)+W(t) \big) + a_2(t) Y(t).
\end{equation}
We note that the constant $C$ contains the value $\int_0^{+\infty}  \| m_1\|_{L^\infty_{\bx,s}}^2 \mathrm dT$ which is bounded according to~\eqref{131}.
\par\noindent
Using an integration by parts we prove that the integral~$\mathscr I$ is non-negative.
More precisely we use the fact that $m_1\big|_{T=0}=0$ and the fact that $M:T\mapsto \|m_1\|_{L^\infty_{\bx,s}}(t,T)$ is a decreasing function.
Indeed, we simply use the following property of the supremum: for any $0<T<T'$ we have
\begin{equation*}\label{143}
M(T') - M(T) \leq \sup_{\bx,s} \big( m_1(t,T',\bx,s) - m_1(t,T,\bx,s) \big).
\end{equation*}
Since $\partial_T m_0 \leq 0$ we have seen that $\partial_T m_1 \leq 0$ (see the discussion page~\pageref{th:local} before the statement of the theorem~\ref{th:unic}). We conclude that $M(T') - M(T) \leq 0$.
\par\noindent
Consequently the estimate~\eqref{142} now reads
\begin{equation}\label{144}
Y' \leq C X+ b_1 W + a_2 Y.
\end{equation}
We also note that $b_1=C$ and, due to the assumptions~\eqref{131}, we have $a_2\in L^1(0,t^\star)$.\\

{\bf Memory estimate --}
We derivate the equation~\eqref{eq-unic:6} with respect to the age~$T$ and next we take the inner product with $m=-\partial_T K$ in $L^2_{\bx,s}$. Using integrations by parts we obtain
\begin{equation*}\label{145}
\begin{aligned}
\frac{1}{2}\partial_t \|m\|_{L^2_{\bx,s}}^2 + \frac{1}{2\We} \partial_T \|m\|_{L^2_{\bx,s}}^2 + & \frac{1}{\We} \|\partial_s m\|_{L^2_{\bx,s}}^2 \leq  
\frac{1}{2} \int_\T \int_{-\frac{1}{2}}^{\frac{1}{2}} (\nabla \bv_1 : \bS_1) m^2 \\
& + \int_\T \int_{-\frac{1}{2}}^{\frac{1}{2}} \Big( \nabla \bv_1: \bS + \nabla \bv: \bS_2 \Big) m_2 \, m \\
& + \int_\T \int_{-\frac{1}{2}}^{\frac{1}{2}} \Big( \nabla \bv_1: \int_0^s \bS + \nabla \bv: \int_0^s \bS_2 \Big) m_2 \, \partial_s m \\
& - \int_\T \int_{-\frac{1}{2}}^{\frac{1}{2}} (\bv \cdot \nabla m_2) m.
\end{aligned}
\end{equation*}
By the Cauchy-Schwarz inequality we obtain the following estimate:
\begin{equation*}\label{146}
\begin{aligned}
\frac{1}{2}\partial_t \|m\|_{L^2_{\bx,s}}^2 + \frac{1}{2\We} \partial_T \|m\|_{L^2_{\bx,s}}^2 + & \frac{1}{\We} \|\partial_s m\|_{L^2_{\bx,s}}^2 \leq  
\frac{1}{2} \|\nabla \bv_1\|_{L^\infty_{\bx}} \|\bS_1\|_{L^\infty_{\bx,s}} \|m\|_{L^2_{\bx,s}}^2 \\
& + \|\nabla \bv_1\|_{L^\infty_{\bx}} \|\bS\|_{L^2_{\bx,s}}  \|m_2\|_{L^\infty_{\bx,s}}  \|m\|_{L^2_{\bx,s}} \\
& + \|\nabla \bv\|_{L^2_{\bx}}  \|\bS_2\|_{L^\infty_{\bx,s}}  \|m_2\|_{L^\infty_{\bx,s}} \|m\|_{L^2_{\bx,s}} \\
& + \|\nabla \bv_1\|_{L^\infty_{\bx}} \|\bS\|_{L^2_{\bx,s}} \|m_2\|_{L^\infty_{\bx,s}} \|\partial_s m\|_{L^2_{\bx,s}} \\
& + \|\nabla \bv\|_{L^2_{\bx}} \|\bS_2\|_{L^\infty_{\bx,s}} \|m_2\|_{L^\infty_{\bx,s}} \|\partial_s m\|_{L^2_{\bx,s}} \\
& + \|\bv\|_{L^{\frac{2d}{d-2}}_{\bx}} \|\nabla m_2\|_{L^d_{\bx,s}} \|m\|_{L^2_{\bx,s}}.
\end{aligned}
\end{equation*}
From the Young inequality $2ab\leq a^2+b^2$ and the Sobolev embedding for the last term we deduce
\begin{equation}\label{147}
\begin{aligned}
\partial_t \|m\|_{L^2_{\bx,s}}^2 + \frac{1}{\We} \partial_T \|m\|_{L^2_{\bx,s}}^2 
\leq \, & \widetilde{C_Z}(t,T) \|m\|_{L^2_{\bx,s}}^2 \\
& + \|\nabla \bv_1\|_{L^\infty_{\bx}} \|m_2\|_{L^\infty_{\bx,s}} \|\bS\|_{L^2_{\bx,s}}^2 \\
& + \|\bS_2\|_{L^\infty_{\bx,s}} \|m_2\|_{L^\infty_{\bx,s}} \|\nabla \bv\|_{L^2_{\bx}}^2 \\
& + \We \|\nabla \bv_1\|_{L^\infty_{\bx}}^2 \|m_2\|_{L^\infty_{\bx,s}}^2 \|\bS\|_{L^2_{\bx,s}}^2 \\
& + \We \|\bS_2\|_{L^\infty_{\bx,s}}^2 \|m_2\|_{L^\infty_{\bx,s}}^2 \|\nabla \bv\|_{L^2_{\bx}}^2 \\
& + C \|\nabla m_2\|_{L^d_{\bx,s}}^2 \|\bv\|_{H^1_{\bx}}^2,
\end{aligned}
\end{equation}
where the function
\begin{equation*}\label{148}
\widetilde{C_Z}(t,T) = \|\nabla \bv_1\|_{L^\infty_{\bx}} \|\bS_1\|_{L^\infty_{\bx,s}} + \|\nabla \bv_1\|_{L^\infty_{\bx}} \|m_2\|_{L^\infty_{\bx,s}} + \|\bS_2\|_{L^\infty_{\bx,s}} \|m_2\|_{L^\infty_{\bx,s}} + 1.
\end{equation*}
Integrating the estimate~\eqref{147} for $T\in (0,T)$ we deduce (we also use the fact that $\|\bS_2\|_{L^2_{\bx,s}} \leq \mathscr S_\infty \|m_2\|_{L^1_T L^\infty_{\bx,s}}$)
\begin{equation}\label{149}
\begin{aligned}
Z'(t) \leq \, C_Z(t) Z(t)
& + \|\nabla \bv_1\|_{L^\infty_{\bx}} \|m_2\|_{L^1_T L^\infty_{\bx,s}} \|\bS\|_{L^2_{\bx,s}}^2 \\
& + \mathscr S_\infty \|m_2\|_{L^1_T L^\infty_{\bx,s}}^2 W(t) \\
& + \We \|\nabla \bv_1\|_{L^\infty_{\bx}}^2 \|m_2\|_{L^2_T L^\infty_{\bx,s}}^2 \|\bS\|_{L^2_{\bx,s}}^2 \\
& + \We \, \mathscr S_\infty^2 \|m_2\|_{L^1_T L^\infty_{\bx,s}}^2 \|m_2\|_{L^2_T L^\infty_{\bx,s}}^2 W(t) \\
& + C \|\nabla m_2\|_{L^2_T L^d_{\bx,s}}^2 \big(X(t)+W(t)\big),
\end{aligned}
\end{equation}
where $\displaystyle C_Z(t) = \sup_{T\in \R^+} \widetilde{C_Z}(t,T)$. Recalling that $\|\bS\|_{L^2_{\bx,s}}^2 \leq \widetilde{\mathscr S_\infty'}^2 Y+ \mathscr S_\infty^2 Z$, we can write the previous estimate~\eqref{149} as follows
\begin{equation}\label{150}
Z' \leq \, a_3 X + a_4 Y + a_5 Z + b_2 W,
\end{equation}
where we have defined
\begin{equation*}\label{151}
\begin{aligned}
a_3(t) = & C \|\nabla m_2\|_{L^2_T L^d_{\bx,s}}^2 , \\
a_4(t) = & C \big( \|\nabla \bv_1\|_{L^\infty_{\bx}} \|m_2\|_{L^1_T L^\infty_{\bx,s}} + \|\nabla \bv_1\|_{L^\infty_{\bx}}^2 \|m_2\|_{L^2_T L^\infty_{\bx,s}}^2 \big) , \\
a_5(t) = & C \big( \|\nabla \bv_1\|_{L^\infty_{\bx}} \|m_1\|_{L^1_T L^\infty_{\bx,s}}
+ \|\nabla \bv_1\|_{L^\infty_{\bx}} \|m_2\|_{L^1_T L^\infty_{\bx,s}} \\
& \qquad + \|\nabla \bv_1\|_{L^\infty_{\bx}} \|m_2\|_{L^\infty_T L^\infty_{\bx,s}}
+ \|m_2\|_{L^1_T L^\infty_{\bx,s}} \|m_2\|_{L^\infty_T L^\infty_{\bx,s}} \\
& \hspace{4.2cm} + \|\nabla \bv_1\|_{L^\infty_{\bx}}^2 \|m_2\|_{L^2_T L^\infty_{\bx,s}}^2 + 1 \big).
\end{aligned}
\end{equation*}
The assumptions~\eqref{131} imply that~$a_3$, $a_4$ and~$a_5$ belong in $L^1(0,t^\star)$.
The last contribution makes appear the function~$b_2$ which is given by
\begin{equation*}\label{152}
b_2(t) = \mathscr S_\infty \|m_2\|_{L^1_T L^\infty_{\bx,s}}^2 + \We \, \mathscr S_\infty^2 \|m_2\|_{L^1_T L^\infty_{\bx,s}}^2 \|m_2\|_{L^2_T L^\infty_{\bx,s}}^2 + C \|\nabla m_2\|_{L^2_T L^d_{\bx,s}}^2.
\end{equation*}
The assumptions~\eqref{131} imply that~$b_2\in L^\infty(0,t^\star)$.\\

{\bf Uniqueness result --}
Finally, we perform the following combination:
\begin{equation*}\label{153}
(1-\omega)\eqref{134} + \eqref{136} + \varepsilon \eqref{144} + \varepsilon \eqref{150}
\end{equation*}
using $\varepsilon>0$ small enough to control the quantity~$W(t)$ of the right side member by the same quantity on the left side member.
This is possible since the functions~$b_1$ and~$b_2$ are bounded. More precisely, we use
\begin{equation*}\label{154}
\varepsilon = \frac{(1-\omega)^2}{1+b_1+ \|b_2\|_{L^\infty(0,t^\star)}}.
\end{equation*}
We deduce an estimate of kind
\begin{equation*}\label{155}
(X+Y+Z)' \leq a \, (C+Y+Z),
\end{equation*}
where the function~$a$ is a linear combination of~$a_i$, $i\in \{1,2,3,4\}$, and consequently it belongs in $L^1(0,t^\star)$.
The classical Gronwall lemma and the initial conditions $X(0)=Y(0)=Z(0)=0$ imply that $X=Y=Z=0$. We deduce that $\bv=\b0$, $\bG=\b0$ and~$K=0$, that concludes the proof.
\cqfd

\section{Global existence}\label{part:proof3}

In this section, we prove the main result of the present paper, that is the global existence of a solution to the System~\eqref{system} in the two dimensional case, see Theorem~\ref{th:global}.
The proof is based on the following fundamental remark: the stress given by the relation~\eqref{eq:3} is ``naturally'' bounded. This bound implies a bound on the velocity field, from which we can deduce a new bound on the stress gradient. The conclusion of the proof consists in to insure that this process is consistent.\\[0.2cm]
\noindent
The main difficulties are the evaluation of the nonlinear terms that can blow the estimates obtained in the proof of local existence. These nonlinearities are present many times in the Doi-Edwards model:
\begin{enumerate}
\item The nonlinearities of the Navier-Stokes equations~\eqref{eq:1}--\eqref{eq:2} are always present in the complete Doi-Edwards model.
On one hand there is the convection term~$\bv\cdot \nabla \bv$ which is discussed in many articles. We handle this first difficulty using the so-called $L^pL^q$ estimates of the heat kernel. On the other hand there is the nonlinear coupling due to the presence of the stress term~$\div \, \bsigma$.
\item Other nonlinearities are present in the evolution equation~\eqref{eq:5} for the deformation tensor~$\bG$. The most restrictive one corresponds to the term~$\bG \cdot \nabla \bv$. Since we will have relatively little information on the velocity gradient, we will prefer to work with the quotient~$\frac{\nabla \bG}{|\bG|}$, making appear the product~$|\bG| \mathscr S(\bG)$. The key point is the fact that this product is well behaved even for large value of~$\bG$.
\item Finally, there are also several nonlinearities in the equation~\eqref{eq:6} describing the evolution of memory~$K$. The most restrictive one comes from the product $\big( \nabla \bv : \int_0^s \bS \big) \partial_s K$. To solve this problem, we introduce a suitable function~$\xi$ in order to obtain a better estimate on~$\xi(\partial_T K)$, and therefore a better estimate on~$\partial_T K$.
\end{enumerate}
\vspace{0.2cm}
\noindent
Let $r>2$, $q>2$, $\gamma>0$ and $\mu>0$.
We consider the initial data $(\bv_0,\bG_0,K_0)$ satisfying
\begin{equation*}\label{156}
\begin{aligned}
& \bv_0 \in \mathscr D_q^r, \\
& \bG_0\in L^\infty_T W^{1,q}_{\bx},
\quad \partial_T \bG_0\in L^\infty_T L^q_{\bx},
\quad \det \bG_0 \geq \gamma, \\
& \mathrm e^{\mu T} \partial_T K_0 \in L^\infty_T L^\infty_{\bx,s},
\quad \mathrm e^{\mu T/2} \partial_T \nabla K_0 \in L^2_T L^q_{\bx,s},
\quad \partial_s K_0 \in L^\infty_T L^2_{\bx,s}, \\
& \text{and} \quad \partial_T K_0 \leq 0.
\end{aligned}
\end{equation*}
It is then possible to use the Theorem~\ref{th:local} and then consider the solution~$(\bu,\bG,K)$ to System~\eqref{system} in $[0,t^\star]$, which satisfies the initial conditions~\eqref{9}.
This solution possesses at least the following regularity:
\begin{equation*}\label{157}
\begin{array}{ll}
\bv \in L^r(0,t^\star;W^{2,q}_{\bx}), \phantom{totototototototototototototo}
& \partial_t \bv \in L^r(0,t^\star;L^q_{\bx}), \\
\bG \in L^\infty(0,t^\star,L^\infty_T W^{1,q}_{\bx}),
& \partial_s \bG, ~\partial_t  \bG \in L^r(0,t^\star;L^\infty_T L^q_{\bx}),\\
K, \partial_T K, \mathrm e^{\mu \, (\We T-t)} \partial_T K \in L^\infty(0,t^\star;L^\infty_T L^\infty_{\bx,s}),
& \partial_t K \in L^2(0,t^\star;L^\infty_T L^2_{\bx,s}),\\
\partial_T \nabla K \in L^\infty(0,t^\star;L^1_T L^q_{\bx,s} \cap L^2_T L^q_{\bx,s}),
& \partial_s K \in L^\infty(0,t^\star;L^\infty_T L^2_{\bx,s}). 
\end{array}
\end{equation*}
\noindent
In the following and as previously we recall that we denote by~$C$ constants that may depend on the initial conditions, on the physical parameters, on the integers~$r$, $q$, on the bounds~$\mathscr S_\infty$ and~$\mathscr S_\infty'$, and on the time~$t^\star$. Note that these constants will always be bounded for bounded~$t^\star$.

\subsection{Maximum principle for the kernel $\boldsymbol{K}$ and their derivatives}

First we resume the proof of the local existence theorem to deduce the following bounds (see the estimates~\eqref{65}, \eqref{69} and~\eqref{72}):
\begin{lemma}\label{lem:2027}
There exists a constant $C$ such that
\begin{equation*}\label{158}
\begin{aligned}
& |K| \leq C,\\
& 0 \leq -\partial_T K \leq C,\\
& 0 \leq -\mathrm e^{\mu \, (\We \, T-t)} \partial_T K \leq C.
\end{aligned}
\end{equation*}
\end{lemma}

\subsection{Additional bounds for the stress tensor and for the velocity field}

The first Proposition that we introduce is one of the fundamental points which give the global existence result. In particular, this kind of result is not proved for other viscoelastic systems like the Oldroyd models.\\[-0.2cm]
\begin{proposition}\label{prop:1}
We have the following $L^\infty$-bound:
\begin{equation*}\label{172.5}
\|\bsigma\|_{L^\infty(0,t^\star;L^\infty_{\bx})} \leq C.
\end{equation*}
\end{proposition}
\proof
Recall that the stress is given by the integral relation (see~\eqref{eq:3}):
\begin{equation*}\label{173}
\bsigma(t,\bx) = \omega \int_{-\frac{1}{2}}^{\frac{1}{2}} \bS(t,\bx,s)\, \mathrm ds.
\end{equation*}
To have a bound on~$\bsigma$ it suffice to have a bound on~$\bS$ which is given by (see~\eqref{eq:4}):
\begin{equation*}\label{174}
\bS(t,\bx,s) = -\int_0^{+\infty} \partial_{T} K(t,T,\bx,s)\, \mathscr S(\bG(t,T,\bx))\, \mathrm dT.
\end{equation*}
Since the function $\mathscr S$ is bounded on $\mathcal L(\R^d) \setminus\{\b0\}$ (see Proposition~\ref{propS}), and since the function~$\bG$ has values in $\mathcal L(\R^d) \setminus\{\b0\}$ (see Lemma~\ref{lemma7}), we deduce that the composed function $\mathscr S(\bG)$ is bounded.

\noindent
Moreover from Lemma~\ref{lem:2027} we know that $-\partial_T K \geq 0$ and we deduce that
\begin{equation}\label{175}
|\bS(t,\bx,s)| \leq - \mathscr S_\infty \int_0^{+\infty} \partial_{T} K(t,T,\bx,s)\, \mathrm dT = \mathscr S_\infty (K\big|_{T=0} - \lim_{T\to +\infty} K).
\end{equation}
Using the fact that $K$ is bounded (see Lemma~\ref{lem:2027} again) we conclude that~$\bS$ is bounded.
\cqfd
By virtue of this proposition~\ref{prop:1}, we can deduce that the velocity field~$\bv$ satisfying the Navier-Stokes equations~\eqref{eq:1}--\eqref{eq:2} have more regular. More precisely we prove the following Proposition (the proof is detailed in~\cite{Chupin4} and~\cite{Constantin-Masmoudi08}).\\[-0.2cm]
\begin{proposition}\label{prop:2}
If the integers~$r$ and~$q$ satisfy $\frac{1}{q} + \frac{1}{r} < \frac{1}{2}$ then for all $t\in (0,t^\star)$ we have
\begin{equation*}\label{176}
\|\nabla \bv\|_{L^\infty(0,t;L^\infty_{\bx})} \leq C \ln(\mathrm e+\|\nabla \bsigma\|_{L^r(0,t;L^q_{\bx})}),
\end{equation*}
\begin{equation*}\label{177}
\|\nabla^2 \bv\|_{L^r(0,t;L^q_{\bx})} \leq C \big( 1 + \|\nabla \bsigma\|_{L^r(0,t;L^q_{\bx})} \big).
\end{equation*}
\end{proposition}

\subsection{Control of the spatial gradient of the orientation tensor~$\bS$}

The goal is now to use the previous bounds (on the velocity field) to obtain a bound on the stress gradient~$\nabla \bsigma$. In fact, the stress gradient is directly linked to the orientation gradient~$\nabla \bS$ (see Equation~\eqref{eq:3}) and we are then interested in this quantity:\\[-0.2cm]
\begin{proposition}\label{prop:3}
For any convex function~$\txi:\R^+\to\R^+$ we have: for all $t\in (0,t^\star)$,
\begin{equation*}\label{178}
\|\nabla \bS\|_{L^r(0,t;L^q_{\bx,s})}^r \leq C \bigg( \Big( \int_0^\infty \frac{\mathrm e^{-2\We \mu T}}{\txi'(0)^2} \, \mathrm dT \Big)^{\frac{r}{2}} \, y_a(t) + y_b(t) \bigg),
\end{equation*}
where $y_a$ and~$y_b$ are defined by
\begin{equation*}\label{179}
y_a(t) =  \int_0^t \Big( \int_0^{+\infty} \|\nabla \tM\|_{L^q_{x,s}}^2(t',T) \, \mathrm dT \Big)^{\frac{r}{2}} \mathrm dt',
\end{equation*}
\begin{equation}\label{180}
y_b(t) = \int_0^t \int_0^{+\infty} \mathrm e^{r \mu t'} \mathrm e^{-\We \mu T} \, \Bigg\| \frac{\nabla \bG}{|\bG|} \Bigg\|_{L^q_{\bx}}^r(t',T) \, \mathrm dT \, \mathrm dt'.
\end{equation}
The quantity $\tM$ is defined by $\tM=\txi(\tm)$ and $\tm=\mathrm e^{\We\, \mu T} m$.
\end{proposition}\\[-0.2cm]
\proof
As previously we denote $m=-\partial_T K$. We derivate the stress tensor given by the relations~\eqref{eq:3} and~\eqref{eq:4} with respect to the spatial variable:
\begin{equation*}\label{181}
\begin{aligned}
\nabla \bS(t,\bx,s)
& =
\underbrace{\int_0^{+\infty} \nabla m(t,T,\bx,s) \otimes \mathscr S(\bG(t,T,\bx))\, \mathrm dT}_{(\nabla \bS)^a(t,\bx,s)}
\\
& \qquad +
\underbrace{\int_0^{+\infty} m(t,T,\bx,s) \, \mathscr S'(\bG(t,T,\bx)) : \nabla \bG(t,T,\bx) \, \mathrm dT}_{(\nabla \bS)^b(t,\bx,s)}.
\end{aligned}
\end{equation*}
We have the estimate
\begin{equation}\label{182}
\|\nabla \bS\|_{L^r(0,t;L^q_{\bx,s})}^r \leq 2^{r-1} \big( \|(\nabla \bS)^a\|_{L^r(0,t;L^q_{\bx,s})}^r + \|(\nabla \bS)^b\|_{L^r(0,t;L^q_{\bx,s})}^r \big) ,
\end{equation}
so that we will independently control $\|(\nabla \bS)^a\|_{L^r(0,t;L^q_{\bx,s})}$ and~$\|(\nabla \bS)^b\|_{L^r(0,t;L^q_{\bx,s})}$.\\

$\checkmark$ We use the definition of $\tM=\txi(\tm)$, $\tm=\mathrm e^{\We\, \mu T} m$ to write the quantity~$(\nabla \bS)^a$ as follows:
\begin{equation*}\label{182.5}
(\nabla \bS)^a(t,\bx,s) =
\int_0^{+\infty} \frac{\mathrm e^{-\We\, \mu T}}{\txi'(\tm(t,T,\bx,s))} \nabla \tM(t,T,\bx,s) \otimes \mathscr S(\bG(t,T,\bx))\, \mathrm dT.
\end{equation*}
We deduce that $\|(\nabla \bS)^a\|_{L^q_{\bx,s}}$ is bounded using the $L^\infty$-bound on~$\mathscr S$ (see the Proposition~\ref{propS}), the bound on $\txi'$ (by the convexity assumption and the positiveness of~$\tm$ we have : $\txi'(\tm) \geq \txi'(0)$) and the triangular inequality:
\begin{equation*}\label{183}
\|(\nabla \bS)^a\|_{L^q_{\bx,s}}(t) \leq \mathscr S_\infty \int_0^{+\infty} \frac{\mathrm e^{-\We\, \mu T}}{\txi'(0)} \| \nabla \tM\|_{L^q_{\bx,s}}(t,T) \, \mathrm dT.
\end{equation*}
We next use the Cauchy-Schwarz inequality and obtain
\begin{equation*}\label{183.5}
\|(\nabla \bS)^a\|_{L^q_{\bx,s}}(t) \leq \mathscr S_\infty \Big( \int_0^\infty \frac{\mathrm e^{-2\We \mu T}}{\txi'(0)^2} \, \mathrm dT \Big)^{\frac{1}{2}} \Big( \int_0^{+\infty} \mathrm \| \nabla \tM\|_{L^q_{\bx,s}}^2(t,T) \, \mathrm dT \Big)^{\frac{1}{2}}.
\end{equation*}
Taking the $L^r$-norm on $(0,t)$ we deduce 
\begin{equation}\label{184}
\begin{aligned}
\|(\nabla \bS)^a\|_{L^r(0,t;L^q_{\bx,s})}^r
& \leq \mathscr S_\infty^r \Big( \int_0^\infty \frac{\mathrm e^{-2\We \mu T}}{\txi'(0)^2} \, \mathrm dT \Big)^{\frac{r}{2}} y_a(t).
\end{aligned}
\end{equation}
$\checkmark$ Due to the Proposition~\ref{propS} we control the last contribution~$(\nabla \bS)^b$ by
\begin{equation*}\label{185}
(\nabla \bS)^b(t,\bx,s)
\leq
\mathscr S'_\infty \int_0^{+\infty} m(t,T,\bx,s) \, \bigg| \frac{\nabla \bG}{|\bG|} \bigg| (t,T,\bx) \, \mathrm dT.
\end{equation*}
By the triangular inequality we deduce that
\begin{equation*}\label{186}
\|(\nabla \bS)^b\|_{L^q_{\bx,s}}(t)
\leq
\mathscr S'_\infty \int_0^{+\infty} \|m\|_{L^\infty_{\bx,s}}(t,T) \Bigg\| \frac{\nabla \bG}{|\bG|} \Bigg\|_{L^q_{\bx}}(t,T) \, \mathrm dT.
\end{equation*}
From Lemma~\ref{lem:2027} there exists a constant~$C$ such that $\|m\|_{L^\infty_{\bx,s}}(t,T) \leq C \, \mathrm e^{\mu(t - \We T)}$. We deduce
\begin{equation*}\label{187}
\|(\nabla \bS)^b\|_{L^q_{\bx,s}}(t)
\leq
C \int_0^{+\infty} \mathrm e^{\mu(t - \We T)} \Bigg\| \frac{\nabla \bG}{|\bG|} \Bigg\|_{L^q_{\bx}}(t,T) \, \mathrm dT.
\end{equation*}
Taking the $L^r$-norm and using the triangular inequality, we obtain
\begin{equation*}\label{188}
\|(\nabla \bS)^b\|_{L^r(0,t;L^q_{\bx,s})}
\leq
C \int_0^{+\infty} \mathrm e^{-\mu \We T} \Bigg( \int_0^t \mathrm e^{r \mu t'} \Bigg\| \frac{\nabla \bG}{|\bG|} \Bigg\|_{L^q_{\bx}}^r(t',T) \, \mathrm dt' \Bigg)^{1/r} \, \mathrm dT.
\end{equation*}
Writing $\mathrm e^{-\mu \We T} = \mathrm e^{-\mu \We T \frac{r-1}{r}} \mathrm e^{-\mu \We T \frac{1}{r}}$ and using the Hölder inequality, we obtain
\begin{equation*}\label{189}
\begin{aligned}
\|(\nabla \bS)^b\|_{L^r(0,t;L^q_{\bx,s})}^r
& \leq
C \Bigg( \int_0^{+\infty} \mathrm e^{-\mu \We T} \, \mathrm dT \Bigg)^{r-1} \\ & \hspace{1.2cm} \times \int_0^{+\infty} \int_0^t \mathrm e^{r\mu t'} \mathrm e^{-\mu \We T} \Bigg\| \frac{\nabla \bG}{|\bG|} \Bigg\|_{L^q_{\bx}}^r(t',T) \, \mathrm dt' \, \mathrm dT.
\end{aligned}
\end{equation*}
This inequality~\eqref{189} corresponds to
\begin{equation}\label{190}
\|(\nabla \bS)^b\|_{L^r(0,t;L^q_{\bx,s})}^r
\leq C \, y_b(t).
\end{equation}

$\checkmark$ Finally, we add the contributions~\eqref{184} and~\eqref{190} and we use the inequality~\eqref{182} to conclude the proof of the Proposition~\ref{prop:3}.
\cqfd

\subsection{Control of $\boldsymbol{y_a}$}

The goal is now to analyze the quantities $y_a$ and $y_b$. This subsection is devoted to the control of~$y_a$; the next subsection will be devoted to the control of~$y_b$:\\[-0.2cm] 
\begin{lemma}\label{lem:1422}
There is a family of functions~$\txi$, parameterized by $k\geq 1$, for which the function~$y_a$ introduced in the Proposition~\ref{prop:3} satisfies
\begin{equation}\label{191}
y_a' \leq  C \big( 1 + \ln (\mathrm e + y) \, y_a \big) + \frac{C}{k} \big( 1 + \ln (\mathrm e + y)^2 \, y^{\frac{2}{r}} \big) y_a^{1-\frac{2}{r}}.
\end{equation}
The function $y$ is defined by $y=k^{-\frac{r}{2}} y_a + y_b$ and the constant~$C$ does not depend on~$k\geq 1$.
\end{lemma}\\[-0.2cm]
\proof
This proof is decomposed into two main steps. In the first step we establish the evolution equation satisfied by the quantity~$\nabla \tM$ which appear in the definition of~$y_a$. In the second step we show that for a "good" choice of the function~$\txi$ we can deduce an interesting estimate for~$\nabla \tM$, that corresponds to the estimate~\eqref{191} announced in the lemma~\ref{lem:1422}.

\checkmark {\bf Step 1: evolution equation for~$\nabla \tM$} -- 
We recall that~$m$ satisfies
\begin{equation}\label{192}
\mathrm d_t m + \frac{1}{\We} \partial_T m + g \partial_s m - \frac{1}{\We} \partial^2_s m = 0,
\end{equation}
with the following boundary and initial conditions:
\begin{equation*}\label{193}
m\big|_{t=0} = -\partial_T K_0, \qquad
m\big|_{T=0} = 0, \qquad
m\big|_{s=-\frac{1}{2}} = m\big|_{s=\frac{1}{2}} = 0,
\end{equation*}
and where the function~$g$ is given by
\begin{equation*}\label{194}
g(t,\bx,s) = \nabla \bv(t,\bx) : \int_0^s \bS(t,\bx,s')\, \mathrm ds'.
\end{equation*}
To obtain the equation satisfied by $\xi(m)$ we multiply~\eqref{192} by~$\xi'(m)$. We deduce
\begin{equation}\label{194.1}
\mathrm d_t \big( \xi(m) \big) + \frac{1}{\We} \partial_T \big( \xi(m) \big) + g \partial_s \big( \xi(m) \big) - \frac{1}{\We} \xi'(m) \partial^2_s m = 0.
\end{equation}
One of the key points of the proof is that we want the function~$\xi$ provides more regularity. After multiplication by~$\xi'(m)$ the regularizing term~$\partial_s^2 m$ has now becomes
\begin{equation*}\label{194.12}
\xi'(m) \partial^2_s m = \partial^2_s \big( \xi(m) \big) - \frac{\xi''(m)}{(\xi'(m))^2} \big( \partial_s \big( \xi(m) \big) \big)^2.
\end{equation*}
Hence the equation~\eqref{194.1} writes
\begin{equation}\label{194.13}
\mathrm d_t \big( \xi(m) \big) + \frac{1}{\We} \partial_T \big( \xi(m) \big) + g \partial_s \big( \xi(m) \big) - \frac{1}{\We} \partial^2_s \big( \xi(m) \big) + \frac{1}{\We} \frac{\xi''(m)}{(\xi'(m))^2} \big( \partial_s \big( \xi(m) \big) \big)^2 = 0.
\end{equation}
Multiplying this equation~\eqref{194.13} by $\mathrm e^{\We \mu T}$ we obtain the following equation on the quantity $\tM=\mathrm e^{\We \mu T} \xi(m)$:
\begin{equation*}\label{194.3}
\mathrm d_t \tM + \frac{1}{\We} \partial_T \tM - \mu \tM + g \, \partial_s \tM - \frac{1}{\We} \partial^2_s \tM + \frac{\mathrm e^{-\We\, \mu T}}{\We} \frac{\xi''(m)}{(\xi'(m))^2} \big( \partial_s \tM \big)^2 = 0.
\end{equation*}
Finally, denoting $\tm=\mathrm e^{\We \mu T} m$ and $\txi(\tm)=\mathrm e^{\We \mu T} \xi(m)$, we obtain
\begin{equation}\label{194.35}
\mathrm d_t \tM + \frac{1}{\We} \partial_T \tM - \mu \tM + g \, \partial_s \tM - \frac{1}{\We} \partial^2_s \tM + \frac{1}{\We} \frac{\txi''(\tm)}{(\txi'(\tm))^2} \big( \partial_s \tM \big)^2 = 0.
\end{equation}
Taking the spatial gradient of~\eqref{194.35} we obtain
\begin{equation}\label{194.4}
\begin{aligned}
& \mathrm d_t \nabla \tM + \nabla \bv. \nabla \tM + \frac{1}{\We} \partial_T \nabla\tM - \mu \nabla\tM + \nabla g \, \partial_s \tM + g \, \partial_s \nabla \tM - \frac{1}{\We} \partial^2_s \nabla \tM \\
& \hspace{1.5cm} + \frac{1}{\We} \Bigg(\frac{\txi''(\tm)}{(\txi'(\tm))^2}\Bigg)' \frac{\nabla \tM}{\txi'(\tm)} \big( \partial_s \tM \big)^2
+ \frac{2}{\We} \frac{\txi''(\tm)}{(\txi'(\tm))^2} \partial_s \tM \partial_s \nabla \tM = 0.
\end{aligned}
\end{equation}
Remark that the first line of the equation~\eqref{194.4} is exactly the same that those obtained during the local proof (see equation~\eqref{76}). The introduction of the function~$\xi$ makes appear the two last terms in~\eqref{194.4} and at this stage we hope that such terms brings more estimates.

\checkmark {\bf Step 2: choice for the function~$\xi$ in order to estimate~$\nabla \tM$} -- 
We first impose that the function $\txi$ is the solution of the following Cauchy problem:
\begin{equation}\label{194.45}
\left\{
\begin{aligned}
& y"(x) = \xi_0 (y(x)+\xi_2)^k (y'(x))^2, \\
& y(0)=0, \quad y'(0)=\xi_1.
\end{aligned}
\right.
\end{equation}
The parameters $\xi_0$, $\xi_1$, $\xi_2$ and~$k$ are assumed to be positive and will be judiciously chosen later. We simply note that these constants do not depend on~$\bx$ and~$s$ but may depend on the times~$t$ and~$T$. Some results about this ordinary differential equation are given in the last preliminary, see Section~\ref{subsection-ode}.

\noindent
With this choice, the equation~\eqref{194.4} becomes
\begin{equation}\label{194.6}
\begin{aligned}
& \mathrm d_t \nabla \tM + \nabla \bv. \nabla \tM + \frac{1}{\We} \partial_T \nabla\tM - \mu \nabla\tM + \nabla g \, \partial_s \tM + g \, \partial_s \nabla \tM - \frac{1}{\We} \partial^2_s \nabla \tM \\
& \hspace{1cm} + \frac{k \, \xi_0}{\We} ( \txi(\tm) + \xi_2 )^{k-1} \big( \partial_s \tM \big)^2 \nabla \tM
+ \frac{2 \, \xi_0}{\We} ( \txi(\tm) + \xi_2 )^k \partial_s \tM \partial_s \nabla \tM = 0.
\end{aligned}
\end{equation}
We now proceeding as in the proof of the local existence theorem (see page~\pageref{77}). We take inner product of this equation~\eqref{194.6} by $q|\nabla \tM|^{q-2}\nabla \tM$, and integrate with respect to $\bx\in \T$ and $s\in (-\frac{1}{2},\frac{1}{2})$. We deduce
\begin{equation}\label{195}
\begin{aligned}
& \partial_t \|\nabla \tM\|_{L^q_{\bx,s}}^q + \frac{1}{\We} \partial_T \|\nabla \tM\|_{L^q_{\bx,s}}^q - q \mu \|\nabla \tM\|_{L^q_{\bx,s}}^q \\
& + \frac{q}{\We} \int_{\T} \int_{-\frac{1}{2}}^{\frac{1}{2}} |\nabla \tM|^{q-2} |\partial_s \nabla \tM|^2 + \frac{q(q-2)}{\We} \int_{\T} \int_{-\frac{1}{2}}^{\frac{1}{2}} |\nabla \tM \cdot \partial_s \nabla \tM|^2 \, |\nabla \tM|^{q-4} \\
& + \frac{q \, k \, \xi_0}{\We} \int_{\T} \int_{-\frac{1}{2}}^{\frac{1}{2}} ( \txi(\tm) + \xi_2 )^{k-1} |\nabla \tM|^q |\partial_s \tM |^2 
\leq A_0 + A_1 + A_2 + A_3,
\end{aligned}
\end{equation}
where
\begin{equation*}\label{196}
\begin{aligned}
& A_0 = - q \int_{\T} \int_{-\frac{1}{2}}^{\frac{1}{2}} (\nabla \bv \cdot \nabla \tM) \cdot \nabla \tM \, |\nabla \tM|^{q-2}, \\
& A_1 = -q \int_{\T} \int_{-\frac{1}{2}}^{\frac{1}{2}} \Big( \nabla g \cdot \nabla \tM \Big) |\nabla \tM|^{q-2} \partial_s \tM, \\
& A_2 = -q \int_{\T} \int_{-\frac{1}{2}}^{\frac{1}{2}} \Big( \partial_s \nabla \tM \cdot \nabla \tM \Big) |\nabla \tM|^{q-2} g, \\
& A_3 = - \frac{2 q \, \xi_0}{\We} \int_{\T} \int_{-\frac{1}{2}}^{\frac{1}{2}} ( \txi(\tm) + \xi_2 )^k \Big( \partial_s \nabla \tM \cdot \nabla \tM \Big) |\nabla \tM|^{q-2} \partial_s \tM.
\end{aligned}
\end{equation*}
The quantities $A_0$ and $A_2$ are estimate by (see~\eqref{79} and~\eqref{87} respectively for the same kind of estimates):
\begin{equation*}\label{197}
|A_0| \leq q \|\nabla \bv\|_{L^\infty_{\bx}} \| \nabla \tM \|_{L^q_{\bx,s}}^q,
\end{equation*}
\begin{equation*}\label{198}
|A_2| \leq \|\partial_s g\|_{L^\infty_{\bx,s}} \, \|\nabla \tM \|_{L^q_{\bx,s}}^q .
\end{equation*}
The proof that we use to obtain global estimate fundamentally differs from the proof presented for the local existence especially in the control of the term~$A_1$.
We control the contribution~$A_1$ without integration by parts: we simply use the Young inequality and then the Cauchy-Schwarz inequality:
\begin{equation}\label{199}
\begin{aligned}
|A_1| 
& \leq 
\frac{q \, k \, \xi_0}{2\We} \int_{\T} \int_{-\frac{1}{2}}^{\frac{1}{2}} ( \txi(\tm) + \xi_2 )^{k-1} |\nabla \tM|^q |\partial_s \tM |^2 \\
& \qquad + \frac{q \, \We}{2 k \, \xi_0} \left\| \frac{1}{\txi(\tm) + \xi_2} \right\|_{L^\infty_{\bx,s}}^{k-1} \int_{\T} \int_{-\frac{1}{2}}^{\frac{1}{2}} |\nabla g|^2 |\nabla \tM|^{q-2} 
\end{aligned}
\end{equation}
Using the fact that the solution~$\txi$ of the Cauchy problem~\eqref{194.45} is increasing (see Section~\ref{subsection-ode} for an explicit expression of the function~$\txi$), and the fact that~$\tm$ is non-negative (and vanishes), we directly estimate
\begin{equation*}\label{199.2}
\left\| \frac{1}{\txi(\tm) + \xi_2} \right\|_{L^\infty_{\bx,s}}
= \left\| \frac{1}{\txi(0) + \xi_2} \right\|_{L^\infty_{\bx,s}}
= \frac{1}{\xi_2}.
\end{equation*}
Using the Hölder inequality, the estimate~\eqref{199} implies
\begin{equation}\label{199.3}
\begin{aligned}
|A_1| 
& \leq 
\frac{q \, k \, \xi_0}{2\We} \int_{\T} \int_{-\frac{1}{2}}^{\frac{1}{2}} ( \txi(\tm) + \xi_2 )^{k-1} |\nabla \tM|^q |\partial_s \tM |^2 \\
& \qquad +
\frac{q \, \We}{2 k \, \xi_0 \, \xi_2^{k-1}} \|\nabla g\|_{L^q_{\bx,s}}^2 \|\nabla \tM\|_{L^q_{\bx,s}}^{q-2}.
\end{aligned}
\end{equation}
The term $A_3$ is new with respect to the estimate introduced in the local proof theorem. This term comes from to the function~$\xi$ introduced here. We must then control this term with the supplementary contribution given by the function~$\xi$ itself.
In practice, the term~$A_3$ is treated similarly as the term $A_1$, that is using the Young inequality and then the Cauchy-Schwarz inequality:
\begin{equation*}\label{199.4}
\begin{aligned}
|A_3|
& \leq 
\frac{q \, k \, \xi_0}{2\We} \int_{\T} \int_{-\frac{1}{2}}^{\frac{1}{2}} ( \txi(\tm) + \xi_2 )^{k-1} |\nabla \tM|^q |\partial_s \tM |^2 \\
& \qquad +
\frac{2 q \, \xi_0}{k \, \We} \| \txi(\tm) + \xi_2 \|_{L^\infty_{\bx,s}}^{k+1} \int_{\T} \int_{-\frac{1}{2}}^{\frac{1}{2}} |\nabla \tM|^{q-2} |\partial_s \nabla \tM|^2.
\end{aligned}
\end{equation*}
Consequently to control $A_3$ with the left hand side member of~\eqref{195} we want to choose the parameters $\xi_0$, $\xi_1$, $\xi_2$ and~$k$ such that
\begin{equation}\label{199.5}
\frac{2 q \, \xi_0}{k \, \We} \| \txi(\tm) + \xi_2 \|_{L^\infty_{\bx,s}}^{k+1}
\leq 
\frac{q}{\We}.
\end{equation}
By lemma~\ref{lem:2027} we know that for all $(t,T,\bx,s)\in (0,t^\star)\times \R^+ \times \T \times (-\frac{1}{2},\frac{1}{2})$ we have $0\leq \tm(t,T,\bx,s) \leq \tm_\infty = C$.
Since the function~$\txi$ is increasing the condition~\eqref{199.5} also reads
\begin{equation*}\label{199.51}
\txi(\tm_\infty) \leq \Big( \frac{k}{2 \xi_0} \Big)^{\frac{1}{k+1}} - \xi_2.
\end{equation*}
Using the preliminary result given by the Proposition~\ref{prop:4}, we have an explicit expression for the solution~$\txi$, making appear a function~$F$. The condition~\eqref{199.5} is then equivalent to  
\begin{equation}\label{199.52}
\xi_1 \leq \frac{1}{\tm_{\infty}}F \Big( \Big( \frac{k}{2 \xi_0} \Big)^{\frac{1}{k+1}} - \xi_2 \Big) \mathrm e^{\frac{\xi_0 \xi_2^{k+1}}{k+1}}.
\end{equation}
In the sequel, we choose~$\xi_0$ and~$\xi_1$ with respect to the parameters~$\xi_2$ and~$k$ as follows:
\begin{equation}\label{199.6}
\xi_0 = \frac{1}{(k \xi_2)^{k+1}}
\quad \text{and} \quad
\xi_1 = \frac{1}{\tm_{\infty}}F \Big( k \Big( \frac{k}{2} \Big)^{\frac{1}{k+1}} \xi_2 - \xi_2 \Big) \mathrm e^{\frac{1}{(k+1)k^{k+1}}}.
\end{equation}
With this choice, the inequality~\eqref{199.52} holds, hence the inequality~\eqref{199.5} holds too.
\\[-0.2cm]

\begin{remark}\label{rem:1}
The choices of~$\xi_0$ and~$\xi_1$ are fundamental since they ensure the validity of~\eqref{199.5}, but they are also for the following three reasons:
\begin{enumerate}
\item First it is important to notice that with such coefficient the real~$\txi(\tm)$ is defined for any $\tm\in [0,\tm_{\infty}]$. Indeed the set of definition of the function~$\xi$ given in the proposition~\ref{prop:4} writes
\begin{equation*}\label{199.310}
\Big[0, \frac{\ell}{\xi_1} \mathrm e^{\frac{\xi_0 \xi_2^{k+1}}{k+1}}\Big[ 
= 
\Big[0, \frac{\ell}{F \Big( k \Big( \frac{k}{2} \Big)^{\frac{1}{k+1}} \xi_2 - \xi_2 \Big)}\tm_{\infty} \Big[
\subset
[0,\tm_{\infty}],
\end{equation*}
since we have $\displaystyle F \Big( k \Big( \frac{k}{2} \Big)^{\frac{1}{k+1}} \xi_2 - \xi_2 \Big) < \ell = \lim_{X\to +\infty} F(X)$.

\item We also note that with this choice the estimate~\eqref{199.3} for~$A_1$ becomes
\begin{equation}\label{199.31}
\begin{aligned}
|A_1| 
\leq 
\frac{q \, k \, \xi_0}{2\We} \int_{\T} \int_{-\frac{1}{2}}^{\frac{1}{2}} & ( \txi(\tm) + \xi_2 )^{k-1} |\nabla \tM|^q |\partial_s \tM |^2 \\
& \qquad + \frac{q \We \xi_2^2}{2k} \|\nabla g\|_{L^q_{\bx,s}}^2 \|\nabla \tM\|_{L^q_{\bx,s}}^{q-2}.
\end{aligned}
\end{equation}
We will see later that the coefficient~$k$ at the denominator of the last term will make this contribution as small as desired (letting~$k$ to~$+\infty$).
\item Using the expression of the function~$F$ given by~\eqref{241.5}, the value of~$\xi_1$ introduce by~\eqref{199.6} is written
\begin{equation*}\label{199.71}
\xi_1
=
\frac{1}{\tm_{\infty}} \mathrm e^{\frac{1}{(k+1)k^{k+1}}} \int_0^{(k(\frac{k}{2})^\frac{1}{k+1} - 1) \xi_2} \mathrm e^{-\frac{\xi_0}{k+1}(x+\xi_2)^{k+1}} \, \mathrm dx.
\end{equation*}
Performing the following change of variable $z=\frac{x+\xi_2}{\xi_2}$ in the integral, and using the expression~\eqref{199.6} for $\xi_0$ we obtain
\begin{equation*}\label{199.72}
\xi_1
=
\frac{1}{\tm_{\infty}} \mathrm e^{\frac{1}{(k+1)k^{k+1}}} \xi_2 \int_1^{k(\frac{k}{2})^\frac{1}{k+1}} \mathrm e^{-\frac{1}{(k+1)k^{k+1}} z^{k+1}} \, \mathrm dz.
\end{equation*}
We can find a lower bound for the integral using the fact that for $z\leq k(\frac{k}{2})^\frac{1}{k+1}$ we have
$$
\mathrm e^{-\frac{1}{(k+1)k^{k+1}} z^{k+1}} \geq \mathrm e^{-\frac{k}{2(k+1)}} \geq \mathrm e^{-\frac{1}{2}}.
$$
We obtain
\begin{equation*}\label{199.73}
\xi_1
\geq
\frac{1}{\tm_{\infty}} \mathrm e^{\frac{1}{(k+1)k^{k+1}}} \xi_2 \Big( k(\frac{k}{2})^\frac{1}{k+1} -1 \Big) \mathrm e^{-\frac{1}{2}}.
\end{equation*}
Finally, since $\mathrm e^{\frac{1}{(k+1)k^{k+1}}} \geq 1$ and $(\frac{k}{2})^\frac{1}{k+1} = \mathrm e^{\frac{1}{k+1} \ln(\frac{k}{2})} \geq 1$ for $k\geq 2$, we deduce
\begin{equation}\label{199.7}
\xi_1 \geq C\, \xi_2\, k ,
\end{equation}
where the constant~$C$ does not depend on~$\xi_2$ nor on~$k$.
\end{enumerate}
\end{remark}

\vspace{0.2cm}
\noindent
Using \eqref{197}, \eqref{198} and~\eqref{199.31}, the estimate~\eqref{195} now write
\begin{equation}\label{200}
\begin{aligned}
\partial_t \|\nabla \tM\|_{L^q_{\bx,s}}^q + \frac{1}{\We}\partial_T \|\nabla \tM\|_{L^q_{\bx,s}}^q 
\leq
& \big( q \|\nabla \bv\|_{L^\infty_{\bx}} + \|\partial_s g\|_{L^\infty_{\bx,s}} + q \mu \big) \|\nabla \tM \|_{L^q_{\bx,s}}^q \\
& \qquad + \frac{q \, \We \, \xi_2^2}{2k} \|\nabla g\|_{L^q_{\bx,s}}^2 \|\nabla \tM\|_{L^q_{\bx,s}}^{q-2}.\end{aligned}
\end{equation}
We multiply this estimate~\eqref{200} by $\frac{2}{q}\|\nabla \tM\|_{L^q_{\bx,s}}^{2-q}$ and deduce
\begin{equation*}\label{201}
\begin{aligned}
\partial_t \|\nabla \tM\|_{L^q_{\bx,s}}^2 + \frac{1}{\We}\partial_T \|\nabla \tM\|_{L^q_{\bx,s}}^2 
\leq
& \big( 2 \|\nabla \bv\|_{L^\infty_{\bx}} + \frac{2}{q}\|\partial_s g\|_{L^\infty_{\bx,s}} + 2 \mu \big) \|\nabla \tM \|_{L^q_{\bx,s}}^2 \\
& \qquad + \frac{\We \, \xi_2^2}{k} \|\nabla g\|_{L^q_{\bx,s}}^2.
\end{aligned}
\end{equation*}
We integrate with respect to $T\in (0,+\infty)$.
Since $m|_{T=0}=0$ the contribution due to the term $\partial_T \|\nabla \tM\|_{L^q_{\bx,s}}^2$ is non negative.
Finally, we obtain
\begin{equation}\label{202}
\begin{aligned}
\partial_t \bigg( \int_0^{+\infty}  \|\nabla \tM\|_{L^q_{\bx,s}}^2 \, \mathrm dT \bigg) 
\leq
& \big( 2 \|\nabla \bv\|_{L^\infty_{\bx}} + \frac{2}{q} \|\partial_s g\|_{L^\infty_{\bx,s}} + 2 \mu \big) \int_0^{+\infty} \|\nabla \tM \|_{L^q_{\bx,s}}^2 \, \mathrm dT  \\
& \qquad + \bigg( \int_0^{+\infty} \frac{\We \, \xi_2^2}{k}  \, \mathrm dT \bigg)  \|\nabla g\|_{L^q_{\bx,s}}^2.
\end{aligned}
\end{equation}
We finally chose $k$ independent of the time $T$, and
\begin{equation}\label{202.5}
\xi_2 = \frac{1}{1+T},
\end{equation}
so that the estimate~\eqref{202} becomes
\begin{equation*}\label{202.6}
\begin{aligned}
\partial_t \bigg( \int_0^{+\infty}  \|\nabla \tM\|_{L^q_{\bx,s}}^2 \, \mathrm dT \bigg) 
\leq
& \big( 2 \|\nabla \bv\|_{L^\infty_{\bx}} + \frac{2}{q} \|\partial_s g\|_{L^\infty_{\bx,s}} + 2 \mu \big) \int_0^{+\infty} \|\nabla \tM \|_{L^q_{\bx,s}}^2 \, \mathrm dT  \\
& \qquad + \frac{\We}{k} \|\nabla g\|_{L^q_{\bx,s}}^2.
\end{aligned}
\end{equation*}
To make appear an equation of evolution on~$y_a$ we multiply the estimate~\eqref{202} by $\frac{r}{2} \big(  \int_0^{+\infty}  \|\nabla \tM \|_{L^q_{\bx,s}} \, \mathrm dT \big)^{\frac{r}{2}-1}$. We deduce
\begin{equation}\label{203}
\begin{aligned}
y_a''
\leq 
& \big( r \|\nabla \bv\|_{L^\infty_{\bx}} + \frac{r}{q} \|\partial_s g\|_{L^\infty_{\bx,s}} + r \mu \big) y_a' + \frac{r \, \We}{2k} \|\nabla g\|_{L^q_{\bx,s}}^2 \big( y_a' \big)^{1-\frac{2}{r}}.
\end{aligned}
\end{equation}
Integrating the Equation~\eqref{203} with respect to time~$t'\in (0,t)$, we obtain
\begin{equation*}\label{204}
\begin{aligned}
y_a'(t) - y_a'(0)
\leq 
& \big( r \|\nabla \bv\|_{L^\infty(0,t;L^\infty_{\bx})} + \frac{r}{q} \|\partial_s g\|_{L^\infty(0,t;L^\infty_{\bx,s})} + r \mu \big) y_a(t) \\
& \qquad + \frac{r \We}{2k} \int_0^t \|\nabla g\|_{L^q_{\bx,s}}^2\!(t') \, y_a'(t')^{1-\frac{2}{r}} \mathrm dt'.
\end{aligned}
\end{equation*}
The last integral is treated using the Hölder inequality again (we also use the fact that $y_a'\geq 0$). We deduce
\begin{equation}\label{205}
\begin{aligned}
y_a'(t) - y_a'(0)
\leq
& \big( r \|\nabla \bv\|_{L^\infty(0,t;L^\infty_{\bx})} + \frac{r}{q} \|\partial_s g\|_{L^\infty(0,t;L^\infty_{\bx,s})} + r \mu \big) y_a(t) \\
& \qquad + \frac{r \We}{2k} \|\nabla g\|_{L^r(0,t;L^q_{\bx,s})}^2 y_a(t)^{1-\frac{2}{r}}.
\end{aligned}
\end{equation}
Using the definition of $g=\nabla \bv:\int_0^s\bS\, \mathrm ds$ we note that
\begin{equation}\label{206}
\begin{aligned}
& \| \nabla g\|_{L^r(0,t;L^q_{\bx,s})}^2 \leq 2 \| \nabla \bv\|_{L^\infty(0,t;L^\infty_{\bx})}^2  \| \nabla \bS\|_{L^r(0,t;L^q_{\bx,s})}^2 + 2 \| \nabla^2 \bv\|_{L^r(0,t;L^q_{\bx})}^2  \| \bS\|_{L^\infty(0,t;L^\infty_{\bx,s})}^2 , \\
& \| \partial_s g\|_{L^\infty(0,t;L^\infty_{\bx,s})} \leq \| \nabla \bv\|_{L^\infty(0,t;L^\infty_{\bx})} \| \bS\|_{L^\infty(0,t;L^\infty_{\bx,s})}.
\end{aligned}
\end{equation}
Moreover, using successively the proposition~\ref{prop:3}, the remark~\ref{rem:1} (more precisely equation~\eqref{199.7}) and the expression~\eqref{202.5} for~$\xi_2$, we have
\begin{equation*}\label{207}
\begin{aligned}
\| \nabla \bS\|_{L^r(0,t; L^q_{\bx,s})}^2
& \leq C \bigg( \Big( \int_0^\infty \frac{\mathrm e^{-2\We \mu T}}{\xi_1(T)} \, \mathrm dT \Big) \, y_a(t)^{\frac{2}{r}} + y_b(t)^{\frac{2}{r}} \bigg) \\
& \leq C \bigg( \frac{1}{k} \Big( \int_0^\infty \frac{\mathrm e^{-2\We \mu T}}{\xi_2(T)} \, \mathrm dT \Big) \, y_a(t)^{\frac{2}{r}} + y_b(t)^{\frac{2}{r}} \bigg) \\
& \leq C \big( \frac{1}{k} y_a(t)^{\frac{2}{r}} + y_b(t)^{\frac{2}{r}} \big).
\end{aligned}
\end{equation*}
Introducing $y=k^{-\frac{r}{2}} y_a + y_b$ and using the inequality $a^{\frac{2}{r}} + b^{\frac{2}{r}} \leq 2^{1-\frac{2}{r}} (a+b)^{\frac{2}{r}}$, holds for $r\geq 2$, $a>0$ and $b>0$, we obtain
\begin{equation}\label{207.5}
\| \nabla \bS\|_{L^r(0,t; L^q_{\bx,s})}^2 \leq C y(t)^{\frac{2}{r}}.
\end{equation}
We also deduce from the Proposition~\ref{prop:2} (and using~\eqref{207.5}):
\begin{equation}\label{208}
\begin{aligned}
\| \nabla^2 \bv\|_{L^r(0,t; L^q_{\bx})}
& \leq C \big( 1 + \| \nabla \bsigma\|_{L^r(0,t; L^q_{\bx})} \big)\\
& \leq C \big( 1 + \| \nabla \bS\|_{L^r(0,t; L^q_{\bx,s})} \big) \\
& \leq C \big( 1 + y(t)^{\frac{1}{r}} \big).
\end{aligned}
\end{equation}
From the Proposition~\ref{prop:2} again we have\footnote{We also use the fact that there exists two constants $C$ and $C'$ (depending on~$r$) such that for any $y>0$ we have $ \ln(\mathrm e + Cy^{\frac{1}{r}}) \leq C' \ln(\mathrm e + y)$. This is due to the fact that the function $\displaystyle \ell:y\mapsto \frac{\ln (\mathrm e + Cy^{\frac{1}{r}})}{\ln(\mathrm e + y)}$ is bounded on $\R^+$ (in fact it is continuous and satisfies $\ell(0)=1$, $\displaystyle \lim_{+\infty}\ell=\frac{1}{r}$).}
\begin{equation}\label{209}
\begin{aligned}
\| \nabla \bv\|_{L^\infty(0,t; L^\infty_{\bx})} 
& \leq C \ln(\mathrm e + \| \nabla \bsigma\|_{L^r(0,t; L^q_{\bx})} ) \\
& \leq C \ln(\mathrm e + Cy(t)^{\frac{1}{r}}) \\
& \leq C \ln(\mathrm e + y(t)).
\end{aligned}
\end{equation}
We have already seen that the orientation tensor is bounded (see Equation~\eqref{175} obtained to get a $L^\infty$-bound on the stress~$\bsigma$):
\begin{equation}\label{210}
\| \bS\|_{L^\infty(0,t;L^\infty_{\bx,s})} \leq C.
\end{equation}
Finally, the initial value $y_a'(0)$ is estimate using the assumption $\mathrm e^{\mu T/2} \nabla m_0 \in L^2_T L^q_{\bx,s}$:
\begin{equation}\label{211}
\begin{aligned}
y_a'(0) 
& = \Big( \int_0^{+\infty} \| \nabla m_0\|_{L^q_{\bx,s}} \, \mathrm dT \Big)^r \\
& \leq \Big( \int_0^{+\infty} \mathrm e^{-\mu T} \, \mathrm dT \Big)^{\frac{r}{2}} 
\Big( \int_0^{+\infty} e^{\mu T} \| \nabla m_0\|_{L^q_{\bx,s}}^2 \, \mathrm dT \Big)^{\frac{r}{2}} \\
& \leq \mu^{-\frac{r}{2}} \| \mathrm e^{\mu T/2} \nabla m_0 \|_{L^2_T L^q_{\bx,s}}^r \leq C.
\end{aligned}
\end{equation}
Gathering the six previous estimates~\eqref{206}--\eqref{211}, we write the inequality~\eqref{205} as follows
\begin{equation*}\label{212}
y_a' \leq  C \big( 1 + \ln (\mathrm e + y) \, y_a \big) + \frac{C}{k} \big( 1 + \ln (\mathrm e + y)^2 \, y^{\frac{2}{r}} \big) y_a^{1-\frac{2}{r}}.
\end{equation*}
\cqfd

\subsection{Control of $\boldsymbol{y_b}$}

We will now control the quantity~$y_b$ in order to obtain the following result:\\[-0.2cm]
\begin{lemma}\label{lem:1037}
The function~$y_b$ introduced in the Proposition~\ref{prop:3} satisfies
\begin{equation*}\label{213}
y_b' \leq C \big( 1+ \ln (\mathrm e + y) \, y_b \big),
\end{equation*}
where we recall that the function $y$ is defined by $y=k^{-\frac{r}{2}} y_a + y_b$ and the constant~$C$ does not depend on~$k\geq 1$.
\end{lemma}\\[-0.2cm]
\proof
Recall that the function~$y_b$ is given by (see~\eqref{180})
\begin{equation*}\label{213.5}
y_b(t) = \int_0^t \int_0^{+\infty} \mathrm e^{r \mu t'} \mathrm e^{-\We \mu T} \, \Bigg\| \frac{\nabla \bG}{|\bG|} \Bigg\|_{L^q_{\bx}}^r(t',T) \, \mathrm dT \, \mathrm dt'.
\end{equation*}
To prove the lemma~\ref{lem:1037} we must write an equation of evolution for~$\frac{\nabla \bG}{|\bG|}$.
The equation satisfied by~$\bG$ reads
\begin{equation}\label{214}
\mathcal D \bG = \bG \cdot \nabla \bv,
\end{equation}
where we recall that~$\mathcal D$ corresponds to the operator~$\mathcal D=d_t+\partial_T$. We take the inner product of Equation~\eqref{214} by $-q|\nabla \bG|^{q}|\bG|^{-q-2}\bG$:
\begin{equation*}\label{215}
|\nabla \bG|^{q} \mathcal D |\bG|^{-q} = -q|\nabla \bG|^{q}|\bG|^{-q-2} (\bG \cdot \nabla \bv):\bG.
\end{equation*}
Using the Cauchy-Schwarz inequality, we deduce
\begin{equation}\label{216}
|\nabla \bG|^{q} \mathcal D |\bG|^{-q} \leq q|\nabla \bG|^{q}|\bG|^{-q} |\nabla \bv|.
\end{equation}
Next we take the spatial derivative of Equation~\eqref{214}.
We obtain the following $3$-tensor equation
\begin{equation}\label{217}
\mathcal D \nabla \bG = \nabla \bG \cdot \nabla \bv + (\bG\cdot \nabla^2\bv)^\dag - \nabla \bv \cdot \nabla \bG.
\end{equation}
More precisely, the component $(i,j,k)$ of this equation reads 
\begin{equation*}\label{218}
\mathcal D \partial_i \bG_{jk} = \partial_i\bG_{j\ell} \partial_\ell \bv_k + \bG_{j\ell} \partial_\ell \partial_i \bv_k - \partial_i\bv_\ell\partial_\ell \bG_{jk}.
\end{equation*}
Taking the inner product of this equation~\eqref{217} by $q|\bG|^{-q}|\nabla \bG|^{q-2} \nabla \bG$ and using the Cauchy-Schwarz inequality we deduce
\begin{equation}\label{219}
|\bG|^{-q} \mathcal D |\nabla \bG|^q \leq 2q|\nabla \bG|^{q}|\bG|^{-q}|\nabla \bv| + q |\nabla \bG|^{q-1} |\bG|^{-(q-1)} |\nabla^2\bv|.
\end{equation}
Adding this inequality~\eqref{219} with inequality~\eqref{216} we deduce
\begin{equation*}\label{220}
\mathcal D \big(|\nabla \bG|^q|\bG|^{-q}\big) \leq 3q|\nabla \bG|^{q}|\bG|^{-q}|\nabla \bv| + q |\nabla \bG|^{q-1} |\bG|^{-(q-1)} |\nabla^2\bv|.
\end{equation*}
Integrating with respect to the spatial variable we obtain
\begin{equation*}\label{221}
\partial_t \Big\|\frac{\nabla \bG}{|\bG|} \Big\|_{L^q_{\bx}}^q + \partial_T \Big\|\frac{\nabla \bG}{|\bG|} \Big\|_{L^q_{\bx}}^q
\leq
3q \int_{\T} \Big|\frac{\nabla \bG}{|\bG|} \Big|^q |\nabla \bv| + q\int_{\T} \Big|\frac{|\nabla \bG|}{|\bG|}\Big|^{q-1} |\nabla^2\bv|.
\end{equation*}
We now use the Hölder inequality to write
\begin{equation}\label{222}
\begin{aligned}
\partial_t \Big\|\frac{\nabla \bG}{|\bG|} \Big\|_{L^q_{\bx}}^q + \partial_T \Big\|\frac{\nabla \bG}{|\bG|} \Big\|_{L^q_{\bx}}^q
\leq
3 q \Big\|\frac{\nabla \bG}{|\bG|} \Big\|_{L^q_{\bx}}^q \|\nabla \bv\|_{L^\infty_{\bx}} 
+ q \Big\|\frac{|\nabla \bG|}{|\bG|}\Big\|_{L^q_{\bx}}^{q-1} \|\nabla^2\bv\|_{L^q_{\bx}}.
\end{aligned}
\end{equation}
We multiply~\eqref{222} by $\displaystyle \frac{r}{q}\Big\|\frac{\nabla \bG}{|\bG|}\Big\|_{L^q_{\bx}}^{r-q}$ to have
\begin{equation*}\label{223}
\begin{aligned}
\partial_t \Big\|\frac{\nabla \bG}{|\bG|}\Big\|_{L^q_{\bx}}^r + \partial_T \Big\|\frac{\nabla \bG}{|\bG|}\Big\|_{L^q_{\bx}}^r
\leq 3 r \Big\|\frac{\nabla \bG}{|\bG|}\Big\|_{L^q_{\bx}}^r \|\nabla \bv \|_{L^\infty_{\bx}}
+ r \Big\|\frac{\nabla \bG}{|\bG|}\Big\|_{L^q_{\bx}}^{r-1} \|\nabla^2 \bv \|_{L^q_{\bx}}.
\end{aligned}
\end{equation*}
Using the Young inequality we obtain
\begin{equation*}\label{224}
\begin{aligned}
\partial_t \Big\|\frac{\nabla \bG}{|\bG|}\Big\|_{L^q_{\bx}}^r + \partial_T \Big\|\frac{\nabla \bG}{|\bG|}\Big\|_{L^q_{\bx}}^r
& \leq 3 r \Big\|\frac{\nabla \bG}{|\bG|}\Big\|_{L^q_{\bx}}^r \|\nabla \bv \|_{L^\infty_{\bx}} \\
& + \Big\|\frac{\nabla \bG}{|\bG|}\Big\|_{L^q_{\bx}}^r + (r-1)^{r-1} \|\nabla^2 \bv \|_{L^q_{\bx}}^r.
\end{aligned}
\end{equation*}
The quantity $\displaystyle \mathcal G(t,T) = \mathrm e^{\mu r t} \Big\|\frac{\nabla \bG}{|\bG|}\Big\|_{L^q_{\bx}}^r(t,T)$ satisfies
\begin{equation*}\label{225}
\partial_t \mathcal G + \partial_T \mathcal G
\leq 3 r \|\nabla \bv \|_{L^\infty_{\bx}} \mathcal G + (1+\mu r) \mathcal G + (r-1)^{r-1} \mathrm e^{\mu r t} \|\nabla^2 \bv \|_{L^q_{\bx}}^r.
\end{equation*}
We multiply by~$\mathrm e^{-\We \mu T}$ and integrate for $T\in (0,+\infty)$. Using the property $\bG\big|_{T=0}=\bdelta$ we deduce that the second term is non negative (it equals $\mu \, y_b'$). We obtain
\begin{equation*}\label{226}
y_b'' \leq 3r\, \|\nabla \bv \|_{L^\infty_{\bx}} y_b' + (1+ \mu r)y_b' + (r-1)^{r-1} \mathrm e^{\mu r t} \|\nabla^2 \bv \|_{L^q_{\bx}}^r.
\end{equation*}
We can write for any $0<t'<t$ the relation
\begin{equation*}\label{227}
y_b''(t') \leq 3r\, \|\nabla \bv \|_{L^\infty(0,t;L^\infty_{\bx})}  y_b'(t') +  (1+ \mu r) y_b'(t') + (r-1)^{r-1} \mathrm e^{\mu r t} \|\nabla^2 \bv \|_{L^q_{\bx}}^r(t').
\end{equation*}
Integrating now with respect to time $t'\in(0,t)$ we deduce
\begin{equation}\label{228}
y_b'(t) - y_b'(0) \leq C \big( \|\nabla \bv \|_{L^\infty(0,t;L^\infty_{\bx})} y_b(t) +  1  + y_b(t) + \|\nabla^2 \bv \|_{L^r(0,t;L^q_{\bx})}^r \big).
\end{equation}
The value of $y_b'(0)$ is given with respect to the initial condition~$\bG_0$:
\begin{equation*}\label{229}
y_b'(0) = \int_0^\infty \mathrm e^{-\We \mu T} \Big\| \frac{\nabla \bG_0}{|\bG_0|} \Big\|_{L^q_{\bx}}^r(T) \, \mathrm dT.
\end{equation*}
We will note that $y_b'(0)$ is bounded since $\bG_0\in L^\infty(\R^+;W^{1,q}_{\bx})$ and $|\bG_0|\geq \widetilde \gamma$:
\begin{equation*}\label{230}
y_b'(0) \leq \frac{1}{\widetilde \gamma^{r} \, \We \, \mu} \|\bG_0\|_{L^\infty(\R^+;W^{1,q}_{\bx})}^r \leq C.
\end{equation*}
Moreover, using the relation~\eqref{208} and~\eqref{209} again, that is 
\begin{equation*}\label{231}
\begin{aligned}
& \| \nabla^2 \bv\|_{L^r(0,t; L^q_{\bx})} \leq C(1 + y^{1/r}),\\
& \| \nabla \bv\|_{L^\infty(0,t; L^\infty_{\bx})} \leq C \ln(\mathrm e + y),
\end{aligned}
\end{equation*}
the inequality~\eqref{228} becomes (note that $(1+y^{1/r})^r \leq 2^{r-1}(1+y)$)
\begin{equation*}\label{232}
y_b' \leq C \big( 1 + \ln (\mathrm e + y) \, y_b \big),
\end{equation*}
that concludes the proof of the lemma~\ref{lem:1037}.
\cqfd

\subsection{End of the proof for the global result}

Adding the results obtained in Lemmas~\ref{lem:1422} and~\ref{lem:1037}, we deduce that $y=k^{-\frac{r}{2}}y_a+y_b$ satisfies
\begin{equation}\label{233}
y' \leq  C \big( 1 + \ln (\mathrm e + y) \, y \big) + \frac{C}{k^2} \big( 1 + \ln (\mathrm e + y)^2 \, y^{\frac{2}{r}} \big) y^{1-\frac{2}{r}}.
\end{equation}
Using the Young inequality we have $y^{1-\frac{2}{r}} \leq C(1+y) \leq C(1+ \ln (\mathrm e + y)^2 \, y)$. The equation~\eqref{233} implies
\begin{equation}\label{233.1}
y' \leq  C \big( 1 + \ln (\mathrm e + y) \, y \big) + \frac{C}{k^2} \big( 1 + \ln (\mathrm e + y)^2 \, y \big).
\end{equation}
Since $y\geq 0$ we have $1\leq  (\mathrm e+y) \ln (\mathrm e + y) \leq  (\mathrm e+y) \ln (\mathrm e + y)^2$ and the equation~\eqref{233.1} implies
\begin{equation*}\label{233.2}
y' \leq  C \Big( (\mathrm e + y) \ln (\mathrm e + y) + \frac{1}{k^2} (\mathrm e + y) \ln (\mathrm e + y)^2 \Big).
\end{equation*}
Dividing by $(\mathrm e + y) \ln (\mathrm e + y)^2$ we get
\begin{equation}\label{233.6}
w' \leq  C \Big( \frac{1}{k^2}  - w \Big)
\quad \text{where} ~~ w=-\frac{1}{\ln (\mathrm e + y)}.
\end{equation}
This linear equation~\eqref{233.6}, combining with the initial value $w(0)=\frac{-1}{\ln (\mathrm e + y_0)}$, is equivalent to
\begin{equation}\label{233.7}
w(t) \leq   \frac{1}{k^2} - \Big( \frac{1}{\ln (\mathrm e + y_0)} + \frac{1}{k^2} \Big) \mathrm e^{-Ct}.
\end{equation}
We note that the right hand side member in~\eqref{233.7} can vanish at time~$t_0$ given by
\begin{equation*}\label{233.8}
t_0=\frac{1}{C} \ln \Big( 1 + \frac{k^2}{\mathrm e + y_0} \Big).
\end{equation*}
It verifies $\displaystyle \lim_{k\to +\infty} t_0 = +\infty$.\\[0.2cm]
But by construction, we have $w<0$ so that the time of existence~$t^{max}$ of the solution~$y$ satisfies $t^{max}\geq t_0$.
Since the constant~$k$ can be chosen arbitrarily large, this time can be choose as large as possible: we deduce that~$y$ is bounded up to time~$t^\star$ if $k$ is large enough.

\noindent
This bound on~$y$ implies the bounds on~$\bS$ ({\it via}~\eqref{207.5} and~\eqref{210}) and the bounds on~$\nabla \bv$, using~\eqref{208} and~\eqref{209}. That completes the proof.
\cqfd

\section{Conclusion and open problems}\label{conclusion}

In this paper we proved that the famous Doi-Edwards model describing the dynamics of flexible polymers in melts and concentrated solutions is mathematically well-posed. In particular, we show that the solution exists for all time in~$2D$, irrespective of the data (not necessarily small).
Such results naturally bring a lot of new open questions:\\[0.2cm]
$\checkmark$ The first one is about the periodic assumption: the results that are proved in the present paper correspond to the case where the spatial domain is periodic. The reason of this choice is purely mathematics and it is certainly possible to extend some of the results to the bounded domain case using Dirichlet condition for the velocity field (we note that such a result - local in time - have been obtained for polymer model with integral law, see~\cite{Chupin3}).\\[0.2cm]
$\checkmark$ The second point which can be a source of interest is the existence of stationary solution. From another point of view, the theorem~\ref{th:global} indicates that the solution exists on $(0,t^\star)$ for any time~$t^\star>0$. But we don't know if there exists a solution bounded on~$(0,+\infty)$, whose the long time limit would be a stationary solution.\\[0.2cm]
$\checkmark$ Finally, as it was specified in the introduction, much progress has been made on the modeling of both linear and branched polymers since the pioneering works of M. Doi and S.F. Edwards. There exists many models, usually derived from the Doi-Edwards model, which take into account more and more complex phenomena. An example of such a model is the so-called pom-pom model introduced by McLeish and Larson~\cite{Leish98}, and later modified by Blackwell et al.~\cite{Blackwell00}. To the best of my knowledge, no mathematical results exists for such problem and a great challenge would be to show that they are globally well-posed.

\newpage

\bibliographystyle{plain}
\bibliography{biblio-global}

\end{document}